\documentclass[10pt,leqno]{article}
\usepackage{amssymb,amsmath}

\numberwithin{equation}{section}

\date{}
\def\<{\left < }
\def\>{\right >}

\begin{document}

{.}\vskip.4in
\centerline{\Huge Geometry} 
\vskip.1in
\centerline{\Huge of}
\vskip.1in
\centerline{\Huge Slant Submanifolds}

\vskip.5in
\centerline{\large by}
\vskip.5in
\centerline{\Large Bang-yen Chen}
\vskip.3in
\centerline{\large Department of Mathematics}
\vskip.1in
\centerline{\large Michigan State University}
\vskip.1in
\centerline{\large East Lansing, Michigan}
\vskip1.2in
\centerline{\huge Katholieke Universiteit Leuven}
\vskip.1in
\centerline{\huge 1990}

 \vfill\eject

{.}\vskip1.4in
\centerline{\huge Dedicated to} 
\vskip.3in
\centerline{\huge Professor Tadashi Nagano}
\vskip.3in
\centerline{\large on the occasion of his sixtieth birthday}

\vfill\eject
\markboth{B. Y. Chen}{Geometry of Slant Submanifolds}

\setcounter{page}{3}
\vskip4in
\centerline  {\bf PREFACE}
\vskip.3in
The present volume is the written version of the series
of lectures the author delivered at the Catholic
University of Leuven, Belgium during the period of
June-July, 1990. The main purpose of these talks is to
present some of  author's recent work and also his 
joint works with Professor T. Nagano and Professor Y.
Tazawa of Japan, Professor P. F. Leung of Singapore and
Professor J. M. Morvan of France on geometry of slant
submanifolds and its related subjects in a systematical
way. 

The main references of the results presented in this
volume are the following  articles:

 {[C2]} Differential geometry of real
submanifolds in a Kaehler manifold, {\sl Monatsh. f\" ur
Math.,} {\bf 91} (1981), 257-274,

 {[C6]}  Slant
immersions, {\sl Bull. Austral. Math. Soc.,} {\bf 41}
(1990), 135-147,

 {[CLN]} Totally geodesic submanifolds of
symmetric spaces, III, {\sl preprint,} 1980,

 {[CM3]} Cohomologie des sous-vari\'et\'es
$\alpha$-obliques, {\sl C. R. Acad. Sc. Paris,} {\bf 314}
(1992), 931--934.

 {[CT1]} Slant surfaces of 
codimension  two, {\sl Ann.\ Fac.\ Sc.\ Toulouse
Math.,} {\bf 11} (1990), 29--43.

\noindent and

 {[CT2]}  Slant submanifolds in
complex Euclidean spaces, {\sl Tokyo J. Math.} {\bf 14}
(1991), 101--120.

The author would like to take this opportunity to express
his appreciation to the Catholic University of Leuven for
their invitation and support. He would like to
express his gratitute to his colleagues at Leuven,
especially to Professor Leopold Verstraelen, Dr. Franki
Dillen and Dr. Luc Vrancken for their hospitality. He also
like to take this opportunity to express his heartfelt
thanks to Dr. F. Dillen, Professor B. Smyth, Professor L.
 Vanhecke, Professor L. Verstraelen and Dr. L.
Vrancken for
  many valuable discussions during his visit. 

 Finally, the  author would
also like to  thank  his colleagues at
the Michigan State University, especially to Professor D.
E. Blair, Professor G. D. Ludden and Professor W. E. Kuan
for their constant encouragement through the  years.
\vskip.4in \hskip2.8in Bang-yen Chen 

\hskip2.8in Summer 1990

\vfill\eject
\vskip.7in

\centerline  {\bf CONTENTS}
\vskip.5in
\noindent PREFACE \hfill 3

\vskip.15in
\noindent CONTENTS \hfill 5
\vskip.2in
\centerline {CHAPTER  I}
\vskip.05in
\centerline {\bf INTRODUCTION}
\vskip.1in
\noindent \S 1. Introduction \hfill 7 
\vskip.2in
\centerline {CHAPTER II}
\vskip.05in
\centerline {\bf GENEFAL THEORY}
\vskip.1in
\noindent \S 1. Prelinimaries \hfill 13 

\noindent \S 2. Some examples \hfill 18

\noindent \S 3. Some properties of $P$ and $F$\hfill 20

\noindent \S 4. Minimal slant surfaces and totally real surfaces
\hfill 29

\vskip.2in
\centerline {CHAPTER III}
\vskip.05in
\centerline {\bf GAUSS MAPS AND SLANT
IMMERSIONS}
\vskip.1in
\noindent \S 1. Geometry of $G(2,4)$\hfill 33

\noindent \S 2. Complex structures on $E^4$\hfill 35

\noindent \S 3. Slant surfaces and Gauss map\hfill 39

\noindent \S 4. Doubly slant surfaces in ${\bf C}^2$ 
\hfill 44

\noindent \S 5. Slant surfaces in  almost Hermitian manifolds
\hfill 46

\noindent {}
\vskip.08in
\centerline {CHAPTER IV} 
\vskip.05in
\centerline {\bf CLASSIFICATIONS OF SLANT
SURFACES}
\vskip.1in
 \noindent \S 1. Slant surfaces with parallel mean
curvature vector \hfill 50

\noindent \S 2. Spherical slant surfaces\hfill 53

\noindent \S 3. Slant surfaces with $rk (\nu)<2$ \hfill 67

\noindent \S 4. Slant surfaces with codimension one 
\hfill 71
\vskip.2in

\centerline {CHAPTER V}
\vskip.05in
\centerline {\bf TOPOLOGY AND STABILITY OF SLANT
SUBMANIFOLDS}
\vskip.1in
 \noindent \S 1. Non-compactnes of proper slant submanifolds \hfill
77

\noindent \S 2. Topology of slant surfaces \hfill 84

\noindent \S 3. Cohomology of slant submanifolds \hfill 89

\noindent \S 4. Stablity and index form \hfill 98

\noindent \S 5. Stability of totally geodesic
submanifolds \hfill 105

\vskip.2in
\noindent REFERENCES \hfill 114
\vskip.15in
\noindent SUBJECT INDEX \hfill 121

\vfill\eject

\vskip1in
\centerline {CHAPTER I}
\vskip.2in
\centerline {\bf INTRODUCTION}
\vskip.5in
\noindent \S 1. INTRODUCTION.
\vskip .2in
The theory of  submanifolds of an almost Hermitian
 manifold, in particular of a Kaehlerian
manifold, is one of the most interesting topics in
differential geometry. According to the behaviour of the
tangent bundle of a submanifold, with respect to the
action of the almost complex structure $J$ of the ambient
manifold, there are two  well-known classes of
submanifolds, namely, the complex submanifolds and the
totally real submanifolds. In this volume we will present
the geometry of another important class of submanifolds,
called slant submanifolds.

The theory of submanifolds of an almost Hermitian manifold
or of a Kaehlerian manifold began as a separate area of
study in the last century with the investigation of
algebraic curves and algebraic surfaces in classical
algebraic geometry. Included among the principal
investigators are Riemann, Picard, Enriques, Castelnuovo,
Severi, and Segre.

It was E. K\" ahler
([Ka1]), J. A. Schouten, and  D. van Dantzig ([SD1] and
[SD2]) who first tried to study complex manifolds from the
viewpoint of Riemannian geometry in the early 1930's. In
their studies a Hermitian space with the so-called
symmetric unitary connection was introduced. A Hermitian
space with such a connection is now known as a Kaehlerian
manifold.

It was A. Weil [W1] who in 1947 pointed out there exists in
a complex manifold a tensor field $J$ of type (1,1) whose
square is equal to the negative of the identity
transformation of the tangent bundle, that is, $J^{2}=-I.$
In the same year, C. Ehresmann introduced the notion of an
almost complex manifold as an even-dimensional
differentiable manifold which admits such a tensor field
$J$ of type (1,1).

An almost complex manifold (respectively, complex
manifold) is called an almost Hermitian manifold
(respectively, a Hermitian manifold) if it admits a
Riemannian structure which is compatible with the almost
complex structure $J$. The theory of almost Hermitian
manifolds, Hermitian manifolds, and in particular
Kaehlerian manifolds has become a very interesting and
important branch of modern differential geometry (see,
for instances, [Ch1], [Ko1], [Ko2], [KN1], [KN2], [Mi1],
[Wl1].)

 The study of complex submanifolds of a
Kaehlerian manifold from the differential geometrical
points of view (that is, with emphasis on the Riemannian
metric) was initiated by E. Calabi and others in the early
of 1950's (cf. [Ca1], [Ca2]). Since then it has
became an active and fruitful field in modern differential
geometry. Many important results on Kaehlerian
submanifolds have been obtained by many differential
geometers in the last three decades.  Two nice survey
articles concerning this subject were given by K. Ogiue in
[O1] and [O2]. 

In terms of the behaviour of the
tangent bundle $TN$ of the submanifold $N$,  complex
submanifolds $N$ of an almost Hermitian manifold $(M,g,J)$
are characterized by the condition:
$$J(T_{p}N) \subseteq T_{p}N\leqno(1.1)$$
for any point $p \in N$. In other words, $N$  is a
complex submanifold of $(M,g,J)$ if and only if for any
nonzero vector $X$ tangent to $N$ at any point $p \in N$,
the angle between $JX$ and the tangent plane $T_{p}N$ is
equal to zero, identically.

 Besides complex submanifolds, there is another
important class of submanifolds, called {\it totally real
submanifolds.\/} 

 A totally real
submanifold $N$ of an almost Hermitian manifold M
or, in particular, of a Kaehlerian manifold,  is a
submanifold such that the (almost) complex structure $J$
of the ambient manifold $M$ carries each tangent vector of
$N$ into the corresponding normal space of $N$ in $M$,
that is, $$J(T_{p}N)\subseteq T_{p}^{\perp}N\leqno(1.2)$$
for any point $p\in N$. In other words, $N$  is a totally
real submanifold of $(M,g,J)$ if and only if for any
nonzero vector $X$ tangent to $N$ at any point $p \in N$,
the angle between $JX$ and the tangent plane $T_{p}N$ is
equal to ${\pi\over 2}$, identically.

 The study of totally real
submanifolds from the differential geometric points of
view was intiated in the early 1970's (see [CO1] and
[YK1]). Since then many differential geometers have
contributed many interesting results in this subject.

In this volume I shall present the third important class
of submanifolds of an almost Hermitian manifold $(M,g,J)$
(in particular, of a Kaehlerian manifold), called
{\it slant submanifolds.\/}  A slant
submanifold is defined in [C6] as  a submanifold of
$(M,g,J)$ such that, for any nonzero vector $X\in T_{p}N$,
the angle $\theta(X)$ between $JX$ and the tangent space
$T_{p}N$ is a constant (which is independent of the choice
of the point $p\in N$ and the choice of the tangent vector
$X$ in the tangent plane $T_{p}N$). It is obvious that
complex submanifolds and totally real submanifolds are
special classes of slant submanifolds. A slant submanifold
is called {\it proper\/} if it is neither a complex
submanifold nor a totally real submanifold.

In the first section of Chapter 2, we present some basic
definitions and  basic formulas for later use. In
particular, for any submanifold $N$ of an almost Hermitian
manifold $(M,g,J)$, the almost complex structure $J$ of the
ambient manifold induces a canonical endomorphism of the
tangent bundle, denoted by  $P$, and a canonical
normal-bundle-valued 1-form on the tangent bundle,
denoted by $F$. It will be shown in later sections that the
endomorphism $P$ plays a fundamental role in our study. 

In Section 2 of Chapter 2,  we give many examples of slant
surfaces in the complex number space ${\bf C}^2$ and give
examples of Kaehlerian slant submanifolds (that is,
proper slant submanifolds such that  the canonical
endomorphism $P$ is parallel) in ${\bf C}^m$.

Section 3 of Chapter 2 is devoted to the fundamental study
of the  endomorphism $P$ and the normal-bundle-valued
1-form $F$. In particular, we prove that every slant
surface in any almost Hermitian manifold is a Kaehlerian
slant submanifold, that is, $\nabla P=0$ (Theorem
3.4). Morevoer, we prove that a proper slant submanifold
$N$ of a Kaehlerian manifold $M$ is Kaehlerian slant if and
only if the Weingarten map of $N$ in $M$ satisfies
$$A_{FX}Y=A_{FY}X\leqno(1.3)$$
for any vectors $X,Y$ tangent to the submanifold. In
particular, by combining this result with Theorem 3.4, we
obtain the important fact that formula (1.3) holds for any
slant surface in any Kaehler manifold. By using this
result we show that the Gauss curvature $G$ and the
normal curvature $G^D$ of every slant surface in  ${\bf
C}^2$ satisfies $$G=G^{D},\leqno(1.4)$$ identically. In
this section we also prove that every proper slant
submanifold of a Kaehlerian manifold with $\nabla F=0$ is
an austere submanifold (Theorem 3.8). We also obtain in
this section two reduction theorems for submanifolds
satisying $\nabla F=0$ (Theorem 3.9 and Theorem 3.10).
Finally we show that for a surface in a real 4-dimensional
Kaehlerian manifold, the parallelism of $F$ implies the
parallelism of $P$ (Theorem 3.11).

In Section 4 of Chapter II, we establish some relations
between minimal slant surfaces and totally real surfaces
in a  Kaehlerian manifold. For instance, we prove
that if a proper slant surface in a real 4-dimensional
Kaehlerian manifold is also totally real with respect to
some compatible complex structure on the ambient manifold
at the same time, then the surface must be a minimal
surface (Theorem 4.2). Some applications of this result
will also be given in this section (Theorem 4.3 and
Theorem 4.4).

In Chapter III, we present some results of the author
and Y. Tazawa [CT1]; in this work slant surfaces in ${\bf
C}^2$ were studied from the viewpoint of the Gauss map.

The first section of Chapter III reviews some basic
geometry of the real Grassmannian $G(2,4)$ which consists
of all oriented 2-planes in $E^4$ which will be used in
later sections. 

The second section is devoted to the detailed study of the
set of compatible complex structures on $E^4$. Several
lemmas are obtained in this section, in particular, Lemma
2.1 and Lemma 2.3  of this section, play some important
roles in the later sections.

In Section 3 we study the following two geometric
problems:
\vskip.1in
{\bf Problem 3.1.} {\it Let $N$ be a surface in ${\bf
C}^{2} = (E^{4},J_{0})$. When is $N$ slant in ${\bf
C}^2$?} \vskip.1in 
{\bf Problem 3.2.} {\it Let $N$ be a surface in $E^4$. If
there exists a compatible complex structure $J$ on
$E^4$ such that $N$ is slant in $(E^{4},J)$. How many other
compatible complex structures ${\tilde J}$ on $E^4$ are
there such that $N$ is slant with resepct to these complex
strucutres?} \vskip.1in
Complete solutions to these two problems are
obtained in this section.

Related to the problems studied in Section 3 
is the notion of doubly slant surfaces. 
We prove that every doubly slant surface in $E^4$ has
vanishing Gauss curvature and vanishing normal curvature
(Theorem 4.1).

In the last section of Chapter III, we study slant
surfaces in almost Hermitian manifolds. In fact we prove
that if a surface in a real 4-dimensional almost
Hermitian manifold has no complex tangent points (that
is, the surface is purely real), then, with respect to
some suitable compatible almost complex structure on the
ambient manifold, the surface
is slant. This result shows that there exist ample
examples of slant surfaces in almost Hermitian manifolds.

Chapter IV is devoted to the classification problems of
slant surfaces in ${\bf C}^{2}$.
In the first section of this chapter, we classify all
slant surfaces in ${\bf
C}^{2}$ with parallel mean curvature vector.
In the second section we classify spherical slant
surfaces in ${\bf C}^{2}$.
 In the third section, we classify slant surfaces in ${\bf
C}^{2}$ whose Gauss map has rank $< 2$ at every point. The
last section gives the classification of slant surfaces
in ${\bf C}^2$ which are contained in a hyperplane of ${\bf
C}^{2}$.

In the last Chapter, we study the topology and cohomology
of slant submanifolds. 

In the first section, we prove
that every compact slant submanifold in ${\bf C}^{m}$
is totally real (Theorem 1.5) and every compact slant
submanifold in an exact sympletic manifold is also
totally real (Theorem 1.7). 

In Section 2, we define a
canonical 1-form $\Theta$ associated with a proper
slant submanifold $N$ in an almost Hermitian manifold (by
formula (2.7) of Chapter V). We prove in this section
that for a slant surface in ${\bf C}^{2}$, this 1-form is
closed, that is, $d\Theta=0$ and if we put
$$\Psi=(2\sqrt{2}\pi)^{-1}(\csc\alpha)\Theta,\leqno
(1.5)$$

\noindent where $\alpha$ denotes the slant angle of $N$ in
${\bf C}^2$, then the 1-form $\Psi$ defines a
 canonical integral class on $N$ (Theorem 2.5): 
$$[\Psi] \in H^{1}(N;{\bf
Z}).\leqno(1.6)$$ 

\noindent In this section we also
prove that if $N$ is a complete, oriented, proper slant
surface of ${\bf C}^{2}$ such that the mean curvature of
$N$ is bounded below by some positive constant, then,
topologically, the surface is either a circular cylinder or
a 2-plane (Theorem 2.6). 

In Section 3 of this chapter, we prove that in fact, for
any $n$-dimensional proper slant submanifold of ${\bf
C}^{n}$, the 1-form $\Theta$ is always closed. Thus for any
$n$-dimensional proper slant submanifold $N$ in 
${\bf C}^{n}$, we have a canonical cohomology class
$[\Theta] \in H^{1}(N;{\bf R})$ (Theorem 3.1). Finally we
prove in this section that every proper slant submanifold
of a Kaehlerian manifold has a canonical sympletic
structure given by the 2-form induced from the canonical
endomorphism $P$ (Theorem 3.4). The last result implies
that if a compact $2k$-dimensional differentiable
manifold $N$ satisfies $H^{2i}(N;{\bf R})=0$ for some $i
\in \{1,\ldots,k\}$, then $N$ cannot be immersed in any
Kaehleian manifold as a proper slant submanifold (Theorem
3.$5'$).

In Section 4, we recall 
some stablity theorems of [CLN] obtained in 1980 and
present some related results concerning 
the index form of
 compact minimal totally real submanifolds, a special
class of slant submanifolds, of  a Kaehlerian manifold.

In the last section, we present a general method
introduced by the author, Leung and Nagano [CLN] for
determining the stability of totally geodesic
submanifolds in compact symmetric spaces. Since every
irreducible totally geodesic submanifolds of a Hermitian
symmetric space is a slant submanifold [CN1], the method
can be used to determine the stability of such
submanifolds.

 Finally, I would like to  mention that   minimal
surfaces of a complex projective space with constant
Kaehlerian angle  were recently studied from a  different
point of view by  Bolton, Jensen,  Rigoli, Woodward,
Maeda, Ohnita and Udagawa  (see, [BJRW], [MU1]  and [Oh1]).
\vfill\eject

\centerline {CHAPTER II}
\vskip.2in
\centerline {\bf GENERAL THEORY}
\vskip .5in
\noindent {\S 1}. PRELIMINARIES.
\vskip.2in
Let $N$ be an $n$-dimensional Riemannian manifold
isometrically immersed in an almost Hermitian manifold
$M$ with (almost) complex structure $J$ and almost
Hermitian metric $g$. We denote by $< , >$ the inner
product for $N$ as well as for $M$.

For any vector $X$ tangent to $N$ we put
$$JX=PX+FX,\leqno(1.1)$$
where $PX$ and $ FX$ are the tangential and the normal
components of $JX$, respectively.  Thus, $P$ is an
endomorphim of the tangent bundle $TN$ and $F$ a
normal-bundle-valued 1-form on $TN$.

For any nonzero vector $X$ tangent to $N$ at a point $x
\in N$, the angle $\theta (X)$  between $JX$ and the
tangent space $T_{x}N$ is called the {\it Wirtinger
angle\/} of $X$. In the following  we call an immersion
$f:N \rightarrow M$ a {\it slant immersion\/} if the
Wirtinger angle $\theta (X)$ is constant (which is
independent of the choice of $x \in N$ and of $X \in
T_{x}N$). Complex and totally real immersions are slant
immersions with $\theta = 0$ and $\theta = \pi /2$,
respectively. Moreover, it is easy to see that slant
submanifolds of an almost Hermitian manifolds are
characterized by the condition: $\, P^{2}=\lambda I,\,$
for some real number $\lambda \in [-1,0],$ where $I$
denotes the identity transformation of the tangent bundle
$TN$ of the submanifold $N$. 

 The Wirtinger angle of a
slant immersion is called the {\it slant angle \/} of the
slant immersion. A slant submanifold is said to be {\it
proper\/} if it is neither complex nor totally real.

A proper slant submanifold is said to be {\it Kaehlerian
slant \/} if the canonical endomorphism $P$ defined
above is parallel, that is, $\nabla P = 0$. A Kaehlerian
slant submanifold of a Kaehlerian manifold with respect to
the induced metric and with the almost complex structure
given by ${\tilde J}= (\sec \theta)P, \theta = $ the slant
angle.

For a submanifold $N$ of an almost Hermitian manifold $M$,
we denote by $\nabla$ and ${\tilde \nabla}$ the Levi-Civita
connections of $N$ and $M$, respectively. Then the Gauss
and Weingarten formulas of $N$ in $M$ are given
respectively by $${\tilde \nabla}_{X}Y=\nabla_{X} +
h(X,Y),\leqno(1.2)$$ $${\tilde\nabla}_{X}\xi =
-A_{\xi}X+D_{X}\xi, \leqno(1.3)$$ for any vector fields
$X,Y$ tangent to $N$ and any  vector field $\xi$ normal
to $N$, where $h$ denotes the second fundamental form, $D$
the normal connection, and $A$ the Weingarten map of the
submanifold $N$ in $M$. The second fundamental form $h$
and the Weingarten map $A$ are related by
$$<A_{\xi}X,Y>=<h(X,Y),\xi >.\leqno(1.4)$$

For any vector field $\xi$ normal to the submanifold $N$,
we put $$J\xi = t\xi + f\xi,\leqno(1.5)$$
where $t\xi$ and $f\xi$ are the tangential and the normal
components of $J\xi$, respectively. Then $f$ is an
endomorphism of the normal bundle and $t$ is a
tangent-bundle-valued 1-form on the normal bundle
$T^{\bot}N$.

For a submanifold $N$ in a Riemannian manifold $M$, the
mean curvature vector $H$ is defined by
$$H = {1\over n}\, tr\,h = {1\over n}\sum_{i=1}^{n}
h(e_{i},e_{i})\leqno(1.6)$$
where $\{e_{1},\ldots,e_{n}\}$ is a local  orthonormal
frame of the tangent bundle $TN$ of $N$. A submanifold
$N$ in $M$ is said to be {\it totally geodesic\/} if the
second fundamental form $h$ of $N$ in $M$ vanishes
identically. 

Let $N$ be an $n$-dimensional submanifold  in an
$m$-diemsional Riemannian manifold $M$. We choose a local
field of orthonoraml  frames

$$e_{1},\ldots,e_{n},e_{n+1},\ldots,e_{m}$$ such that,
restricted to $N$, the vectors $e_{1},\ldots,e_{n}$ are
tangent to $N$ and hence $e_{n+1},\ldots,e_{m}$ are normal
to $N$. 
We shall make use of the following convention on the
ranges of indices unless mentioned otherwise:
$$1\leq A,B,C,\ldots \leq m;\,\,\,\,\,\, 1 \leq
i,j,k,\ldots \leq n;$$
$$n+1 \leq  r,s,t,\ldots \leq  m. $$

With respect to the frame field of $M$ chosen above, let
$$\omega^{1},\ldots,\omega^{n},\omega^{n+1},\ldots,\omega^{m}$$
be the field of dual frames. Then the structure equations
of $M$ are given by
$$d\omega^{A}=-\sum_{B} \omega_{B}^{A}\wedge\omega^{B},
\hskip.5in \omega_{B}^{A}+\omega_{A}^{B}=0,\leqno(1.7)$$
and
 $$d\omega_{B}^{A}=-\sum_{C} \omega_{C}^{A}\wedge
\omega_{B}^{C} + \Phi_{B}^{A},\hskip.2in
\Phi_{B}^{A}={1\over 2}\sum_{C,D}
K^{A}_{BCD}\omega^{C}\wedge\omega^{D},\leqno(1.8)$$
$$K^{A}_{BCD}+K^{A}_{BDC}=0.$$

We restrict these forms to $N$. Then
$$\omega^{r}=0.$$
Since
$$0=d\omega^{r}=-\sum_{i}\omega_{i}^{r}\wedge\omega^{i},$$
by Cartan's lemma we have
$$\omega^{r}_{i} = \sum_{j} h_{ij}^{r}
\omega^{j}, \,\,\,\,\,  h^{r}_{ij}=
h^{r}_{ji}.\leqno(1.9)$$ 
Formula (1.9) is equivalent to
$$\omega_{i}^{r}(X)=<A_{e_{r}}e_{i},X>\leqno(1.9)'$$
for any vector $X$ tangent to $N$.
From these formulas we obtain
$$d\omega^{i}=
-\sum_{j} \omega_{j}^{i}\wedge\omega^{j},\,\,\,\,
\omega_{j}^{i}+\omega_{i}^{j}=0,\leqno(1.10)$$
$$d\omega_{j}^{i}=-\sum_{k}\omega_{k}^{i}\wedge\omega_{j}^{k}
+ \Omega_{j}^{i},\,\,\, \Omega_{j}^{i}={1\over
2}\sum_{k,\ell}R^{i}_{jk\ell}\omega^{k}\wedge\omega^{\ell},
\leqno(1.11)$$
$$R^{i}_{jk\ell}=K^{i}_{jk\ell}+\sum_{r}(h^{r}_{ik}h^{r}_{j\ell}
-h^{r}_{i\ell}h^{r}_{jk}),\leqno(1.12)$$
$$d\omega_{i}^{r}=-\sum_{j,k}
h^{r}_{jk}\omega_{i}^{j}\wedge\omega^{k} -\sum_{j,s}
h^{s}_{ij}\omega^{j}\wedge\omega^{r}_{s} +{1\over
2}\sum_{j,k}K^{r}_{ijk}\omega^{j}\wedge\omega^{k},
\leqno(1.13)$$
$$d\omega_{s}^{r}=-\sum_{t}\omega_{t}^{r}\wedge\omega_{s}^{t}
+\Omega_{s}^{r}, \,\,\,\, \Omega_{s}^{r}= {1\over
2}\sum_{k,\ell}R^{r}_{sk\ell}\omega^{k}\wedge
\omega^{\ell}, \leqno(1.14)$$
and
$$R^{r}_{sk\ell}=K^{r}_{sk\ell}+\sum_{i}
(h^{r}_{ik}h^{s}_{i\ell}-h^{r}_{i\ell}h^{s}_{ik}).
\leqno(1.15)$$

For 
any vector field $X$ tangent to the submanifold $N$,
these forms are also given by $${\tilde
\nabla}_{X}e_{i}=\sum_{j=1}^{n}\omega_{i}^{j}(X) e_{j}+
\sum_{r=n+1}^{m}\omega_{i}^{r}(X)e_{r},\leqno(1.16)$$

$${\tilde
\nabla}_{X}e_{r}=\sum_{j=1}^{n}\omega_{r}^{j}(X) e_{j}+
\sum_{s=n+1}^{m}\omega_{r}^{s}(X)e_{s}.\leqno(1.17)$$

 These 1-forms $\omega_{i}^{j}$, $\omega_{i}^{r}$ and
$\omega_{r}^s$ are called the {\it connection forms\/} of
$N$ in $M$.

Denote by $R$ and ${\tilde R}$ the Riemann curvature
tensors of $N$ and $M$, respectively, and by $R^D$ the
curvature tensor of the normal connection $D$. Then the
{\it equation of Gauss\/} and the {\it equation of
Ricci\/} are given respectively by $${\tilde
R}(X,Y;Z,W)=R(X,Y;Z,W)+<h(X,Z),h(Y,W)>\leqno(1.18)$$
$$-<h(X,W),h(Y,Z)>,$$ 
$$R^{D}(X,Y;\xi,\eta)={\tilde R}(X,Y;\xi,\eta)+<[A_{\xi},
A_{\eta}](X),Y>\leqno(1.19)$$ for vectors $X,Y,Z,W$ tangent
to $N$ and $\xi,\eta$ normal to $N$.

For the second fundamental form $h$, we define the
covariant derivative ${\bar \nabla}h$ of $h$ with respect
to the connection in $TN \oplus T^{\perp}N$ by
$$({\bar\nabla}_{X}h)(Y,Z)=D_{X}(h(Y,Z))-h(\nabla_{X}Y,Z)
-h(Y,\nabla_{X}Z).\leqno(1.20)$$

The {\it equation of Codazzi\/} is given by
$$({\tilde R}(X,Y)Z)^{\perp}=({\bar\nabla}_{X}h)(Y,Z)-
({\bar\nabla}_{Y}h)(X,Z),\leqno(1.21)$$
where $({\tilde R}(X,Y)Z)^{\perp}$ denotes the normal
component of ${\tilde R}(X,Y)Z$.

A submanifold $N$ of a Riemannian manifold $M$ is called
a {\it parallel submanifold\/} if the second fundamental
form $h$ is parallel, that is, ${\bar\nabla}h=0$,
identically.

 \vfill\eject

\vskip.3in
\noindent \S 2. SOME EXAMPLES.
\vskip.2in
In the following, $E^{2m}$ denotes the Euclidean
$2m$-space with the standard metric. An almost complex
structure $J$ on $E^{2m}$ is said to be compactible if
$(E^{2m},J)$ is complex analytically isometric to the
complex number space ${\bf C}^m$ with the standard flat
Kaehlerian metric. We denote by $J_0$ and $J_{1}^{-}$
(when $m$ is even) the compatible almost complex structures
on $E^{2m}$ defined respectively by
$$J_{0}(a_{1},\ldots,a_{m},b_{1},\ldots,b_{m})=(-b_{1},
\ldots,-b_{m},a_{1},\ldots,a_{m})\leqno(2.1)$$
and 
$$J_{1}^{-}(a_{1},\ldots,a_{m},b_{1},\ldots,b_{m})\leqno(2.2)$$
$$=(-a_{2},a_{1},\ldots,-a_{m},a_{m-1},b_{2},-b_{1},
\ldots,b_{m},-b_{m-1}).$$

 In this
section we give some examples of proper slant surfaces in
${\bf C}^{2}=(E^{4},J_{0})$ and some examples of
Kaehlerian slant submanifolds of higher dimension in
${\bf C}^{m}$.  
\vskip.1in 
{\bf Example 2.1.} For
any $\alpha >0$, $$x(u,v) = (u\cos \alpha,v, u\sin \alpha,
v,0)$$ defines a slant plane with slant angle $\alpha$ in
${\bf C}^2$.

\vskip.1in
{\bf Example 2.2.} Let $N$ be a complex surface in 
${\bf C}^{2}=(E^{4},J_{0})$. Then for any constant
$\alpha$, $0<\alpha\leq \pi /2, N$ is slant surface in
$(E^{4},J_{\alpha})$ with slant angle $\alpha$, where
$J_{\alpha}$ is the compatible almost complex structure on
$E^{4}$ defined by 
$$J_{\alpha}(a,b,c,d)
=(\cos\alpha)(-c,-d,a,b)+(\sin\alpha)(-b,a,d,-c).$$
This example shows that there exist infinitely many
proper slant minimal surfaces in ${\bf C}^{2}=
(E^{4},J_{0})$.
\vskip.1in
The following example provides us  some non-minimal proper
slant surfaces in ${\bf C}^{2}=(E^{4},J_{0})$. \vskip.1in
{\bf Example 2.3.} [GVV] For any positive constant $k$,
$$x(u,v) =(e^{ku}\cos u\cos v, e^{ku}\sin u\cos
v,e^{ku}\cos u\sin v, e^{ku}\sin u\sin v)$$ defines a
complete, non-minimal, pseudo-umbilical proper slant
surface with slant angle $\theta=\cos
^{-1}(k/\sqrt{1+k^{2}})$ and with non-constant mean
curvature given by $|H|=e^{-ku}/\sqrt{1+k^{2}}$. \vskip.1in
{\bf Example 2.4.} For any positive number $k$,
$$x(u,v)=(u,k\cos v,v,k\sin v)$$
defines a complete, flat, non-minimal and
non-pseudo-umbilical, proper slant surface with slant
angle $\cos^{-1}(1/\sqrt{1+k^{2}})$ and constant mean
curvature $k/2(1+k^{2})$ and with non-parallel mean
curvature vector.
\vskip.1in
{\bf Example 2.5.} Let $k$ be any positive number and
$(g(s),h(s))$  a unit speed plane curve. Then
$$x(u,s) =(-ks\sin u, g(s),ks\cos u, h(s))$$ defines a
non-minimal, flat, proper slant surface with slant angle
$k/\sqrt{1+k^{2}}$.
\vskip.1in
{\bf Example 2.6} [CT1].  For any nonzero real numbers
$p$ and $q$, we consider the following immersion form
${\bf R}\times (0,\infty )$ into ${\bf C}^2$ defined by
$$x(u,v)=(pv\sin u,pv\cos u,v\sin qu, v\cos qu).$$
Then the immersion $x$ gives us a complete flat slant
surface in ${\bf C}^2$.
\vskip.1in
{\bf Example 2.7.} For any $k>0$, $$x(u,v,w,z)=(u,v,k\sin
w,k\sin z, kw,kz,k\cos w,k\cos z)$$ defines a Kaehlerian
slant submanifold in ${\bf C}^4$ with  slant angle
$\cos^{-1}k.$
\vskip.1in
{\bf Example 2.8.} Let $N$ be a complex submanifold of the
complex number space ${\bf C}^{2m}=(E^{4m},J_{0})$. For
any constant $\alpha$ we define $J_{\alpha}$ by
$$J_{\alpha}=(\cos\alpha)J_{0}+(\sin\alpha)J_{1}^{-}.\leqno(2.3)$$
Then $J_{\alpha}$ is a compatible complex structure on
$E^{4m}$ and $N$ is a Kaehlerian slant submanifold with
slant angle $\alpha$ in $(E^{4m},J_{\alpha})$.

\vfill\eject

\vskip.3in

\noindent \S 3. SOME PROPERTIES OF $P$ AND $F$.
\vskip.2in
Let $f: N \rightarrow M$ be an isometric immersion of an
$n$-dimensional Riemannian manifold into an almost
Hermitian manifold. Let $P$ and $F$ be the endomorphism
and the normal-bundle-valued 1-form on the tangent bundle
defined by (1.1). Since $M$ is almost Hermitian, we have
$$<PX,Y>=-<X,PY>$$ for any vectors X, Y tangent to $N$.
Hence, if we put $Q=P^2$, then $Q$ is a self-adjoint
endomorphim of $TN$. Therefore, each tangent space $T_{x}N$
of $N$ at $x \in N$ admits an orthogonal direct
decompostion of eigenspces of $Q$:
$$T_{x}N = {\mathcal D}_{x}^{1} \oplus\dots\oplus {\mathcal
D}_{x}^{k(x)}.$$
Since $P$ is skew-symmetric and $J^{2}=-I, $ each
eigenvalue $\lambda_i$ of $Q$ lies in $[-1,0]$ and,
moreover, if $\lambda_{i} \not= 0,$ then the
corresponding eigenspace ${\mathcal D}_{x}^{i}$ is of even
dimension and it is invariant under the endomorphism $P$,
that is $P({\mathcal D}_{x}^{i})={\mathcal D}_{x}^{i}.$
Furthermore, for each  $\lambda_{i} \not= -1,$ dim
$F({\mathcal D}_{x}^{i})=$ dim ${\mathcal D}_{x}^{i}$ and the
normal subspaces $F({\mathcal D}_{x}^{i})$, 
$i=1,\ldots,k(x),$ are mutually perpendicular. From these
we have $${\rm dim}\, M \geq  2 {\rm dim}\, N - {\rm dim
}\,{\mathcal H}_{x}$$ where ${\mathcal H}_{x}$ denotes the
eigenspace of $Q$ with eigenvalue $-1$.

In this section we mention some results given in [C2] and 
[C6] concerning the endomorphism $P$ and the
normal-bundle-valued 1-form $F$ associated with the
immersion $f : N \rightarrow M$. 

The following Lemma 3.1 follows from the definition of
$\nabla Q$ which is defined by
$$(\nabla_{X}Q)Y=\nabla_{X}(QY)-Q(\nabla_{X}Y)\leqno(3.1)$$
for $X$ and $Y$ tangent to $N$.

\vskip.1in
{\bf Lemma 3.1.} {\it Let $N$ be a submanifold of an
almost Hermitian manifold $M$. Then the self-adjoint
endomorphism $Q \,(= P^{2})$ is parallel, that is, $\nabla
Q=0,$ if and only if}

(1) {\it each eigenvalue $\lambda_i$ of $Q$ is constant on
$N$;}

(2) {\it each distribution ${\mathcal D}^i$ (associated with
the eigenvalue $\lambda_i$) is completely integrable and}

(3) {\it $N$ is locally the Riemannian product
$N_{1}\times\ldots\times N_k$ of the leaves of the
distributions.} \vskip.1in
{\bf Proof.} Since $Q$ is a self-adjoint endomorphism of
the tangent bundle $TN$, there exist $n$ continuous
functions $\lambda_{1}\leq \lambda_{2} \leq \ldots \leq
\lambda_{n}$ such that $\lambda_{i},\,\,i=1,\ldots,n,\,\,$
are the eigenvalues of $Q$ at each point points $p\in N$.
Let $e_{1},\ldots,e_n$ be a local orthonormal frame
given by eigenvectors of $Q$. If $Q$ is parallel, then
(3.1) implies
$$\nabla_{X}(\lambda_{i}e_{i})=Q(\nabla_{X}e_{i}),\,\,\,i=1,\ldots,
n$$ for any vector $X$ tangent to $N$. Thus we have
$$(X\lambda_{i})e_{i}+\lambda_{i}(\nabla_{X}e_{i})=
Q(\nabla_{X}e_{i}).$$
Since both $\nabla_{X}e_i$ and $Q(\nabla_{X}e_{i})$ are
perpendicular to $e_i$, we conclude that each eigenvalue
of $Q$ is constant on $N$. This proves statement (1).

For statements (2) and (3) we let
$\lambda_{1},\ldots,\lambda_k$ denote the distinct
eigenvalues of $Q$. For each $i;\,i=1,\ldots,k\,$, let
${\mathcal D}^i$ denote the distribution given by the
eigenspaces of $Q$ with eigenvalue $\lambda_i$. For any two
vector fields $X,Y$ in the distribution ${\mathcal D}^i$,
(3.1) and statement (1) imply
$$Q(\nabla_{X}Y)=\lambda_{i}(\nabla_{X}Y),$$ from which we
conclude that $\nabla_{X}Y \in {\mathcal D}^i$ for any $X,Y$
in ${\mathcal D}^i$. Therefore, each distribution ${\mathcal D}^i$
is completely integrable and each maximal integrable
submanifold of ${\mathcal D}^i$ is totally geodesic in $N$.
Consequently, $N$ is locally the Riemannian product
$N_{1}\times\ldots\times N_k$ of the leaves of these
distributions. 

The converse of this is easy to verify.

 \vskip.1in By using Lemma 3.1 we
have the following characterization of subma-nifolds with
$\nabla P=0$. \vskip.1in 
\eject
{\bf Lemma 3.2.} {\it Let $N$ be a
submanifold of an almost Hermitian manifold $M$. Then
$\nabla P=0$ if and only if $N$ is locally the Riemannian
product $N_{1}\times\ldots \times N_{k}$, where each $N_i$
is either a complex submanifold, a totally real
submanifold, or a Kaehlerian slant submanifold of $M$.}
\vskip.1in
{\bf Proof.} Under the hypothesis, if $P$ is parallel,
then $Q=P^2$ is parallel. Thus, by applying Lemma 3.1, we
see that $N$ is locally the Riemannian product
$N_{1}\times\ldots\times N_k$ of leaves of distributions
defined by eigenvectors of $Q$ and moreover each
eigenvalue $\lambda_i$ is constant on $N$. If an
eigenvalue $\lambda_{i}$ is zero, the corresponding
leaf $N_i$ is totally real. If $\lambda_{i}$ is $-1$,
then $N_i$ is a complex submanifold. If
$\lambda_{i}\not=0,-1$, then because ${\mathcal D}^i$ is
invariant under the endomorphism $P$ and $<PX,Y>=- \lambda_{i}
<X,Y>$ for any $X,Y$ in ${\mathcal D}^i$, we have
$|PX|=\sqrt{-\lambda_{i}}\,|X|$. Thus the Wirtinger angle
$\theta(X)$ satisfying
$\cos\theta(X)=\sqrt{-\lambda_{i}}$, which is a constant
$\not=0,-1.$ Therefore, $N_i$ is a proper slant
submanifold.

 Assume $\lambda_{i}\not= 0$. We put
$P_{i}=P_{\,| TN_{i}}$. Then $P_i$ is nothing but the
endomorphism of $TN_i$ induced from the almost complex
structure $J$. Let $\nabla^i$ denote the Riemannian
connection of $N_i$. Since $N_i$ is totally geodesic in
$N$, we have $$(\nabla^{i}_{X}P_{i})Y=(\nabla_{X}P)Y=0$$
for any $X,Y$ tangent to $N_i$. This shows that if $N_i$
is  a complex submanifold, $N_i$ is a Kaehlerian
manifold. And if $N_i$ is proper slant, then $N_i$ is a
Kaehlerian slant submanifold of $M$ by definition.

The converse can be verified directly.
\vskip.1in

 From Lemma 3.2 we
may obtain the following 
\vskip.1in 
{\bf Proposition 3.3.}
{\it Let $N$ be an irreducible submanifold of an almost
Hermitian manifold $M$. If $N$ is neither complex nor
totally real, then $N$ is a Kaehlerian slant submanifold
if and only if the endomorphism $P$ is parallel, that is,
$\nabla P=0.$}
 \vskip.1in 
{\bf Theorem 3.4.} {\it Let $N$
be a surface in an almost Hermitian manifold $M$. 
Then the
following three statements are equivalent:}

(1) {\it $N$ is neither totally real nor complex
in $M$ and $\nabla P =0$, that is, $P$ is parallel;}

(2) $N$ {\it is a Kaehlerian slant surface;}

(3) $N$ {\it is a proper slant surface.}
\vskip.1in
{\bf Proof.} Since every proper slant submanifold is of
even dimension, Lemma 3.2 implies that if the
endomorphism $P$ is parallel, then $N$ is a Kaehlerian
surface, or a totally real surface, or a Kaehlerian slant
surface. Thus, if $N$ is  neither totally
real nor complex, then statements (1) and (2) are
equivalent by definition. It is obvious that (2)
implies (3). Now, we prove that (3) implies (2).
Let $N$ be a proper slant surface in $M$ with slant angle
$\theta$. If we choose an orthonormal frame $e_{1},e_2$
tangent to $N$ such that

$$Pe_{1}=(\cos\theta)e_{2},\,\,\,\,Pe_{2}=
-(\cos\theta)e_{1}.$$

\noindent then we have

$$(\nabla_{X}P)e_{1}=\cos\theta(\omega_{2}^{1}(X)+\omega_{1}^{2}
(X))e_{1}.$$

\noindent Since $\omega_{1}^{2}=\omega_{2}^{1},$ we obtain
$\nabla P=0$.

 \vskip.1in For submanifolds of a
Kaehlerian manifold we have the following general lemma.
 \vskip .1in {\bf Lemma 3.5} {\it Let $N$ be a
submanifold of a Kaehlerian manifold $M$. Then }
\vskip.1in
(i) {\it For any vectors $X,Y$ tangent to $N$, we have

$$(\nabla_{X}P)Y =
th(X,Y)+A_{FY}X.\leqno(3.2)$$
 
 Hence $\nabla P = 0$ if and only if 
$A_{FX}Y =A_{FY}X$
for any $X,Y$ tangent of $N$.}
\vskip.1in
(ii) {\it For any  vectors $X,Y$ tangent to $N$, we have

$$(\nabla_{X}F)Y =
fh(X,Y)-h(X,PY)\leqno(3.3)$$
 
Hence $\nabla F = 0$ if and only if $A_{f\xi}X =
-A_{\xi}(PX).$
for any normal vector $\xi$ and tangent vector $X$.}

\vskip.1in

{\bf Proof.} Since $M$ is Kaehlerian, $J$ is parallel.
Thus, by applying the formulas of Gauss and Weingarten
and using formulas (1.1) and (1.5), we may obtain
$$(\nabla_{X}P)Y =
\nabla_{X}(PY)-P(\nabla_{X}Y)=th(X,Y)+A_{FY}X,\leqno(3.2)'$$
and $$(\nabla_{X}F)Y=D_{X}(FY)-F(\nabla_{X}Y)=
fh(X,Y)-h(X,PY).\leqno(3.3)'$$
Thus, $P$ is parallel if and only if 
 $$<th(X,Y)+A_{FY}X,Z>=0$$
which is equivalent to 
$$<A_{FY}X,Z>=-<th(X,Y),Z>=<h(X,Y),FZ>=$$
$$=<A_{FZ}X,Y>=
<A_{FZ}Y,X>.$$
This proves statement (i).

Statement (ii) follows easily from $(3.3)'$.
\vskip.1in
{\bf Remark 3.1.} If $N$ is either a totally real or
complex submanifold of a Kaehlerian manifold, then
$\nabla P = \nabla F =0$, automatically.
 \vskip.1in
Combining Theorem 3.4 and Lemma 3.5 we obtain the following
charaterization of slant surfaces in terms of Weingarten 
map.
 \vskip.1in {\bf Corollary 3.6.} {\it Let $N$ be a
surface in a Kaehlerian manifold $M$. Then $N$ is slant if
and only if $A_{FY}X=A_{FX}Y$ for any $X,Y$ tangent to
$N$.} \vskip.1in
Let $N$ be a slant surface in the complex number space
${\bf C}^2$ with slant angle $\theta$. For a unit tangent
vector field $e_1$ of $N$, we put 
$$e_{2}=(\sec \theta)Pe_{1},\hskip.2in e_{3}=(\csc
\theta)Fe_{1},\hskip.2in
e_{4}=(\csc\theta)Fe_{2}.\leqno(3.4)$$
Then $e_{1}=-(\sec\theta)Pe_{2}$, and
$e_{1},e_{2},e_{3},e_{4}$ form an orthonormal  frame such
that $e_{1},e_2$ are tangent to $N$ and $e_{3},e_4$ are
normal to $N$. As before we put
$$h_{ij}^{r}=<h(e_{i},e_{j}),e_{r}>,\hskip.3in
i,j=1,2;\hskip.2in r=3,4.\leqno(3.5)$$

Let $G$ and $G^D$ denote the {\it Gauss curvature\/} and
the {\it normal curvature\/} of $N$ in ${\bf C}^2$,
respectively. Then we have
$$G=h_{11}^{3}h_{22}^{3}-(h_{12}^{3})^{2}+h_{11}^{4}h_{22}^{4}
-(h_{12}^{4})^{2}\leqno(3.6)$$ and
$$G^{D}=h_{11}^{3}h_{12}^{4}+h_{12}^{3}h_{22}^{4}-h_{12}^{3}
h_{11}^{4}-h_{22}^{3}h_{12}^{4}.\leqno(3.7)$$

From Corollary 3.6 we obtain the following 
\vskip.1in
{\bf Theorem 3.7.} {\it If $N$ is a slant surface in
 ${\bf C}^2$, then $G=G^D$, identically.}
\vskip.1in
{\bf Proof.} Let $N$ be a slant surface in ${\bf C}^2$.
Then Corollary 3.6 implies $A_{FY}X=A_{FX}Y$ for any
vectors $X,Y$ tangent to $N$. Let $e_{1},e_{2},e_{3},e_{4}$
be an orthonormal frame satisfying (3.4). Then we have
$$h_{12}^{3}=h_{11}^{4},\,\,\,\,\,\,h_{22}^{3}=
h_{12}^{4}.$$
Therefore, by (3.6) and (3.7), we obtain $G=G^D$.

\vskip.1in
In the remaining part of this section we mention some
properties of the normal-bundle valued 1-form $F$.
In order to do so, we recall the following definition.
\vskip.1in
{\bf Definition 3.1.} Let $N$ be a submanifold of a
Riemannian manifold $M$. Then $N$ is called a {\it minimal
submanifold\/} if $\,tr\, h = 0,\,$ identically. And it is
called {\it austere\/} (cf. [HL1]) if for each normal
vector $\xi$ the set of eigenvalues of $A_{\xi}$ is
invariant under multiplication by $-1$; this is equivalent
to the condition that all the invariants of odd order of
the Weingarten map at each normal vector of $N$ vanish
identically. \vskip.1in
 Of course every austere submanifold is a minimal
submanifold. 
\vskip.1in
{\bf Theorem 3.8.}  {\it Let $N$ be a proper slant
submanifold of a Kaehlerian manifold $M$. If $\nabla F
=0$, then $N$ is autere.}
\vskip.1in
{\bf Proof.} Let $N$ be a proper slant submanifold of a
Kaehlerian manifold $M$. If $\nabla F=0$, then we have
from formula (3.3)
 $$fh(X,Y)=h(X,PY).$$
Let $X$ be any unit eigenvector of $Q=P^2$ with
eigenvalue $\lambda \not= 0$.
Then $X_{*}=PX/\sqrt{-\lambda}$ is a unit vector
perpendicular to $X$. Thus, we have
$$h(X,X)=h(PX,PX)/\lambda = -h(X_{*},X_{*})$$
which implies that $N$ is autere.
\vskip.1in
If the ambient space $M$ is a complex-space-form, then we
have the following reduction theorem.
\vskip.1in
 {\bf Theorem 3.9.} {\it Let $N$ be an n-dimensional proper
slant submanifold of a complex
$m$-dimensional complex-space-form $M^{m}(c)$ with constant
holomorphic sectional curvature $c$. If $\nabla F =0$, then
$N$ is contained in a complex $n$-dimensional complex
totally geodesic submanifold of  $M^{m}(c)$ 
as an austere submanifold.} 
 \vskip.1in  
{\bf Proof.} Let $N$ be an  $n$-dimensional
proper slant submanifold of ${\bf C}^m$. Assume that
$\nabla F=0$. Then the normal bundle $T^{\perp}N$ has the
following orthogonal direct decomposition:
$$T^{\perp}N=F(TN)\oplus \nu,\,\,\,\nu_{p}\perp F(T_{p}N)$$
for any point $p\in N$. For any vector field $\xi$
in $\nu$ and any vector fields $X,Y$ in $TN$, we have
$$<A_{J\xi}X,Y>=<h(X,Y),J\xi>=<{\tilde\nabla}_{X}Y,J\xi>=$$
$$=-<PY,A_{\xi}X>+<FY,D_{X}\xi>,$$
from which we find
$$<D_{X}(FY),\xi>=-<A_{\xi}(PY)+A_{J\xi}Y,X>.\leqno(3.8)$$

On the other hand, for any $\xi$ normal to $N$, if we
denote by $t\xi$ and $f\xi$ by using (1.5), then Lemma
3.5 gives
$$A_{f\xi}Y+A_{\xi}(PY)=0.\leqno(3.9)$$
Since $f=J$ on the normal subbundle $\nu$, formulas (3.8)
and (3.9) imply
$<D_{X}(FY),\xi>=0$ for any $\xi$ in $\nu$. From this we
conclude that  the normal subbundle $F(N)$ is a
parallel normal subbundle.

Next, we claim that the first normal subbundle $\,Im\,h\,$
is contained in $F(TN)$. This can be proved as follows.

Since $\nabla F=0$, statement (ii) of Lemma 3.5 implies
$$<h(X,Y),J\xi>=-<h(X,PY),\xi>$$
for any normal vector $\xi$ in $\nu$. Thus, for any
eigenvector $Y$ of the self-adjoint endomorphism $Q$
with eigenvalue $\lambda$ and any normal vector $\xi$ in
$\nu$, we have
$$<h(X,Y),\xi>=-\lambda <h(X,Y),\xi>.$$
Since $N$ is  a proper slant submanifold, $-1<\lambda <
0$. Thus, we obtain $Im\,h \subset F(TN)$. Consequently,
by applying the reduction theorem, we
obtain the result.
\vskip.1in
In particular if the ambient space $M$ is the complex
number space ${\bf C}^m$, then we have the following
\vskip.1in 
 {\bf Theorem 3.10.} {\it Let $N$ be an n-dimensional
proper slant submanifold of ${\bf C}^m$. If $\nabla F =0$,
then $N$ is contained in a complex linear subspace ${\bf
C}^n$ of ${\bf C}^m$ as an austere submanifold.}  
\vskip.1in For
surfaces in a real 4-dimensional Kaehlerian manifold we
have the following 
\vskip.1in 
{\bf Theorem 3.11.} {\it Let
$N$ be a surface in a real 4-dimensional Kaehlerian
manifold $M$. Then $\nabla F=0$ if and only if either $N$
is a complex surface, or a totally real surface, or a
minimal proper slant surface of $M$.} \vskip.1in
{\bf Proof.} ($\Rightarrow$) Let $N$ be a surface in $M$
with $\nabla F=0$. Assume that $N$ is neither complex nor
totally real.  Then both $P: T_{x}N \rightarrow T_{x}N$
and $F: T_{x}N \rightarrow T_{x}^{\bot}N$ are surjective.
Denote by $\theta$ the Wirtinger angle. Define
$e_{1},e_{2}, e_{3},e_{4}$ by (3.4). Then, by using
$J^{2}=-I$, we have
$$te_{3}=-\sin \theta e_{1},\, te_{4}=- \sin \theta
e_{2},\, fe_{3}=-\cos\theta e_{4},\,  fe_{4}=\cos\theta
e_{3}.\leqno(3.10)$$
Since $F$ is parallel, Lemma 3.5 implies
$A_{f\xi}X=-A_{\xi}(PX)$. Therefore, we find
$$A_{Fe_{1}}e_{2}=(\csc \theta)^{-1}\sec \theta
A_{e_{3}}(Pe_{1})=$$ $$ =(\tan\theta)A_{e_{3}}(Pe_{1})=
-\tan \theta A_{fe_{3}}e_{1}$$
$$ = \sin \theta A_{e_{4}}e_{1}
= A_{Fe_{2}}e_{1}.$$ 

\noindent Therefore, by applying Corollary 3.6,
$N$ is slant. Furthermore, by applyinf Theorem 3.8, we know
that $N$ is minimal.

($\Leftarrow$) It is clear that if $N$ is a complex or
totally real surface in $M$, then $F$ is parallel.
Therefore, we may assume that $N$ is a minimal proper slant
surface in $M$. We choose  $e_{1},  e_{2},  e_{3}, e_{4}$
according to (3.4). Then by Corollary 3.6, we have 
$$h_{11}^{3}=-h_{22}^{3}=-h_{12}^{4},\,\,\,h_{12}^{3}=
h_{11}^{4}=-h_{22}^{4}.\leqno(3.11)$$
From (3.11) and direct computation we may prove that
$A_{f\xi}X=-A_{\xi}(PX)$ for any tangent vector $X$ and
normal vector $\xi$ of $N$. Therefore, by applying Lemma
3.5, we conclude  $F$ is parallel.

\vfill\eject
\vskip.3in

\noindent \S 4. MINIMAL SLANT SURFACES AND TOTALLY REAL
SURFACES.
\vskip.2in
In this section we want to establish some relations between
minimal slant surfaces and totally real surfaces in a
Kaehlerian manifold, in parti-cular, in ${\bf C}^2$.

Let $N$ be a proper slant surface with slant angle
$\theta$ in a real 4-dimen-sional Kaehlerian manifold $M$.
Let $e_1$ be a local vector field tangent to $N$. We choose
a canonical orthonormal local frame  $e_{1},e_{2},
e_{3},e_{4}$  defined by
$$e_{2}=(\sec \theta)Pe_{1},\hskip.2in e_{3}=(\csc
\theta)Fe_{1},\hskip.2in
e_{4}=(\csc\theta)Fe_{2}.\leqno(4.1)$$

\noindent We call such an orthonormal frame $e_{1},e_{2}, e_{3},
e_{4}$ an {\it adapted slant frame.\/} 

For an adapted
slant frame we have
$$te_{3}=-\sin \theta e_{1},\, te_{4}=- \sin \theta
e_{2},\, fe_{3}=-\cos\theta e_{4},\,  fe_{4}=\cos\theta
e_{3}.\leqno(4.2)$$

As before we put $$De_{r}=\sum_{s} \omega_{r}^{s}\otimes
e_{s},\,\,\,\nabla e_{i} =\sum_{j}
\omega_{i}^{j}\otimes e_{j},$$
$$h=\sum_{r} h^{r}e_{r},\,\,\,\,i,j=1,2;\,\,\,r,s=3,4.$$

We have the following 
 \vskip.1in
{\bf Lemma 4.1.} {\it Let $N$ be a proper slant surface
in a real 4-dimensional Kaehlerian manifold $M$. Then,
with respect to an adapted slant frame, we have}
$$\omega_{3}^{4}-\omega_{1}^{2}=-\cot \theta \{(tr
h^{3})\omega^{1}+(tr h^{4})\omega^{2}\},\leqno(4.3)$$
{\it where $\omega^{1},$ $\omega^2$ is the dual frame
of $e_{1},e_{2}$.}
\vskip.1in
 {\bf Proof.} Since $J$ is parallel, we
have $$D_{X}(FY)-F(\nabla_{X}Y)=fh(X,Y)-h(X,PY).$$

\noindent Thus, we find
$$D_{e_{1}}e_{3} =D_{e_{1}}(\csc \theta Fe_{1}) =
(\csc\theta)De_{e_{1}}(Fe_{1})=$$
$$=(\csc\theta)\{F(\nabla_{e_{1}}e_{1})+fh(e_{1},
e_{1})-h(e_{1},Pe_{1})\}=$$ 

$$=(\csc\theta)\{ \omega_{1}^{2}(e_{1})Fe_{2}
+h_{11}^{3}fe_{3}+$$

$$+h_{11}^{4}fe_{4}-\cos\theta
(h_{12}^{3}e_{3}+h_{12}^{4}e_{4})\}=$$

$$=\omega_{1}^{2}(e_{1})e_{4} - (\cos\theta)(tr
h^{3})e_{4}.$$ 

\noindent This implies
$$\omega_{3}^{4}(e_{1})-\omega_{1}^{2}(e_{1}) =
-(\cot\theta)(tr h^{3}).$$

\noindent Similarly, we may obtain
$$\omega_{3}^{4}(e_{2})-\omega_{1}^{2}(e_{2}) =
-(\cot\theta)(tr h^{4}).$$

\noindent These prove the lemma.
\vskip.1in
{\bf Theorem 4.2.} {\it Let $N$ be a proper slant surface
in a real 4-dimensional Kaehlerian manifold $(M,J,g)$. If
there exists a compatible  complex structure $J_1$ such
that $N$ is totally real with respect to the Kaehlerian
manifold $(M,J_{1},g)$, then $N$ is minimal in $M$.}
\vskip.1in
{\bf Proof.} Since $N$ is assumed to be totally real in
$(M,J_{1},g)$, there exists a function $\varphi$ such that
$$e_{3}=(\cos\varphi)J_{1}e_{1}+(\sin\varphi)J_{1}e_{2}
\leqno(4.4)$$ and
$$e_{4}=-(\sin\varphi)J_{1}e_{1}+(\cos\varphi)J_{1}e_{2}.\leqno(4.5)$$

Since $J_{1}$ is parallel, these imply
$$\omega_{3}^{4}(X)=<{\tilde \nabla}_{X}e_{3},e_{4}>=$$
$$=<-(\sin\varphi)(X\varphi)J_{1}e_{1}+(\cos\varphi)
(X\varphi)
J_{1}e_{2}+(\cos\varphi)\omega_{1}^{2}(X)J_{1}e_{2}+$$ $$
+(\sin\varphi)\omega_{2}^{1}(X)J_{1}e_{1}, -(\sin\varphi)
J_{1}e_{1}+(\cos\varphi)J_{1}e_{2}>.$$

\noindent Therefore, we have
$$\omega_{3}^{4}(X)=\sin^{2}\varphi (X\varphi)+$$
$$+\cos^{2}\varphi (X\varphi) + \cos^{2}\varphi
\omega_{1}^{2}(X)-\sin^{2}\varphi \omega_{2}^{1}(X).$$

\noindent This implies
$$\omega_{3}^{4}-\omega_{1}^{2} =d\varphi.\leqno(4.6)$$

\noindent Combining (4.4) and Lemma 4.1 we obtain
$$\cot\theta \,\{(tr h^{3})\omega^{1}+(tr
h^{4})\omega^{2}\} = -d\varphi.$$

\noindent Also from (4.4) and (4.5) we find
$$h_{11}^{3}=-<{\tilde \nabla}_{e_{1}}e_{3},e_{1}>=$$
$$=-<(\cos \varphi){\tilde \nabla}_{e_{1}}(J_{1}e_{1})+
(\sin \varphi){\tilde
\nabla}_{e_{2}}(J_{1}e_{2}),e_{1}>=$$
$$<\cos\varphi)h(e_{1},e_{1})+(\sin\varphi)h(e_{1},e_{2}),
J_{1}e_{1}>= $$
$$=<(\cos\varphi)h(e_{1},e_{1})+(\sin\varphi)h(e_{1},e_{2}),
\cos\varphi e_{3}-\sin\varphi e_{4}>=$$
$$=\cos^{2}\varphi h_{11}^{3}-\sin^{2}\varphi h_{22}^{3}.$$

\noindent This implies
$$\sin^{2}\varphi (h_{11}^{3}+h_{22}^{3})=0.$$

\noindent Similarly, we have
$$\sin^{2}\varphi (h_{11}^{4}+h_{22}^{4})=0.$$

Let $U=\{x\in N: H(x) \not= 0 \}.$ Then $U$ is an open
subset of $N$. If $U \not= \emptyset$, then $\varphi \equiv
0$ (mod $\pi$) on $U$. Thus, $$\cos\theta \{(tr
h^{3})\omega^{1} +(tr h^{4})\omega^{2}\}=-d\varphi =0$$ on
$U$, which implies $\cos \theta =0$. This implies that $N$
is totally real in  $(M,J,g)$ which is a contradiction.
Thus $U = \emptyset,$ that is, $N$ is minimal.
\vskip.1in
{\bf Theorem 4.3.} {\it Let $N$ be a proper slant surface
in ${\bf C}^2$. Then $N$ is
minimal if and only if 
there exists a compatible almost
complex structure $J_{1}$ on $E^4$ such that $N$ is
totally real in $(E^{4},J_{1})$.} 
\vskip.1in
This Theorem follows from Example 2.2, Theorem 4.3 and
the following [C6, Theorem 5.2]

\vskip.1in
{\bf Theorem 4.4.} {\it Let $N$ be a proper slant
surface in ${\bf C}^2$. Then $N$ is minimal if and
only if there exists a compatible almost complex
structure $J_{2}$ on $E^4$ such that $N$ is a complex
surface in $(E^{4},J_{2})$.}
\vskip.1in

This theorem was proved by studying the relation between
Gauss map and slant immersions. In the next chapter we
will treat these problems by using the notion of Gauss map.

\vfill\eject
\centerline {CHAPTER III}
\vskip.2in
\centerline {\bf GAUSS MAP AND SLANT IMMERSIONS}
\vskip.3in
In this chapter we present some results of the author
and Y. Tazawa [CT1]; in this work  slant surfaces in ${\bf
C}^2$ were studied  from the point of view  of the Gauss
map.  \vskip.3in
 \noindent \S 1. GEOMETRY OF $G(2,4)$.
\vskip.2in
In this section we review the geometry of the real
Grassmannian $G(2,4)$ which consists of all oriented
2-planes in $E^4$.

Let $\{\varepsilon_{1},\varepsilon_{2},\varepsilon_{3},
\varepsilon_{4}\}$ be the canonical orthonormal basis of
$E^4$. Then
$\Psi_{0}=:\varepsilon_{1}\wedge\dots\wedge\varepsilon_{4}$
gives the canonical orientation of $E^4$. Let
$\wedge^{2}(E^{4})^{*}$ denote the 6-dimensional real
vector space with innner product, also denoted by $<\, ,\,
>,$ defined by $$<X_{1}\wedge X_{2},Y_{1}\wedge Y_{2}> =
det\,(<X_{i},Y_{j}>)\leqno(1.1)$$
and extended bilinearly. The two vector spaces
$\wedge^{2}(E^{4})^{*}$ and $(\wedge^{2}E^{4})^{*}$
are identified in a natural way by
$$\Phi(X_{1}\wedge X_{2})=:\Phi(X_{1},X_{2})\leqno(1.2)$$
for any $\Phi \in \wedge^{2}(E^{4})^{*}$. The
Grassmannian $G(2,4)$ was identified with the set
$D_{1}(2,4)$ which consists of all unit decomposible
2-vectors in  $\wedge^{2}E^{4}$ via $\phi:G(2,4)
\rightarrow D_{1}(2,4)$ given by $\phi(V)=X_{1}\wedge
X_{2}$, for any positive orthonormal basis
$\{X_{1},X_{2}\}$ of $V \in G(2,4)$.

The Hodge star operator $* :\wedge^{2}E^{4}
\rightarrow \wedge^{2}E^{4}$ is defined by
$$<*\xi,\eta>\Psi_{0}=\xi\wedge\eta,\leqno(1.3)$$
for any $\xi, \eta \in \wedge^{2}E^{4}$. So, if we
regard an oriented 2-plane $V \in G(2,4)$ as an element in
$D_{1}(2,4)$ via $\phi$, then we have $*V = V^{\bot}$,
where $V^{\bot}$ denotes the oriented orthogonal complement
of the oriented 2-plane $V$ in $E^4$. 

Since $*^{2}=1$ and
$*$ is a self-adjoint endomorphism of $\wedge^{2}E^4$, we
have the following orthogonal decomposition: $$\wedge
^{2}E^{4}= \wedge ^{2}_{+}
E^{4}\oplus\wedge ^{2}_{-}E^{4}\leqno(1.4)$$
of eigenspaces of $*$ with eigenvalues 1 and $-1$,
respectively.

Denote by $\pi_{+}$ and $\pi_{-}$ the natural projections:
$\pi_{\pm} : \wedge ^{2}E^{4} \rightarrow \wedge
^{2}_{\pm}E^{4},$
respectively.

We put
$$\,\,\,\,\,\,\,\,\,\,\eta_{1}= {1\over
\sqrt{2}}(e_{1}\wedge e_{2}+ e_{3} \wedge
e_{4}),\,\,\,\,\eta_{4}= {1\over \sqrt{2}}(e_{1}\wedge
e_{2}- e_{3} \wedge e_{4}),$$ $$\eta_{2}= {1\over
\sqrt{2}}(e_{1}\wedge e_{3}- e_{2} \wedge
e_{4}),\,\,\,\,\eta_{5}= {1\over \sqrt{2}}(e_{1}\wedge
e_{3}+ e_{2} \wedge e_{4}),\leqno(1.5)$$
$$\,\,\,\,\,\,\,\,\,\,\,\eta_{3}= {1\over
\sqrt{2}}(e_{1}\wedge e_{4}+ e_{2} \wedge
e_{3}),\,\,\,\,\eta_{6}= {1\over \sqrt{2}}(e_{1}\wedge
e_{4}- e_{2} \wedge e_{3}).$$ Then
$\{\eta_{1},\eta_{2},\eta_{3}\}$ and $\{\eta_{4},
\eta_{5},\eta_{6}\}$ form canonical orthonormal bases of
$\wedge ^{2}_{+}E^{4}$ and $\wedge ^{2}_{-}E^{4}$,
respectively. We shall orient the spaces $\wedge
^{2}_{+}E^{4}$ and $\wedge ^{2}_{-}E^{4}$ such that these
two bases are positive, that is, they give positive
orientations for the oriented spaces $\wedge
^{2}_{+}E^{4}$ and $\wedge ^{2}_{-}E^{4}$. 

For any $\xi \in D_{1}(2,4)$ we have
$$\pi_{+}(\xi)= {1 \over 2}(\xi +*\xi),
\,\,\,\pi_{-}(\xi) = {1 \over 2}(\xi -*\xi).\leqno(1.6)$$

Denote by $S_{+}^{2}$ and 
$S_{-}^{2}$ the 2-spheres in $\wedge ^{2}_{+}E^{4}$ and $\wedge
^{2}_{-}E^{4}$ centered at the origin with
radius $1/\sqrt{2}$, respectively. Then we have
$$\pi_{+} : D_{1}(2,4) \rightarrow S_{+}^{2},\,\,\, 
\pi_{-} : D_{1}(2,4) \rightarrow S_{-}^{2}\leqno(1.7)$$ and
$$D_{1}(2,4) = S_{+}^{2}\times S_{-}^{2}.\leqno(1.8)$$
\vfill\eject

\vskip.3in
\noindent \S 2. COMPLEX STRUCTURES ON $E^4$.
\vskip.2in
Let ${\bf C}^{2}=(E^{4},J_{0})$ be the complex 2-plane with
the canonical complex structure $J_0$ and the canonical
metric. It is well-known that $J_0$ is an orientation
preserving isomorphism. We denote by $\mathcal J$ the set of
all almost complex structures (or simply called complex
structures) on $E^4$ which are compatible with the inner
product $<\, ,\, >$, that is, $${\mathcal J} = \{J : E^{4}
\rightarrow E^{4} : J \,\,{\rm is \,\,linear,} J^{2} = -I,
\,\,{\rm and} $$ $$<JX,JY>=<X,Y>, \,{\rm for}\,\,{\rm
any}\,\, X,Y \in E^{4}\}.$$

An orthonormal basis $\{e_{1},e_{2},e_{3},e_{4}\}$ on
$E^4$ is called a $J$-{\it basis\/} if we have
$Je_{1}=e_{2}, Je_{3}=e_{4}$. Two $J$-bases of the same
$J$ have the same orientation. Using the canonical
orientation on $E^4$ we divide $\mathcal J$ into two disjoint
subsets of $\mathcal J$: $${\mathcal J}^{+} = \{ J \in {\mathcal J}:
J{\rm - bases \,\, are \,\, positive }\},$$
$${\mathcal J}^{-} = \{ J \in {\mathcal J}: J{\rm -
bases \,\, are \,\, negative }\}.$$

For any $J \in {\mathcal J}$, there exists a unique 2-vector
$\zeta_{J} \in \wedge^{2}E^4$ defined as follows:
$$<\zeta_{J},X\wedge
Y>=-\Omega_{J}(X,Y)=:-<X,JY>,\leqno(2.1)$$
for any $X, Y \in E^4$. In other words, $\zeta_{J}$ is
nothing but the metrical dual of $-\Omega_{J}$, where 
$\Omega_{J}$ is the {\it Kaehler form\/}  associated with
$J$.
\vskip.1in
{\bf Lemma 2.1.} {\it The mapping $\zeta : {\mathcal J}
\rightarrow \wedge^{2}E^{4}$ defined by $ \zeta
(J)=\zeta_{J}$, gives rise to two bijections:
$$\zeta^{+} : {\mathcal J}^{+} \rightarrow
S_{+}^{2}(\sqrt{2}),\,\,\, \zeta^{-} : {\mathcal J}^{-}
\rightarrow S_{-}^{2}(\sqrt{2}),\leqno(2.2)$$
where $S_{+}^{2}(\sqrt{2})$ and $S_{-}^{2}(\sqrt{2})$ are
the 2-spheres centered at the origin with radius
$\sqrt{2}$ in $\wedge^{2}_{+}E^{4}$ and
$\wedge^{2}_{-}E^{4}$, respectively.}
\vskip.1in
{\bf Proof.} If $J\in {\mathcal J}^{+}$ and
$e_{1},e_{2},e_{3}, e_{4}$ is a $J$-basis, then
$e_{1},e_{2},e_{3},e_{4}$ is positive. From (2.1) we have
$$\zeta_{J}=e_{1}\wedge e_{2}+e_{3}\wedge e_{4}.$$
Thus, we have $\zeta_{J}\in S_{+}^{2}(\sqrt{2}).$
Similarly, if $J\in {\mathcal J}^-$, then we have
$\zeta_{J}\in S_{-}^{2}(\sqrt{2}).$

If $J$ and $J'$ are two distinct compatible complex
structures on $E^4$, then their corresponding Kaehler
forms $\Omega_J$ and $\Omega_{J'}$ are distinct. Thus, 
 $\zeta_{J}\not= \zeta_{J'}.$ This proves the
injectivity of $\zeta$.

Conversely, for any element $\xi \in S_{+}^{2}(\sqrt{2})$,
we have ${1\over 2}\xi \in S_{+}^2$. Since
$\pi_{+}:D_{1}(2,4) \rightarrow S_{+}^2$ is isomorphic,
there exists an element $V \in D_{1}(2,4)$ such that
$\pi_{+}(V)= {1\over 2}\xi$. Let $e_{1},e_{2},e_{3},e_{4}$
be a positive orthonormal basis of $E^4$ such that
$e_{1}\wedge e_{2}=V$. Define $J\in {\mathcal J}^+$ such that
$Je_{1}=e_{2},Je_{3}=e_{4}.$ Then we have $\zeta_{J}=\xi$.

This completes the proof of the Lemma.
 \vskip.1in
 By applying Lemma 2.1 we may  make the following
identifications via $\zeta, \zeta^{+}$ and $\zeta^-$,
respectively: $${\mathcal J}^{+} {\cong}
S_{+}^{2}(\sqrt{2}),\,\,
 {\mathcal J}^{-} {\cong} S_{-}^{2}(\sqrt{2}),\,\, 
{\mathcal J} {\cong} S_{+}^{2}(\sqrt{2}) \cup 
S_{-}^{2}(\sqrt{2})$$

For any $V \in G(2,4)$ and for any $J \in {\mathcal J}$, we
choose a positive orthnormal basis $\{e_{1},e_{2}\}$ of
$V$ and put
$$\alpha_{J}(V) =\cos^{-1}(<Je_{1},e_{2}>).\leqno(2.3)$$
Then $\alpha_{J}(V) \in [0,\pi ]$ and $\theta(X)=min \,\{
\alpha_{J}(V),\pi - \alpha_{J}(V)\}.$ A 2-plane $V \in
G(2,4)$ is said to be {\it $\alpha$-slant\/} if
$\alpha_{J}(V)=\alpha$, identically.

If $N$ is an oriented surface in ${\bf C}^2$, then $N$
has a unique complex structure determined by its
orientation and its induced metric. With respect to the
angle $\alpha_{J},\, (J \in {\mathcal J}),$ we have
$$N \,\,{\rm is \,\, holomorphic \,}\Longleftrightarrow
\alpha_{J}(TN) \equiv 0,$$
$$N \,\,{\rm is \,\, antiholomorphic \,}\Longleftrightarrow
\alpha_{J}(TN) \equiv \pi,$$
$$N \,\,{\rm is\,\, totally \,\, real} \Longleftrightarrow
\alpha_{J}(TN) \equiv {\pi \over 2}.$$

The following lemma obtained in [CT1] establishes the
fundamental relations between slant angle and the
projections $\pi_+$ and $\pi_-$.
\vskip.1in
{\bf Lemma 2.2.} 

(i) {\it If $J \in {\mathcal J}^{+}$, then
$\alpha_{J}(V)$ is the angle between $\pi_{+}(V)$ and
$\zeta_J$ and }

(ii) {\it  If $J \in {\mathcal J}^{-}$, then
$\alpha_{J}(V)$ is the angle between $\pi_{-}(V)$ and
$\zeta_J$}.
\vskip.1in
{\bf Proof.} (i) If $J\in {\mathcal J}^+$, then we have
$$\cos (\alpha_{J}(V))=-\Omega_{J}(V)=<\zeta_{J},V>=$$
$$=<\zeta_{J},\pi_{+}(V)+\pi_{-}(V)>=<\zeta_{J},
\pi_{+}(V)>.$$

\noindent Since $||\zeta_{J}||=\sqrt{2}$ and
$||\pi_{+}(V)||=1/\sqrt{2}$, this implies $$\alpha_{J}(V)=
\angle (\pi_{+}(V),\zeta_{J}).$$

(ii) can be proved in a similar way.
\vskip.1in 
For any $a \in [0,\pi ]$ and for any $J \in
{\mathcal J}$, we define

$$G_{J,a}= \{ V \in G(2,4) :
\alpha_{J}(V)=a\},\leqno(2.4)$$

\noindent  which is equivalently to
say that $G_{J,a}$ is the set consisting of all oriented
$a$-slant oriented 2-planes in $E^4$ with respect to the
complex structrure $J$. 

For any $a \in [0,\pi]$ and any $V \in G(2,4)$, we
define $${\mathcal J}_{V,a}= \{ J \in {\mathcal J} :
\alpha_{J}(V)=a\}.\leqno(2.5)$$ 

We put $${\mathcal
J}_{V,a}^{+} = {\mathcal J}_{V,a} \cap {\mathcal J}^{+},\,\, {\mathcal
J}_{V,a}^{-} = {\mathcal J}_{V,a} \cap {\mathcal J}^{-}.$$ 

From
Lemma 2.2 we may obtain [CT1]
\vskip.1in
{\bf Lemma 2.3.} 

(1) {\it If $J \in {\mathcal J}^+$, then
$G_{J,a}= S_{J,a}^{+} \times S_{-}^2$, where $S_{J,a}^+$
is the circle on $S_{+}^2$ consisting of 2-vectors
which makes constant angle $a$ with $\zeta_J$.}

(2) {\it If $J \in {\mathcal J}^-$, then
$G_{J,a}= S_{+}^{2}  \times S_{J,a}^{-}$, where
$S_{J,a}^-$ is the circle on $S_{-}^2$ consisting of
2-vectors which makes constant angle $a$ with $\zeta_J$.}

 (3) {\it Via the identification given by Lemma 2,1,
${\mathcal J}_{V,a}^{+}$ is the circle on $S_{+}^{2}(\sqrt{2})$
consisting of 2-vectors in $S_{+}^{2}(\sqrt{2})$ which
makes constant angle $a$ with $\pi_{+}(V)$, and 
${\mathcal J}_{V,a}^{-}$ is the circle on $S_{-}^{2}(\sqrt{2})$
consisting of 2-vectors in $S_{-}^{2}(\sqrt{2})$ which
makes constant angle $a$ with $\pi_{-}(V)$.}
\vskip.1in
This Lemma can be regarded as a generalization of
Proposition 2 of [CM2].
\vskip.1in
 For later use we give the following
\vskip.1in
{\bf Notation.} Let $V$ be an oriented 2-plane in
$G(2,4)$.  $J_{V}^+$ and $J_{V}^-$ are defined
by 
$$J_{V}^{+}= (\zeta^{+})^{-1}(\pi_{+}(V)) \in {\mathcal
J}^{+},\,\,\,\,J_{V}^{-}= (\zeta^{-})^{-1}(\pi_{-}(V)) \in
{\mathcal J}^{-}.\leqno(2.6)$$
\vskip.1in

{\bf Remark. 2.1.} It is easy to see that $J_{V}^+$
(respectively, $J_{V}^-$) defined by (2.6) is the complex
sturcture in ${\mathcal J}^+$  (respectively, in ${\mathcal J}^-$)
such that $V$ is a holomorphic plane with respect to
the complex structure $J$.

\vfill\eject

\vskip.3in
\noindent \S 3. SLANT SURFACES AND GAUSS MAP.
\vskip.2in 
In this section we study the following problems:

\vskip.1in
{\bf Problem 3.1.} {\it Let $N$ be a surface in ${\bf
C}^{2} = (E^{4},J_{0})$. When is $N$ slant in ${\bf
C}^2$?} \vskip.1in 
{\bf Problem 3.2.} {\it Let $N$ be a surface in $E^4$. If
there exists a compatible complex structure $J$ on
$E^4$ such that $N$ is slant in $(E^{4},J)$. How many other
compatible complex structures ${\tilde J}$ on $E^4$ are
there such that $N$ is slant with resepct to these complex
strucutres?} \vskip.1in
Let $f : N \rightarrow E^4$ be an immersion from an
oriented surface $N$ into $E^4$. Denote by $\nu : N
\rightarrow G(2,4)$ be the Gauss map associated with the
immersion $f$ defined by $\nu(p)=T_{p}N,\,\,p\in N$ (or,
equivalently, by $\nu(p)=(e_{1}\wedge e_{2})(p)$). 

We put
$$\nu_{+}
=\pi_{+}\circ\nu,\,\,\,\,\nu_{-}=\pi_{-}\circ\nu.\leqno(3.1)$$
Then we have $$\nu_{\pm}: N \rightarrow G(2,4) \rightarrow
S_{\pm}^{2}.\leqno(3.2)$$ 

We give the following [CT1]
\vskip.1in
{\bf Proposition 3.1.}  {\it Let $f: N \rightarrow
E^4$ be an immersion of an oriented surface $N$ into
$E^4$. Then}

(1) {\it $f$ is slant with respect to a complex structure
$J \in {\mathcal J}^+$ (respectively, $J \ in \,{\mathcal J}^-$) if
and only if $\nu_{+}(N)$ (respectively, $\nu_{-}(N)$) is
contained in a circle on $S_{+}^2$ (respectively, on
$S_{-}^2$).}

(2) {\it  $f$ is $\alpha$-slant with respect to a complex
structure $J \in {\mathcal J}^+$ (respectively, $J \in
\, {\mathcal J}^-$) if and only if $\nu_{+}(N)$ (respectively,
$\nu_{-}(N)$) is contained in a circle $S^{+}_{J,\alpha}$
on $S_{+}^2$ (respectively, $S_{J,\alpha}^{-}$ on
$S_{-}^2$), where $S_{J,\alpha}^{+}$ (respectively,
$S_{J,\alpha}^{-}$) is the circle on $S_{+}^2$
(respectively, $S_{-}^2$) consisting of all 2-vectors
which makes constant angle $\alpha$ with $\zeta_J$. }
\vskip.1in
{\bf Proof.} If $f : N \rightarrow E^4$ is
$\alpha$-slant with respect to a compatible complex
structure $J \in {\mathcal J}^{+}$, then, by (2.5) and Lemma
2.3, we have $\alpha_{J}(T_{p}N) \in S_{J,\alpha}^{+}
\times S_{-}^2$ for any point $p \in N$. Thus,
$\nu_{+}(N)$ is contained in the circle $S_{J,\alpha}^{+}$
on the 2-sphere $S_{+}^2$ consisting of 2-vectors in
$S_{+}^2$ which makes constant angle $\alpha$ with
$\zeta_J$.

Conversely, if $f: N \rightarrow E^4$  is an immersion
such that $\nu_{+}(N)$ is contained in a circle $S^1$ on
the 2-sphere $S_{+}^2$. Let $\eta$ be a vector of length
$\sqrt{2}$ in $\wedge_{+}^{2}E^4$ perpendicular
to the 2-plane in $\wedge_{+}^{2}E^4$ containing $S^1$ .
Then $\eta \in S^{2}_{+}(\sqrt{2})$. By Lemma 2.1, there
is a unique complex structure $J \in {\mathcal J}^+$ such that
$\zeta_{J}=\eta$. It is clear that $S^1$ is a circle
$S_{J,\alpha}^+$ for some constant angle $\alpha$.
Therefore, by Lemma 2.3, the immersion $f$ is
$\alpha$-slant with respect to this compatible complex
structure $J \in {\mathcal J}^+$. 

Similar argument applies to the other cases.
\vskip.1in
Now we may give the main result of this section [CT1].

\vskip.1in
{\bf Theorem 3.2.} 

(1)  {\it Let $f : N \rightarrow E^4$ be a minimal
immersion. If there exists a compatible complex
structure ${\hat J}\in {\mathcal J}^+$ (respectively, ${\hat
J} \in {\mathcal J}^-$) such that the immesion $x$ is slant
with respect to $\,\,\hat J$, then} 

\hskip.3in (1-a) {\it for any $\alpha \in [0,\pi]$, there
is a compatible complex structure $J_{\alpha}\in {\mathcal
J}^+$ (respectively, $J_{\alpha}\in {\mathcal J}^-$) such that
$f$ is $\alpha$-slant with respect to the
complex structure $J_{\alpha}$.}

\hskip.3in (1-b) {\it the immersion $f$ is slant with
respect to any complex structure $J \in {\mathcal J}^+$
(respectively, $J \in {\mathcal J}^-$).}

(2) {\it If $f : N \rightarrow E^4$ is a non-minimal
immersion, then there exist at most two complex
structures $\,\,\pm J^{+} \in {\mathcal J}^+$ and at most two
complex structures $\,\,\pm J^{-}\in {\mathcal J}^-$ such that
the immersion $f$ is slant with respect to them.}
\vskip.1in
{\bf Proof.} (a) Assume that $\, f : N \rightarrow E^{4}
\,$ is a minimal immersion. Then both $\nu_{+}$ and
$\nu_{-}$ are anti-holomorphic (cf. Lawson's book [L1]).
So, both $\nu_{+}$ and $\nu_-$ are open maps if they are
not constant maps. If the immesion $f$ is slant with
respect to a complex structrure $J \in {\mathcal J}^+$
(respectively,  $J \in {\mathcal J}^-$), then $\nu_+$
(respectively, $\nu_-$) cannot be an open map by
Proposition 3.1. Thus $\,\nu_{+}(N)$ (respectively,
$\,\nu_{-}(N)$) is a singleton. since a singleton is
contained in every circle on the 2-sphere $S_{+}^2$
(respectively, $S_{-}^2$), the immersion $f$ is slant with
respect to every complex structure $J \in {\mathcal J}^+$
(respectively, $J \in {\mathcal J}^-$). So, for any constant
$\alpha \in [0,\pi]$, there is a complex structure
$J_{\alpha}$ which makes the immersion $f$ to be
$\alpha$-slant.

(b)  Assume that $f : N \rightarrow E^4$ is a
non-minimal immersion. If the immersion $f$ is
$\alpha$-slant with respect to a complex structure 
$J \in {\mathcal J}^+$ (respectively, $J \in {\mathcal
J}^-$), then $\nu_{+}(N)$ (respectively, $\nu_{-}(N)$) is
contained in the circle $S_{J,\alpha}^+$ (respectively,
$S_{J,\alpha}^-$) and $\nu_{+}(N)$ (respectively,
$\nu_{-}(N)$) must contain an arc of the
circle $S_{J,\alpha}^+$ (respectively, $S_{J,\alpha}^-$),
since otherwise the immersion $f$ is holomophic with
respect to some compatible complex structures on $E^4$
which implies that $N$ is minimal in $E^4$. Therefore, the
complex structures $\pm J^+$ and $\pm J^-$ are the only
possible complex structures on $E^4$ which may make the
immersion  $f$ slant. \vskip.1in
From Theorem 3.2 we obtain the following corollaries of
[CT1] immediately.
\vskip.1in
{\bf Corollary 3.3.} 

(a) {\it  If $f : N \rightarrow {\bf
C}^{2} = (E^{4},J_{0})$ is  holomorphic, then the
immersion $f$ is slant with respect to every complex
structure $J \in {\mathcal J}^+$.}

(b) {\it If $f : N \rightarrow {\bf
C}^{2} = (E^{4},J_{0})$ is  anti-holomorphic, then
the immersion $f$ is slant with respect to every complex
structure $J \in {\mathcal J}^-$.}
\vskip.1in
{\bf Corollary 3.4.} {\it Let $f : N \rightarrow E^4$ be
a minimal immersion. Then the immersion $f$ is slant with
respect to some complex structure $J \in {\mathcal J}^+$ if
and only if $f$ is holomorphic (respectively,
anti-holomorphic) with respect to some complex structure 
structure $J \in {\mathcal J}^+$ (respectively, $J \in {\mathcal
J}^-$). }
\vskip.1in
{\bf Corollary 3.5.} {\it If $f : N \rightarrow E^3$ is a
non-totally geodesic minimal immersion, then $f : 
N \rightarrow E^{3} \subset E^4$ is not slant with
respect to every compatible complex structure on $E^4$.}
\vskip.1in
{\bf Proof.} If $f : N \rightarrow E^{3} \subset E^4$ is
slant, then, by the minimality of $f$,
Theorem 3.2 implies that the immersion $f$ must be a
proper slant immersion with respect to some compatible
complex structure $J$ on $E^4$. Therefore, by Theorem 3.7
of Chapter II, we have $G = G^D=0$, identically. Because
the only flat minimal surfaces in a Euclidean space are
totally geodesic ones, this is impossible.
\vskip.1in

{\bf Definition 3.1.} An immersion $f : N \rightarrow
E^4$ is called {\it doubly slant\/} if it is slant with
respect to some complex structure $J^{+}\in {\mathcal J}^+$ and
at the same time it is slant with respect to another
complex structure $J^{-} \in {\mathcal J}^-$.
\vskip.1in
From Theorem 3.2 we have the following [CT1]
\vskip.1in
{\bf Corollary 3.6.}  {\it Every non-minimal
immersion $f : N \rightarrow E^4$ which is slant with
respect to more than two complex structures are doubly
slant.}
\vskip.1in
{\bf Remark 3.1.} From Theorem 3.2 we also know that {\it
for any immersion $f: N \rightarrow E^4$, exactly one of
the following four cases occurs:}

(a) {\it $f$ is not slant with respect to every compatible
complex structure on $E^4$.}

(b) {\it $f$ is slant with respect to infinitely many
compatible complex structures on $E^4$.}

(c) {\it $f$ is slant with respect to exactly two
compatible complex  structures on $E^4$.}
 
(d) {\it $f$ is slant with respect to exactly four
compatible complex  structures on $E^4$.}

\vskip.1in
Corollary 3.5 shows that if $N$ is a non-totally geodesic
minimal surface in $E^3$, then, by regarding $E^{3}$ as a
linear subspace of $E^4$, $N$  is not slant with respect
to every compatible complex structures on $E^4$. This 
provides us
many examples for case (a) of Remark 3.1. 
\vskip.1in
Here we remark that every totally real immersion of a
2-sphere $S^2$ into ${\bf C}^2$ gives us an example of
surface in $E^4$ which is slant with respect to exactly
two compatible complex structures on $E^4$. This is due the
fact that the Gauss curvture of any Riemannian metric on 
$S^2$ is non-flat.

\vskip.1in
{\bf Example 3.1.} Let $f : E^{3} \rightarrow E^4$ be the
map from $E^3$ into $E^4$ defined by
$$f(x_{0},x_{1},x_{2}) = (x_{1},x_{2},2x_{0}x_{1},
2x_{0}x_{2}).\leqno(3.3)$$
Then $f$ induces an immersion ${\hat f} : S^{2} \rightarrow
E^4$ from the unit 2-sphere $S^2$ into $E^4$, called the
{\it Whitney immersion\/} which has a unique
self-intersection point ${\hat f}(-1,0,0)={\hat
f}(1,0,0)$. It is know that this immersion ${\hat f} :
S^{2} \rightarrow E^4$ is a totally real immersion with
respect to two suitable compatible complex structures on
$E^4$. Moreover, since the surface is non-flat, ${\hat f}$
is a slant immersion 
 with respect to only two
compatible complex structures (cf. Theorem 4.1.)
\vskip.1in 
In Section 1 of Chapter V we will prove that there exist
no compact proper slant submanifolds in any complex
number space ${\bf C}^m$.
 \vskip.1in {\bf Example 3.2.} Let
$N$ be the surface in $E^4$ defined by
$$x(u,v) = (u,v,k\cos v, k\sin v)\leqno(3.4)$$
Then $N$ is the Riemannian product of a line and a
circular helix in a hyperplane $E^3$ of $E^4$. Let $J_{1},
J_{2}$ be the compatible complex structures on $E^4$
defined respectively by $$J_{1}(a,b,c,d) =
(-b,a,-d,c),\,\,\, J_{2}(a,b,c,d) =(b,-a,-d,c).$$
Then $J_{1} \in {\mathcal J}^+$ and $J_{2} \in {\mathcal J}^-$.
Moreover, by direct computation, we can prove that the
surface $N$ is slant with respect to the following
four complex structures: $J_{1}, -J_{1}, J_{2}, -J_{2}$,
with slant angles given by $$\cos^{-1}({{1}\over
{\sqrt{1+k^{2} } } }), \cos^{-1}({{-1} \over
{\sqrt{1+k^{2}}}}),  \cos^{-1}({{-1}\over
{\sqrt{1+k^{2}}}}),  \cos^{-1}({{1}\over
{\sqrt{1+k^{2}}}}),$$ respectively.  \vskip.1in
{\bf Remark 3.2.} In view of Corollary 3.6, it is
interesting to point out that {\it the only doubly slant
minimal immersion from a surface  into a complex
2-plane is the totally geodesic one.}

\vfill\eject
{\sl III-4. DOUBLY SLANT SURFACES IN ${\bf C}^2$}
\noindent{}

\vskip.3in
\noindent \S 4. DOUBLY SLANT SURFACES IN ${\bf C}^2$.
\vskip.2in
As we defined in Section 3 an immersion $f : N
\rightarrow E^4$ is called doubly slant if it is slant
with respect to a complex structure in ${\mathcal J}^+$ and
at the same time it is slant with respect to another
complex structure in ${\mathcal J}^-$. Equivalently, the
immersion $f$ is doubly slant if and only if there exists
an oriented 2-plane $V \in G(2,4)$ such that $f$ is slant
with respect to both $J_{V}^+$ and $J_{V}^-$, where  $J_{V}^+$ 
and $J_{V}^-$ are defined by (2.6).

In this section we give the following [CT1]
\vskip.1in
{\bf Theorem 4.1.} {\it If $\, f : N
\rightarrow E^4$ is a doubly slant 
immersion, then    $$G = G^{D}=0\leqno(4.1)$$
identically.} \vskip.1in
{\bf Proof.} If $f$ is doubly slant, then, by Proposition
3.1, we know that both $\nu_{+}(N)$ and $\nu_{-}(N)$ lie
in some circles on $S_{+}^2$ and $S_{-}^2$, respectively.
Thus, both $(\nu_{+})_*$ and $(\nu_{-})_*$ are singular
maps at every point $p \in N$. Therefore, the result
follows from the following
\vskip.1in
{\bf Lemma 4.2.} {\it For any immersion $f : N
\rightarrow E^4$ of an oriented surface $N$ into $E^4$ we
have} $$det\,(\nu_{+})_{*}={1 \over 2}(G + G^{D}),\,\,\,
det\,(\nu_{-})_{*}={1 \over 2}(G - G^{D}).\leqno(4.2)$$
\vskip.05in
{\bf Proof.} Let  $x : N
\rightarrow E^4$ be an immersion of an oriented surface $N$
into $E^4$.  Denote by $\{e_{1},e_{2}\}$ be a positive
orthonormal basis of $N$. Then at each point $p \in N$
the Gauss map $\nu$ is given by
$$\nu(p) = (e_{1}\wedge e_{2})(p).\leqno(4.3)$$
Thus, for any vector $X$ tangent to $N$, we have 
$$\nu_{*}X = ({\tilde \nabla}_{X}e_{1})\wedge e_{2} +
e_{1}\wedge ({\tilde \nabla}_{X}e_{2})$$
$$= \omega_{1}^{3}(X)e_{3}\wedge e_{2} +
\omega_{1}^{4}(X)e_{4}\wedge e_{2} 
+ e_{1} \wedge \omega_{2}^{3}(X) e_{3} + 
e_{1} \wedge \omega_{2}^{4}(X) e_{4}  $$
$$= -\omega_{1}^{3}(X)e_{2}\wedge e_{3} -
\omega_{1}^{4}(X)e_{2}\wedge e_{4} 
+ \omega_{2}^{3}(X) e_{1} \wedge  e_{3} + 
 \omega_{2}^{4}(X)e_{1} \wedge e_{4} $$
$$={1 \over 2}\{(\omega_{1}^{4}+\omega_{2}^{3})(X)(e_{1}
\wedge e_{3} - e_{2}\wedge e_{4}) +
(-\omega_{1}^{3}+\omega_{2}^{4})(X)(e_{1}
\wedge e_{4} - e_{2}\wedge e_{3}) $$
$$+(-\omega_{1}^{4}+\omega_{2}^{3})(X)(e_{1}
\wedge e_{3} + e_{2}\wedge e_{4}) +
(\omega_{1}^{3}+\omega_{2}^{4})(X)(e_{1}
\wedge e_{4} - e_{2}\wedge e_{3})\}$$
$$={1\over
\sqrt{2}}\{(\omega_{1}^{4}+\omega_{2}^{3})(X)\eta_{2}
+(-\omega_{1}^{3}+\omega_{2}^{4})(X)\eta_{3}+$$
$$+(-\omega_{1}^{4}+\omega_{2}^{3})(X)\eta_{5}
+(\omega_{1}^{3}+\omega_{2}^{4})(X)\eta_{6}\}.$$
Therefore we have
$$(\nu_{+})_{*}= {1\over \sqrt{2}}\{(\omega_{1}^{4}+
\omega_{2}^{3})\eta_{2}
+(-\omega_{1}^{3}+\omega_{2}^{4})\eta_{3}\},\leqno(4.4)$$
and
 $$(\nu_{-})_{*}= {1\over \sqrt{2}}\{(-\omega_{1}^{4}+
\omega_{2}^{3})\eta_{5}
+(\omega_{1}^{3}+\omega_{2}^{4})\eta_{6}\}.\leqno(4.5)$$
This proves Lemma 4.2.
\vskip.1in
{\bf Remark 4.1.} Examples 2.1-2.6 given in Section 2 of
Chapter II are examples of doubly slant surfaces.
\vskip.1in 
{\bf Remark 4.2.} Lemma 4.2 can be found in [HO1]. 
Our proof of Lemma 4.2 is different from theirs. 
\vfill\eject

\vskip.3in
\noindent \S 5. SLANT SURFACES IN  ALMOST
HERMITIAN MANIFOLDS.
\vskip.2in
Let $f: N \rightarrow (M,J)$ be an immersion of a
differentiable manifold $N$ into an almost complex manifold
$(M,J)$. Then a point $p \in N$ is called a {\it complex
tangent point\/} if the tangent plane of $N$ at $p$ is
invariant under the action of the almost complex structure
$J$.

 The
purpose of this section is to prove the following [CT1]
\vskip.1in
{\bf Theorem 5.1.} {\it Let $f: N \rightarrow (M,g,J)$ be
an imbedding of an oriented surface $N$ into a real
4-dimensional almost Hermitian manifold $(M,g,J)$. If the
immersion $f$ has no complex tangent points, then, for any
prescribed angle $\alpha \in [0,\pi]$, thre exists an
almost complex structure ${\tilde J}$ on $M$ satisfying
the following conditions:}

(a) $(M,g,{\tilde J})$ {\it is an almost Hermitian manifold
and}

(b) {\it the immersion $f$ is $\alpha$-slant with respect
to ${\tilde J}$.}
\vskip.1in
{\bf Proof.} ($M,g,J$) has the natural orientation
determined by the almost complex structure $J$ and, at
each point $p \in M$, the tangant space $T_{p}M$
together with the metric $g_p$ is a Euclidean 4-space.
So we may apply the argument given in Sections 1 and 2 of
this chapter.

According to (1.4) the vector bundle $\wedge^{2}(M)$ of
2-vectors on the ambient manifold $M$ is the direct sum of
two vector subbundles:

$$\wedge^{2}(M)=\wedge^{2}_{+}(M)\oplus\wedge^{2}_{-}(M).
\leqno(5.1)$$

 We define two sphere-bundles over $M$ by
$$S^{2}_{+}(M)=\{\xi\in \wedge^{2}_{+}(M) : ||\xi
||={1\over \sqrt{2}}\},$$
$${\bar S}^{2}_{+}(M)=\{\xi\in \wedge^{2}_{+}(M) : ||\xi
||=\sqrt{2}\}.$$

By using Lemma 2.1 we can identify a cross-section
$$\gamma : M \rightarrow {\bar S}^{2}_{+}(M)\leqno(5.2)$$ 
with an almost complex structure $J_\gamma$ on $M$ such
that $(M,g,J_{\gamma})$ is an almost Hermitian manifold.

In the following we denote by $\rho$ the cross-section
corresponding to the almost complex structure $J$ and we
want to construct another cross-section ${\tilde \sigma}$
to obtain the desired almost complex structure ${\tilde
J}$ on the ambient manifold $M$.

We consider the pull-backs of these bundles via the
imbedding $f : N \rightarrow M$, that is,
$$\wedge^{2}_{+}(N)=f^{*}(\wedge^{2}_{+}(M)),\,\,\,\,\,\,
S^{2}_{+}(N)=f^{*}(S^{2}_{+}(M)),\leqno(5.3)$$
$${\bar S}^{2}_{+}(N)=f^{*}({\bar S}^{2}_{+}(M)).$$

The tangent bundle $TN$ of $N$ determines a cross-section
$\tau : N \rightarrow S^{2}_{+}(N)$ defined by
$$\tau(p)=\pi_{+}(T_{p}N)\leqno(5.4)$$
for any point $p \in N$, where $\pi_+$ denotes the natural
projection from $\wedge^{2}(T_{p}M)$ onto
$\wedge_{+}^{2}(T_{p}M)$. Notice that $2\tau$ is a
cross-section of ${\bar S}^{2}_{+}(N)$:

$$2\tau : N \rightarrow {\bar S}^{2}_{+}(N).\leqno(5.5)$$

We denote  $f^{*}\rho$ also by $\rho$ for simplicity.
We have the following cross-section:
$$\rho=f^{*}\rho : N \rightarrow {\bar
S}^{2}_{+}(N).\leqno(5.6)$$

Since the imbedding $f$ is assumed to have no complex
tangent points, 

$$\rho(p)\not= \pm 2\tau(p)\leqno(5.7)$$

\noindent for any point $p \in N.$ Therefore, $\rho(p)$ and
$2\tau(p)$ determine a 2-plane in
$\wedge^{2}_{+}(T_{p}M)$ which intersects the circle
$({\mathcal J}^{+}_{\tau,a})_{p}$ at two points. Where 
$({\mathcal J}^{+}_{\tau,a})_{p}$ is a circle on $({\bar
S}^{2}_{+}(N))_p$  defined in Lemma 2.3 with
$V=T_{p}N$. 

Let $\sigma(p)$ be one of these two points
which lies on the half-great-circle on $({\bar
S}^{2}_{+}(N))_p$ starting from $2\tau(p)$ and passing
through $\rho(p)$. Since $\rho$ and $\tau$ are
differentiable, so is $\sigma$. Thus we obtain a third
cross-section:
$$\sigma : N \rightarrow {\bar S}^{2}_{+}(N)\leqno(5.8)$$
and we want to extend $\sigma$ to a cross-section
${\tilde \sigma}$ of ${\bar S}^{2}_{+}(M).$

For each point $p \in N$, we choose an open neighborhood
$U_p$ of $p \in M$ such that $\sigma_{| U_{p}\cap N}$
can be extended to a cross-section of ${\bar
S}^{2}_{+}(M)$ on $U_p$:
$$\sigma_{p} : U_{p} \rightarrow {\bar S}^{2}_{+}(M)_{|
U_{p}}.\leqno(5.9)$$
Here we identify  the manifold $N$ with its image $f(N)$
of $N$ under the imbedding $f$. We put
$${\mathcal U} = \cup_{p\in N}\, U_{p}\leqno(5.10)$$
and pick  a locally finitely countable refinement
$\{U_{i}\}$ of the open covering $\{U_{p}: p \in N\}$ of
$U$. 

For each $i$ we pick a point $p \in N$ such that
$U_i$ is contained in $U_p$ and we put
$$\sigma_{i}=\sigma_{p | U_{i}}.\leqno(5.11)$$

Let $\phi_{i}$ be a differentiable partition of unity on
${\mathcal U}$ subordinate to the covering $\{ U_{i}\}$. We
define a cross-section ${\bar \sigma}$ of
$\wedge^{2}_{+}(M)_{| {\mathcal U}}$ by
$${\bar \sigma}: {\mathcal U} \rightarrow
\wedge^{2}_{+}(M)_{| {\mathcal U}};
\,\,\,\,{\bar\sigma}=\sum_{i}
\phi_{i}\sigma_{i}.\leqno(5.12)$$

From the constructions of $\sigma_i$ and ${\bar\sigma}$
we have
$${\bar \sigma}_{| N}=\sigma .\leqno(5.13)$$

Since the angle $\angle \,
({\bar\sigma(p)},\rho(p))$ between ${\bar \sigma}(p)$ and
$\rho(p)$ is less that $\pi$ for any point $p \in N$, we
have $${\bar \sigma}(p)\not= 0,\,\,\,\angle \,
({\bar\sigma(p)},\rho(p)) < \pi$$
for any point $p \in N$. By continuity of ${\bar \sigma}$
we can pick an open neighborhood $W$ of $N$ contained in 
${\mathcal U}$ such that
$${\bar\sigma}(q)\not= 0,\,\,\,\,\angle \,
({\bar\sigma(q)},\rho(q))<\pi$$
for any point $q \in W.$ We define a cross-section ${\hat
\sigma}$ of ${\bar S}^{2}_{+}(M)_{| W}$ by

$${\hat \sigma} : W \rightarrow {\bar S}^{2}_{+}(M)_{| W};
\,\,\,\,\,{\hat \sigma}={\bar\sigma}/\sqrt{2}||{\bar
\sigma}||.$$

\noindent Then we have $\angle
\,({\hat\sigma}(q),\rho(q))<\pi$ for any point $q\in M$,
too. 

Finally, we consider the open
covering $\,\{W,M-N\}\,$ of $M$ and local cross-section
${\hat\sigma}$ and $\rho_{| M-N}$ and repeat the same
argument using an partition of unity subordinate to $\{W,
M-N\}$ to obtain a cross-section ${\tilde \sigma}:M
\rightarrow {\bar S}^{2}_{+}(M)$ satisfying ${\tilde
\sigma}_{| N}=\sigma$. Now it is clear that the almost
complex structure ${\tilde J}$ corresponding to ${\tilde
\sigma}$ is the desired almost complex structure. This
completes the proof of the theorem.

\vfill\eject

\vskip1in
\centerline {CHAPTER IV}
\vskip.2in
 \centerline {\bf CLASSIFICATIONS OF SLANT SURFACES}
\vskip.5in
\noindent \S 1.$\,\,$ SLANT SURFACES WITH
PARALLEL MEAN CURVATURE 

VECTOR.
\vskip.2in
The main purpose of this section is to present the
following classification theorem of [C6].
\vskip.1in
{\bf Theorem 1.1.} {\it Let $N$ be a surface in ${\bf
C}^2$. Then $N$ is a slant surface with parallel mean
curvature vector, that is, $DH=0$, if and only if $N$ is
one of the following surfaces:}

(a) {\it an open portion of the product surface of two
plane circles, or}

(b) {\it an open portion of a circular cylinder
which is contained in a hyperplane of ${\bf C}^2$, or}

(c) {\it a minimal slant surface in  ${\bf C}^2$.}

\noindent {\it Moreover, if either case {\rm (a)} or case
{\rm (b)} occurs, then $N$ is a totally real surface in 
${\bf C}^2$.} \vskip.1in
{\bf Proof.} Let $N$ be a slant surface in  ${\bf C}^2$
with parallel mean curvature vector. Then the length of
the mean curvature vector $H$ is constant. If the length
is zero, $N$ is a minimal slant surface. So Case (c)
occurs.  Now assume that $N$ is not minimal in  ${\bf
C}^2$. Then one may choose a unit tangent vector $e_1$
such that $e_{3}=(\csc\theta)Fe_1$ is in the direction of
$H$, where $\theta$ denotes the slant angle of $N$ in 
${\bf C}^2$. Such an $e_1$ can be chosen if we choose
$e_1$ to be in the direction of $-tH$. Since the mean
curvature vector is parallel, $\omega_{3}^{4}=0$. Thus the
normal curvature $G^{D}=0$, identically. Hence, by applying
Theorem 3.7 of Chapter II, we have $G=0$, that is, $N$ is
flat. Let $V=\{ \, p\in N : A_{e_{4}}\not=
0\,\,at\,\,p\,\}$. Then $V$ is an open subset of $N$.

 {\bf Case (i):}  $V=\emptyset$. In this case we have
$Im\,h\subset Span\{ H\}$. Thus, by applying the condition
$DH=0$, we conclude that the surface $N$ lies in a
hyperplane $\,E^{3}\,$ of  ${\bf C}^2$. Since $N$ is a flat
surface in $E^3$ with nonzero constant mean curvature, $N$
is an open portion of a circular cylinder (cf. Proposition
3.2 of [C1, p.118]). 

{\bf Case (ii):}  $V\not=\emptyset$. In this case, let $W$
be a connected component of $V$. Then $e_4$ is a parallel
minimal non-geodesic section on $W$ (cf. [C1]). Since $W$
is flat, Proposition 5.4 of [C1, p.128] implies that $W$ is
an open piece of the product surface of two plane circles.
Since $det\,(A_{e_{4}})$ is a nonzero constant on $W$, by
continuity we have $V=N$. Thus the whole surface is an
open portion of the product surface. 

Finally, if $N$ is an open portion of the product surface
of two plane circles or an open portion of a circular
cylinder contained in a hyperplane of  ${\bf C}^2$, then
$N$ is totally real in $E^4$ with respect to some
compatible complex structure, say $J'$, on $E^4$.
Therefore, by applying Theorem 4.3 of Chapter II, we know
that $N$ must be totally real in  ${\bf C}^2$, since $N$
is non-minimal. This completes the proof of the theorem.
\vskip.1in
By applying Theorem 1.1 we may obtain the following
classification  of parallel slant surfaces, that is,
slant surfaces with parallel second fundamental form [C4].
\vskip.1in
{\bf Theorem 1.2.} {\it Let $N$ be a surface in  ${\bf
C}^2$. Then $N$ is a slant surface in  ${\bf C}^2$ with
parallel second fundamental form, that is, ${\bar
\nabla}h=0,$ if and only if $N$ is one of the following
surfaces:}

(a) {\it an open portion of the product surface of two
plane circles;}

(b) {\it an open portion of a circular cylinder which is
contained in a hyperplane of  ${\bf C}^2$;}

(c) {\it an open protion of a plane in  ${\bf C}^2$}.

\noindent {\it Moreover, if either case {\rm (a)} or case
{\rm (b)} occurs, then $N$ is totally real in  ${\bf
C}^2$.} \vskip.1in
{\bf Proof.} If $N$ is a surface in ${\bf C}^2$ with
parallel second fundamental form, then $N$ has parallel
mean curvature vector. Thus, by Theorem 1.1, it suffices
to prove that slant planes are the only minimal slant
surfaces with parallel second fundamental form. But this
follows from the facts that every surface in  ${\bf
C}^2$ with parallel second fundamental form has constant
Gauss curvature and every minimal surface in  ${\bf
C}^2$ with constant Gauss curvature is totally geodesic.
\vskip.1in
By using Theorem 1.1, we obtain
\vskip.1in
{\bf Corollary 1.3.} {\it Let $N$ be a slant surface in 
 ${\bf C}^2$
with constant mean curvature. Then $N$ is spherical if
and only if $N$ is an open portion of the product surface
of two plane circles.}
\vskip.1in
{\bf Proof.} If $N$ is a spherical surface with constant
mean curvature, then the mean curvature vector of $N$ in 
 ${\bf C}^2$ is parallel. Thus, by applying Theorem 1.1,
$N$ is one of the surfaces mentioned in Theorem 1.1.
Among them, surfaces of type (a) are the only spherical
surfaces in  ${\bf C}^2$. 

The converse is obvious.
\vskip.1in
Similarly we may prove the following
\vskip.1in
{\bf Corollary 1.4.} {\it Let $N$ be a slant surface in
${\bf C}^2$ with constant mean curvature. Then $N$ lies
in a hyperplane of ${\bf C}^2$ if and only if $N$ is
either an open portion of a 2-plane or an open portion
of a circular cylinder.}
  \vskip.1in {\bf Remark 1.1.} In
views of Theorem 1.1, Corollary 1.3 and Corollary 1.4, it
seems to be interesting to propose the following open
problem: \vskip.1in {\bf Problem 1.1.} {\it Classify all
slant surfaces in ${\bf C}^2$ with nonzero constant mean
curvature.}

\vfill\eject

\vskip.3in
\noindent \S 2.  SPHERICAL SLANT SURFACES.
\vskip.2in
In the remaining part of this chapter we want present
some classification theorems obtained in [CT2].

  In this section we want to classify slant surfaces of 
 ${\bf C}^2$ which lie in a hypersphere of  ${\bf C}^2$. 
In order to do so, we need to review the geometry of the
unit hypersphere $S^{3}=S^3(1)$ in  ${\bf C}^2$ centered at
the origin.

 It is known that $S^3$ is the Lie group consisting of all 
unit quaternions $\{ u = a +{\bf i} b+{\bf j} c+{\bf k} d
: ||u|| = 1 \} $ 
 which can also be regarded as a subgroup of the
orthogonal group $O(4)$ in a natural way. Let 1 denote the
identity element of the Lie group $S^3$ given by 
$$1=(1,0,0,0) \in S^{3} \subset E^{4}.\leqno(2.1)$$  

We put
$$X_{1}=(0,1,0,0),\,\,\,X_{2}=(0,0,1,0),\,\,\,
X_{3}=(0,0,0,1) \in T_{1}S^{3}. \leqno(2.2)$$ 

 We denote by
${\tilde X}_{i}, i = 1, 2, 3$, the left-invariant vector
fields obtained from the  extensions of $X_{i}, i = 1, 2
,3$, on $S^3$, respectively. Let $\phi :S^{3} \rightarrow
S^{3}$ be the orientation-reversing isometry defined by
$$\phi (a,b,c,d)=(a,b,d,c).\leqno(2.3)$$

 We recall that the
left-translation $L_p$ and the right-translation $R_p$ on
$S^3$ are isometries which are analogous to the parallel
translations on $E^3$ and they are given by
$${}^t(L_{p}q)=\begin{pmatrix}a&-b&-c&-d\\ b&a&-d&c\\
c&d&a&-b\\ d&-c&b&a\end{pmatrix}\begin{pmatrix}x\\ y\\
z\\ w\end{pmatrix}\leqno(2.4)$$ and
\eject 
$${}^t(R_{p}q)= \begin{pmatrix}a&-b&-c&-d\\ b&a&d&-c\\
c&-d&a&b\\ d&c&-b&a\end{pmatrix}\begin{pmatrix}x\\ y\\ z\\ w
\end{pmatrix}\leqno(2.5)$$

\noindent for $\, p = (a,b,c,d),\,\, q=(x,y,z,w) \in S^{3}
\subset E^{4},\,$ where $\,{}^{t}A \,$ denotes the
transpose of $A$. 

Let $\eta$ denote the unit outer normal of $S^3$ in $E^4$
and $J_{1}$ and  $J_{1}^{-}$ the complex structures on
$E^4$  defined respectively by 
$$J_{1}(a,b,c,d)=(-b,a,-d,c),\leqno(2.6)$$
$$J_{1}^{-}(a,b,c,d)=(-b,a,d,-c).\leqno(2.7)$$

\vskip.1in
{\bf Lemma 2.1.} {\it For any $q \in S^{3}$, we have}
$$(J_{1} \eta) (q)=R_{q*}X_{1},\leqno(2.8)$$
$$(J_{1}^{-}\eta) (q)=L_{q*}X_{1}= {\tilde
X}_{1}(q).\leqno(2.9)$$ {\it Hence $J_{1}\eta$ and
$J_{1}^{-}\eta$ are right-invariant and left-invariant
vector fields on $S^3$, respectively.}
\vskip.1in
{\bf Proof.} Let $q=(a,b,c,d) \in S^3$. Then
$$\eta(q)=(a,b,c,d)\in  T_{q}^{\perp}S^{3},\leqno(2.10)$$
$$(J_{1}\eta)(q)=(-b,a,-d,c),\,\,\,\,\,(J_{1}^{-}\eta)
(q)=(-b,a,d,-c).$$

Let $\gamma \subset S^3$ be a curve on $S^3$ parametrized
by arclength $s$ given by 
$$\gamma(s)=(\cos s,\sin s,0,0).\leqno(2.11)$$
Then we have
$$R_{q^{*}}X_{1}={d\over {ds}}(R_{q}(\gamma(s)))_{|
s=0} =$$ $$= {d\over {ds}}((\cos s + {\bf i}\sin s)
(a+{\bf i}b+{\bf j}c+{\bf k}d))_{| s=0}$$
$$=(-b,a,-d,c)=(J_{1}\eta)(q).$$

Similarly we may prove $L_{q^{*}}X_{1}=(J_{1}^{-}\eta)(q).$
\vskip.1in
 As
in Section 2 of Chapter III, we denote by ${\mathcal J}$ the
set of all complex structures on $E^4$ compatible with the
inner product $<\,,\,>$.  By using the
natural orientation of $E^4$ we  divide $\mathcal J$ into two
disjoint subsets:  $${\mathcal J}^{+} = \{{J\in {\mathcal J}\,|\,
J{\rm -bases\,\, are\,\, positive}}\},$$   $${\mathcal J}^{-}
= \{{J \in {\mathcal J}\,|\, J{\rm -bases \,\, are \,\,
negative}}\}.$$ 

\noindent Thus we have ${\mathcal J}= {\mathcal J}^{+} \cup
{\mathcal J}^{-}$ (cf. Section 2 of Chapter III).
\vskip.1in
{\bf Lemma 2.2.} {\it Let $W\in G(3,4)$ and $V\in G(2,4)$
such that $V\subset W$. Then $V$ is $\alpha$-slant with
respect to a complex structure $J\in {\mathcal J}^{+}$
(respectively, $J\in {\mathcal J}^{-}$) if and only if 
$$<\xi_{V},J\eta_{W}> = -\cos\alpha \hskip.15in 
(respectively, <\xi_{V},J\eta_{W}> =
\cos\alpha),\leqno(2.12)$$
where $\, \xi_V$ and $\,\eta_W$ are positive unit normal
vectors of $V$ in $W$  and of $\, W$ in $E^4$,
respectively.}  \vskip.1in

{\bf Proof.} We choose an orthonormal
$J$-basis $\{e_{1},\ldots,e_{4}\}$ of $E^4$ such that
$$ e_{1}, e_{2} \in W\cap JW,\,\,\,\,
e_{4}=Je_{3}=\eta_{W}.\leqno(2.13)$$
We also choose a positive orthonormal basis
$\{X_{1},X_{2}\}$ of $V$. Let $\zeta_J$ be the 2-vector
defined as before as the metrical dual of $-\Omega_{J} \in
(\wedge^{2}E^{4})^*$, that is, $$<\zeta_{J},X\wedge Y>
= -\Omega_{J}(X,Y) $$ for any $X, Y \in E^4$.  Then,
by formula  (2.3) of Section 2 of Chapter III, we see that
the slant angle $\alpha_{J}(V)$ of $V$ with respect to $J$
satisfies  
$$\cos\alpha_{J}(V) = <\zeta_{J},X_{1}\wedge
X_{2}>$$
$$ = <e_{1}\wedge e_{2} + e_{3}\wedge e_{4},
X_{1}\wedge X_{2}>= <e_{1}\wedge e_{2}, X_{1}\wedge
X_{2}>$$
 $$ =< \pm
e_{3},\xi_{V}> = \mp <J\eta_{W},\xi_{V}>$$
for $J \in {\mathcal J}^{\pm}$. This proves the lemma.

\vskip.1in

Let $f : N \rightarrow S^{3} \subset E^4$ be a spherical
immersion of an oriented surface $N$ into $S^3$ and $\xi$
the positive unit normal of $f(N)$ in $S^3$. 

It is
easy to see that every spherical surface in ${\bf C}^2$
 is non-minimal.
Hence, every spherical surface in ${\bf C}^2$  cannot 
be a complex surface of ${\bf C}^2$. 

 \vskip.1in

{\bf Lemma 2.3.} {\it Let $f: N \rightarrow S^{3}
\subset E^4$ be an immersion of an oriented surface $N$
into the hypersphere $S^3$ of $E^4$. Then the following
three statements hold.}

(i) {\it The immersion $f$ is $\alpha$-slant with respect
to the complex structure $J_1$ if and only if }
$$<\xi(p),J_{1}\eta (f(p))>= -\cos\alpha\hskip.2in
 for \,\,any\,\, p \in N.\leqno(2.14)$$

(ii) {\it The immersion $f$  is $\alpha$-slant with respect
to the complex structure $J_{1}^-$ if and only if }
$$<\xi(p),{\tilde X}_{1}(f(p))> = \cos \alpha \hskip.2in
for \,\,any\,\, p \in M.\leqno(2.15)$$

(iii) {\it The immersion $f$  is $\alpha$-slant with
respect to the complex structure $J_1$ if and only if the
composition $ \,\phi\circ f \,$ is $\alpha$-slant with
respect to the complex structure  $J_{1}^{-}$.} \vskip.1in

	{\bf Proof.} Statement (i) follows from
Lemma 2.2. Statement (ii) follows from Lemma 2.1 and
Lemma 2.2. Finally, the last statement follows from
statements (i) and (ii) and from the fact that $\phi$ is
an isometric involution which reverses the orientation of
$E^4$.
 \vskip.1in
 We define two maps $g_{+}$ and $g_{-}$ from $N$ into the
unit sphere $S^2$ in $T_{1}S^3$ by
$$g_{+}(p)=(L_{\phi(f(p))*})^{-1}(\phi_{*}\xi(p)),\,\,\,\,
 g_{-}(p)=(L_{f(p)*})^{-1}(\xi (p)) \leqno(2.16)$$
for $p \in N$. In fact, $g_{+}$ and $g_{-}$ are the
analogues of the classical Gauss map of a surface in $E^3$
in which  the parallel translations in $E^3$  are
replaced by  the left-translations $L_q$ on $S^3$. 

We
also define a circle $S_{\alpha}^{1}$ for $\alpha \in
[0,\pi ]$ on the unit sphere $S^2$ in $T_{1}S^3$ by

$$S_{\alpha}^{1}= \{ X \in T_{1}S^{3} \, \mid \, \|X\| =
1, \,<X,X_{1}>= - \cos \alpha \}. \leqno(2.17)$$ 

Now we prove the following result of [CT2] which
characterizes spherical slant surfaces in  {\bf C}$^2$. 

\vskip.1in
{\bf Proposition 2.4.} {\it Let  $\,f : N \rightarrow
S^{3} \subset E^4$ be an immersion of an oriented surface
$N$. Then we have}

(i) $f$ {\it is $\alpha$-slant with respect to the
complex structure $J_{1}$ if and only if }
$$g_{+}(N) \subset S_{\alpha}^{1} \subset
T_{1}S^{3}.\leqno(2.18)$$

 (ii) $f$ {\it is $\alpha$-slant with respect to the
complex structure  $J_{1}^{-}$ if and only if} $$g_{-}(N)
\subset S_{\pi -\alpha}^{1} \subset
T_{1}S^{3}.\leqno(2.19)$$  \vskip.1in
 {\bf Proof.} Proposition 2.4 
follows from Lemma 2.1, Lemma 2.3 and the definitions of
$g_+$ and $g_-$. 
\vskip.1in
 Proposition 2.4 can be regarded as the spherical version
of Proposition 3.1 in Section 3 of Chapter III.
\vskip.1in
Concerning the images of $N$ under the spherical Gauss
maps $g_+$ and $g_-$, we give
 the following two simple examples (cf. [CT2]).
\vskip.1in

{\bf Example 2.1.} If $N = S^{1}\times S^1$ is the flat
torus in $E^4$ defined by
$$f(u,v) = {1\over {\sqrt{2}}}(\cos u,\sin u,\cos v,\sin
v),\leqno(2.20)$$ then  the images of $N$ under the
spherical Gauss maps $g_+$ and $g_-$ are the great circles
perpendicular to $X_{1}= (0,1,0,0).$ \vskip.1in

{\bf Example 2.2.}   If $N = S^2$ is the totally geodesic
2-sphere of $S^3$ parametrized by
$$f(u,v) = (\cos u\cos v,\sin u\cos v,\sin v,0),$$ then 
 $$g_{+}(u,v) = (0,-\sin v,-\cos u\cos v,\sin u\cos
v),$$
$$g_{-}(u,v) = (0,\sin v,\sin u\cos v,-\cos u\cos v).$$
Hence, both $g_+$ and $g_-$ are isometries.
\vskip.1in

In order to describe slant surfaces in $S^3$
geometrically, we give the following

\vskip.1in
{\bf Definition 2.1.} Let $c(s)$ be a curve in $S^3$
 parametrized by  arclength $s$ and let
$$c'(s) = \sum_{i=1}^{3} f_{i}(s){\tilde
X}_{i}(c(s)).\leqno(2.21)$$

We call the curve $c(s)$ a {\it helix in\/} $S^3$
{\it with  axis vector field} ${\tilde X}_1$ if 
$$f_{1}(s)=b,\,\,\, f_{2}(s)=a\,\cos (ks+s_{0}), \,\,\,
f_{3}(s)= a\,\sin (ks+s_{0})\leqno(2.22)$$

\noindent for some constants $a, b, k$ and $s_0$ 
satisfying
$$a^{2}+b^{2}=1.\leqno(2.23)$$

We call the curve $c(s)$ a {\it generalized helix in\/}
$S^3$ {\it with axis vector field} ${\tilde X}_1$  if
$$<c'(s),{{\tilde X}_{1}}(c(s))> = constant.\leqno(2.24)$$
\vskip.1in

Helices in $S^3$ defined above are the analogues of
Eucliden helices in the Euclidean 3-space $E^3$.

\vskip.1in 

{\bf Definition 2.2.} We call an immersion $\, f : D
\rightarrow S^3$ of a domain $D$ around the origin
$(0,0)$ of $R^2$ into $S^3$ a {\it helical cylinder in}
$S^3$ if
 $$f(s,t)=\gamma (t)\cdot c(s),\leqno(2.25)$$

\noindent for some   helix $c(s)$ in $S^3$ with axis ${\tilde
X}_1$ satisying $k = - 2/b$ and $ab < 0$, and for some
curve $\gamma (t)$ in $S^3$ which is either a geodesic or
 a curve of constant torsion 1 parametrized by  arclength
such that 

(i) $c(0)=\gamma (0)$, and 

(ii) the osculating
planes of $c(s)$ and of $\gamma (t)$ coincide at $t=s=0$.

We note that the binormal of $c(s)$ is normal to
$f(D)$ in $S^3$. Here we orient the curve $c$ in such a
way that the binormal of $c(s)$ is the positive unit
normal of $f(D)$. 
\vskip.1in

The main purpose of this section is to prove the following
classification theorem for spherical slant surfaces [CT2].
 \vskip.1in

{\bf Theorem 2.5.} {\it Let $f : N \rightarrow S^{3}
\subset {\bf C}^{2} = (E^{4},J_{1})$ be a spherical
immersion of an oriented surface $N$ into the complex
2-plane   {\bf C}$^2$ $=(E^{4},J_{1})$. Then $f$ is a
proper slant immersion if and only if
 $f(N)$ is locally of the form $\{\phi
(\gamma (t)\cdot c(s))\}$ where $\phi$ is the isometry on
$S^3$ defined by (2.3) and $\{\gamma (t) \cdot c(s)\}$ is
a helical cylinder in $S^3$ (cf. Definition 2.2).}

\vskip.1in
 In order to prove this theorem so we need several lemmas.

First, we note that curves in $S^3$ can be described in
terms of the orthonormal left-invariant vector fields
$\{{\tilde X}_{1},{\tilde X}_{2},{\tilde X}_{3}\}$. 
 Let $I$ be an open interval
containing 0 and $c : I \rightarrow S^3$   a curve
parametrized by  arclength $s$. Let ${\bf t}(s), {\bf
n}(s), {\bf b}(s), \kappa(s)$, and $\tau(s)$ be the unit
tangent vector, the unit principal normal vector, the unit
binormal vector, the curvature, and the torsion of $c$ in
$S^3$, respectively. We put
$${\bf t}(s) = \sum_{i=1}^{3} f_{i}(s){\tilde
X}_{i}(c(s)).\leqno(2.26)$$ Then 
$$(f_{1}(s))^{2} + (f_{2}(s))^{2} + (f_{3}(s))^{2} =
1.\leqno(2.27)$$ 
\vskip.1in
Conversely, we have the following
\vskip.1in

{\bf Lemma 2.6.} {\it Let $f_{i}(s), i = 1, 2, 3,$ be
differentiable functions on $I$ satisfying (2.27). Then,
for any point $p_{0} \in S^3$, there exists a curve
$c(s)$ in $S^3$ defined on an open subinterval $I'$
of $I$ containing $0$ and  satisfying (2.26) and the
initial condition $c(0) = p_0$.} \vskip.1in

{\bf Proof.} Considering the curve $L_{p_{0}}^{-1}\circ
c$ instead of $c$, if necessary, we can assume without loss
of generality that $p_{0}=1$. The solution of the following
system of the first order linear differential equations:
$$\begin{pmatrix}x'\cr y'\cr z'\cr w'\end{pmatrix} = 
\begin{pmatrix}x&-y&-z&-w\\
y&x&-w&z\\ z&w&x&-y\\ w&-z&y&x\end{pmatrix}\begin{pmatrix}0\\
 f_{1}\\
f_{2}\\ f_{3}\end{pmatrix} \leqno(2.28)$$

\noindent with the initial condition $(x(0),y(0),z(0),w(0))
= (1,0,0,0)$ satisfies  the condition: $$xx' + yy' + zz' +
ww' = 0$$ and the curve $c(s) = (x(s),y(s),z(s),w(s))$ is
in fact the desired one.
\vskip.1in
Lemma 2.6 guarantees the existence of helices in $S^3$.
\vskip.1in

{\bf Lemma 2.7.} {\it The following two statements are
equivalent:}
\vskip.05in
(i) {\it The curve $c(s)$ is a helix in $S^3$ with axis
vector ${\tilde X}_1$ satisfying}
$$f_{1}(s) = b,\leqno(2.29)$$
$$f_{2}(s) = a\cos (-{2\over b}s + s_{0}),\leqno(2.30)$$
$$f_{3}(s) = a\sin (-{2\over b}s + s_{0}),\leqno(2.31)$$
$$a^{2} + b^{2} = 1,\,\,\,\,\, ab < 0.\leqno(2.32)$$

(ii) {\it The curve $c(s)$ is a curve in $S^3$ satisfying}

$$\tau (s) \equiv -1,\leqno(2.33)$$
$$<{\bf b}(s),{\tilde X}_{1}(c(s))> \equiv a.\leqno(4.9)$$
$$a \not= \pm 1, 0.\leqno(2.34)$$
\vskip.1in

{\bf Proof.} (ii) $\Rightarrow$ (i): 
Suppose $c$ is a curve in $S^3$ parametrized by
arclength and  the unit tangent vector {\bf t} of
$c$ is given by (2.26). Let $$g_{1}= f_{2}f'_{3} -
f_{3}f'_{2},\,\,\,\, g_{2}= f_{3}f'_{1} - f_{1}f'_{3},$$
$$ g_{3}= f_{1}f'_{2} - f_{2}f'_{1}.$$  

\noindent By Frenet-Serret formulas and  (2.33) we have 
$$({{g_{i}}\over {\kappa}})' = {2f'_{i}\over
{\kappa}}, \,\,\,\, i = 1, 2, 3.\leqno(2.35)$$

\noindent By using (2.34) and the identity ${\bf b} = {\bf
t}\times{\bf n}$, we may obtain
 $$a = {{g_{1}}\over {\kappa}}.\leqno(2.36)$$ 

\noindent Hence, we find
$2f'_{1}/{\kappa} = a' = 0$. Let $b$ denote the constant
 $f_{1}$.  Then, from (2.36) and 
${\bf b} = {\bf t}\times{\bf n}$, we may find
$${\bf b} = a{\tilde X}_{1} - ({{bf'_{3}}\over
{\kappa}}){\tilde X}_{2} + ({{bf'_{2}}\over
{\kappa}}){\tilde X}_{3},\leqno(2.37)$$
 $$\kappa^{2} = (f'_{2})^{2} + (f'_{3})^{2},\leqno(2.38)$$
$${\bf n} = ({{f'_{2}}\over {\kappa}}){\tilde X}_{2} +
({{f'_{3}}\over {\kappa}}){\tilde X}_{3}.\leqno(2.39)$$

\noindent Since $\|{\bf b}\| = 1, \|{\bf n}\| = 1,$ (2.37) and
(2.39) imply $a^{2} + b^{2} = 1.$ Thus, from $\|{\bf t}\| 
= 1$ and (2.26), we get $f_{2}^{2} + f_{3}^{3} = a^2$. So
we may put $$f_{2} = a\cos \theta, \,\,\, f_{3} = a\sin
\theta,\,\,\,\, \theta = \theta (s).\leqno(2.40)$$

\noindent Thus, by applying the definition of $g_{1}, g_{2}, g_{3}$,
we have
$$g_{1} = a \kappa,\,\,\, g_{2} = - bf'_{3},\,\,\, g_{3} =
bf'_{2}.\leqno(2.41)$$

\noindent By using (2.35), (2.38), (2.40), and $\tau \not= 0$, we get
$$\kappa = \,\mid a\theta'\mid\, \not= 0.\leqno(2.42)$$

From (2.26), (2.40), (2.41) and (2.42) we find
$\sin\theta (b\theta' + 2) = 0$. Since $\sin\theta(s)$
has only isolated zeros by (2.42), $b\theta' + 2 = 0$.
Thus, $b
\not= 0$. So $\theta = - {{2}\over {b}}s + s_{0}, s_{0} =$
const. Hence, by (2.37) and $\kappa > 0$, we get $ab <
0$.
\vskip.1in
$(i) \Rightarrow (ii)$ follows from straight-fordward
computation. 
\vskip.1in
{\bf Lemma 2.8.} {\it A helical cylinder $f(N) =
\{\gamma(t)\cdot c(s)\}$ in $S^3$ is a proper slant
surface with respect to the complex structure $J_{1}^-$
with the slant angle equal to $\cos^{-1}a$, where $a$ is
the constant given by (2.22) of Definition 2.1.}
\vskip.1in

{\bf Proof.} Let $\xi$ be the positive unit normal of
$f(N)$ in $S^3$ and ${\bf b}$  the binormal vector of
$c$ in $S^3$. Then we have
$$\xi(\gamma(t)\cdot c(s)) = L_{\gamma(t)^*}{\bf
b}(s).\leqno(2.43)$$  Lemma 2.8 now follows from Lemma
2.3 and Lemma 2.7.

\vskip.1in
{\bf Lemma 2.9.} {\it For any point $\, p_{0}\in S^3$ and
any oriented 2-plane $P_{0} \subset T_{p_{0}}S^{3}
\subset E^4$  which is proper slant with respect to the
complex structure $J_{1}^-$, there exist helical cylinders
in $S^3$ passing through $p_{0}$ and whose tangent planes
at $p_0$ are $P_0$.}

\vskip.1in
{\bf Proof.} Let $\xi$ be the positive unit normal of
$P_0$ in $T_{p_{0}}S^3$ and $\alpha$ the slant angle of
$P_0$ with respect to $J_{1}^-$. Put 
$$a = \cos \alpha
\,\,\,(\not= 0, \pm1), \,\,\,\, b =
\pm(1-a^{2})^{{1}\over {2}},\leqno(2.44)$$

\noindent  where $\pm$ was
chosen so that $ab < 0.$ Pick $s_{0} \in [0,2\pi)$ such
that  $$\cos s_{0} = -{{1}\over {b}}<\xi,{\tilde
X}_{2}(p_{0})>,\,\,\,\, \sin s_{0} = -{{1}\over
{b}}<\xi,{\tilde X}_{3}(p_{0})>.$$

\noindent We define $f_i$ by (2.29)-(2.52). Then they satisfy
(2.27) and we can choose a curve $c(s)$ satisfying the
conditions mentioned in Lemma 2.6.

Let $\gamma(t)$ be either a geodesic in $S^3$ safisfying
$$\gamma(0) = p_{0},\,\,\,\gamma'(0) \in
P_{0},\,\,\,\gamma'(0) \not= c'(0),$$ or a curve in $S^3$
satisfying this condition and also the condition  $\tau
\equiv 1$ (see Theorem 3 of [Sp1, p.35] for the existence
of such curves). Then we can verify that $\{\gamma(t)\cdot
c(s)\}$ is a desired surface. 
\vskip.1in
{\bf Proof of Theorem 2.5.} First, we note that the
isometry $\,\phi\,$ of $\, S^3$ has the following
properties: 
$$\phi(p\cdot q) = \phi(q)\cdot\phi(p),
\,\,\,\, {\rm for}\,\,\,\forall p, q \in
S^{3},\leqno(2.45)$$ 
\begin{align}&X \in {\mathcal X}(S^{3})
\,\,\,\,{\rm is \,\, left{\rm -}\, (repsectively,\,
right{\rm -})\,invariant}\tag{2.46}\\ &\Longleftrightarrow \phi_{*}X
\,\,\, \rm is\,\, right{\rm -} \, (respectively,\,
left{\rm -})\, invariant,\notag\end{align}
$$\tau_{\phi\circ c} = -\tau_{c} \,\,\,{\rm for\, a\,
curve}\, c\,\, {\rm in}\, S^{3},\leqno(2.47)$$ $${\bf
b}\,\, {\rm is\, the\, binormal\,  of\, a \, curve\,}\,
c\, {\rm in}\, S^{3} \Longleftrightarrow \leqno(2.48)$$
$$\hskip.9in -\phi_{*}{\bf b}\, {\rm is\, the\,
binormal\,of }\,\, \phi\circ c \,\, {\rm in}\,
S^{3},$$ 

\noindent where $\tau_c$ denotes the torsion of
the curve $\, c \,$ in $\, S^3$.

Let $\alpha$ be the slant angle of $f(N)$ with respect to
$J_1$. Since $f(N)$ is spherical, the normal curvature
$G^D$ of the slant immersion $f$ vanishes. Thus, by 
Theorem 3.2 in Section 3 of Chapter II, $N$ is a flat
surface in $S^3$. Therefore, $f(N)$  is locally a flat
translation surface $f(N) = \{c(s)\cdot\gamma(t)\}$ (cf.
[Sp1, pp.149-157]), where $c$ and $\gamma$ are curves in
$S^3$ parametrized by arclength satisfying one of the
following conditions: $$\tau_{c} \equiv +1 \,\,\,\,{\rm
and}\,\,\,\,\tau_{\gamma} \equiv -1,\leqno(i)$$
$$\tau_{c} \equiv +1\,\,\,\,{\rm and}\,\,\gamma\,\,{\rm
is\,\, a\,\, geodesic,}\leqno(ii)$$
$$c \,\,{\rm is\,\, a \,\, geodesic\,\, and}\,\,
\,\,\tau_{\gamma} \equiv -1,\leqno(ii')$$
$$c\,\,\,{\rm and}\,\,\,\gamma \,\,\,{\rm are\,\,
distinct\,\, geodesics.}\leqno(iii)$$
\vskip.1in

 {\bf Cases (i) and (ii):} Let {\bf
b} be the binormal of $c$. With a suitable choice of the
orientation, {\bf b} is the positive unit normal of
$f(N)$ in $S^3$. By Lemma 2.3, Lemma 2.7, (2.47), and
(2.48), $\phi\circ c$ is a helix in $S^3$ with $a$ and
$b$ in (2.29)-(2.32) determined by 
$$a = \cos\alpha,\,\,\,\,\, b = \pm\sin\alpha,\,\,\, ab <
0,\leqno(2.49)$$ where either $\tau_{\phi\circ\gamma}
\equiv +1$ or $\phi\circ\gamma$ is a geodesic. So, by
(2.45), $(\phi\circ f)(N)$ is a helical cylinder in $S^3$.

The converse is given by Lemma 2.8. Moreover, Lemma 2.9
guarantees the existence of such surfaces.

\vskip.1in
Next, we want to show that both cases (ii$'$) and (iii) 
cannot occur. Without loss of generality we can assume
$$c(0) = \gamma(0) = 1 \in S^{3},\leqno(2.50)$$
because Lemmas 3.1 and 3.3 imply that the slantness of a
surface in $S^3$ with respect to $J_1$ is
right-invariant, that is, if $f$ is $\alpha$-slant with
respect to $J_1$, so is $R_{q}\circ f$ for any $q \in
S^3$, and hence we can replace $f, c$ and $\gamma$ by
$R_{c(0)\cdot\gamma(0)}\circ f, R_{c(0)}\circ c$ and
$L_{c(0)}\circ R_{\gamma(0)^{-1}\cdot
c(0)^{-1}}\circ\gamma$, respectively, if necessary.
\vskip.1in 

{\bf Case (ii$'$):} Let {\bf b} be the binormal of
$\gamma$. We can choose the orientation so that
$$\xi(c(s)\cdot\gamma(t)) = L_{c(s)^{*}}{\bf
b}(t),\,\,\,\,\,\,{\rm for\,\,any\,\,} s,
t.\leqno(2.51)$$ So, by Lemmas 2.1 and 2.3 and
also(2.50), we have $$<L_{c(s)^{*}}{\bf b}(0),
R_{c(s)^{*}}X_{1}> = -\cos\alpha\,\,\,\,\,{\rm for}
\,\,{\rm any\,\,\,} s.\leqno(2.52)$$

 Put
$$c'(0) = (0,a_{1},a_{2},a_{3}),\,\,\, {\bf b}(0) =
(0,b_{1},b_{2},b_{3}) \in T_{1}S^{3} \subset
E^{4}.\leqno(2.53)$$ Then 
$$c(s) = (\cos s, a_{1}\sin s, a_{2}\sin s, a_{3}\sin
s).\leqno(2.54)$$

Putting $s = 0$ in (2.52), we find
$$b_{1} = -\cos\alpha \not= 0, \pm1,\leqno(2.55)$$

\noindent  since
$f(N)$ is properly slant. On the other hand, by (2.4),
(3.5), (2.45), and (2.46), we have
$$<L_{c(s)^{*}}{\bf b}(0), R_{c(s)^{*}}X_{1}> = b_{1}\cos
2s + (-a_{3}b_{2} + a_{2}b_{3})\sin 2s.\leqno(2.56)$$

\noindent So, from (2.52) and (2.56), we get $b_{1} = 0$ which
contradicts to (2.55). Consequently, case (ii$'$) cannot
occur.

\vskip.1in
{\bf Case (iii):} Let $\{c(s)\cdot\gamma(t)\}$ be defined
by using two distinct geodesics $c$ and $\gamma$ and
assume 
$$f : I_{1}\times I_{2} \rightarrow S^{3}\,\,
;\,\, (s,t) \mapsto c(s)\cdot\gamma(t)\leqno(2.57)$$ 

\noindent is
properly slant. Since the geodesics $c$ and $\gamma$ can be
extended for all $\, s, t \in {\bf R}$, we can extend
the immersion $x$ to a global mapping: 
$$y : {\bf R}^{2}
\rightarrow S^{3}\,\, ;\,\, (s,t) \mapsto
c(s)\cdot\gamma(t).\leqno(2.58)$$ 

Now, we claim that  $y$
is in fact an immersion and it is properly slant. To see
this, we recall  (2.50) and put
 $$c'(0) =
(0,a_{1},a_{2},a_{3}),\,\,\,\gamma'(0) =
(0,b_{1},b_{2},b_{3}) \in T_{1}S^{3}.\leqno(2.59)$$ 

Let ${\tilde X}, {\tilde Y}$ be the
vector fields along $y({\bf R}^{2})$ defined by 
$${\tilde
X}(s,t) = R_{\gamma(t)^{*}}c'(s),\,\,\,{\tilde Y}(s,t) =
L_{c(s)^{*}}\gamma'(t).\leqno(2.60)$$ 

\noindent Then it follows from
(2.4) and (2.5) that
$$<{\tilde X}(s,t),{\tilde Y}(s,t)>
\,\, {\rm is\,\, a\,\, polynomial \,\, of}\,\, \sin s, \,
\cos s,\, \sin t\,\, {\rm and}\, \cos t.\leqno(2.61)$$

On the other hand, since the $s$-curves and the $t$-curves
on $f(I_{1}\times I_{2})$ intersect at a constant angle
(cf. [Sp1, p.157]), we have
 $$<{\tilde X}(s,t),{\tilde
Y}(s,t)> = const. \not= \pm1,\,\,\,\, {\rm for \,\,any\,}
\,(s,t) \in I_{1}\times I_{2}.\leqno(2.62)$$ 

\noindent From (2.61), we see that
(2.62) holds for all $(s,t) \in {\bf R}^2$ and hence $y$
is an immersion. Since 
$$\xi(c(s)\cdot\gamma(t)) =
\|{\tilde X}(s,t)\times {\tilde Y}(s,t)\|^{-1}({\tilde
X}(s,t)\times{\tilde Y}(s,t)),$$

\noindent  where $\times$ denotes
the usually vector product in
$\, T_{c(s)\cdot\gamma(t)}S^{3} \,$ determined by the
metric and the orientation, so, by (2.4), (2.5) and (2.6),
we know that
$<\xi(c(s)\cdot\gamma(t)),J_{0}\eta(c(s)\cdot\gamma(t))>$
is a polynomial  of $\,\sin s, \cos s, \sin t\,$ and $\,
\cos t.$ By Lemma 2.3 we conclude that this polynomial is
a constant on $I_{1}\times I_2$ and hence $y$ is a proper
slant immersion defined globally on ${\bf R}^2$. Now, by
the double periodicity, $y$ induces a proper slant
immersion: 
$${\tilde y} : T^{2} = ({\bf R}/2\pi {\bf
Z})\times({\bf R}/2\pi {\bf Z}) \rightarrow {\bf C}^{2} =
(E^{4}, J_{0})$$

 \noindent of a torus into ${\bf
C}^2$, which contradicts to Theorem 1.5 of Chapter V.
Consequently, case (iii) cannot occur. 

This completes the
proof of the theorem.

\vfill\eject
\noindent{}  
\vskip.3in
\noindent \S 3.  SLANT SURFACES WITH
$rk(\nu)<2.$
\vskip.2in
For an immersion $f:N \rightarrow {\bf C}^{m}$, the Gauss
map $\nu$ of the immersion $f$ is given by 
$$\nu:N\rightarrow G(n,2m)\equiv
D_{1}(n,2m) \subset S^{K-1} \subset
\wedge^{n}(E^{2m}),\leqno(3.1)$$
$$ \nu (p) =
e_{1}(p)\wedge\ldots\wedge {e_{\ell}} (p),\,\,\, p
\in N,$$

\noindent where $n = dim\, N,\,K={{2m}\choose {n}},\,
D_{1}(n,m)\,$ is the set of all unit decomposable
$\,n$-vectors in $\wedge^{n}E^{2m}$, identified with
the real Grassmannian $G(n,2m)$ in a natural way, and $
S^{K-1}$ is the unit hypersphere of $\wedge^{n}(E^{2m})$
centered at the origin, and $\{e_{1},\ldots,e_{2m}\}$
is a local adapted orthonormal tangent frame along $f(N)$.

The main purpose of this section is to prove the following
classification theorem [CT2]. 

\vskip.1in
{\bf Theorem 3.1.} {\it If $f : N \rightarrow {\bf
C}^{2}  =(E^{4},J_{1})$ is a  slant
immersion    such that the
rank of its Gauss map is less than 2, then the image
$f(N)$ of $f$ is a union of some flat ruled surfaces in
$E^4$. Therefore, locally, $f(N)$ is a cylinder, a cone
or a tangential developable surface in ${\bf C}^2$.
Furthermore,}

(i) {\it A cylinder in {\bf C}$^2$ is a slant
surface if and only if it is of the form $\{ c(s)+te\},$
where $e$ is a fixed unit vector and $c(s)$ is a
(Euclidean) generalized helix  with axis  $J_{1}e$
contained in a hyperplane  of $E^4$ and with $e$ as
its hyperplane normal,}

(ii) {\it A cone in {\bf C}$^2$ is a  slant
surface if and only if, up to translations,  it is of the
form $\{tc(s)\},$ where $(\phi\circ c)(s)$ is a
generalized helix in $S^3$ with axis ${\tilde X}_1$
(cf. Definition 2.1), and}

(iii) {\it A tangential developable surface $\{
c(s)+(t-s)c'(s)\}$ in {\bf C}$^2$ is a  slant
surface if and only if, up to ridid motions, $(\phi\circ
c')(s)$ is a generalized helix in $S^3$ with axis  
${\tilde X}_1$.}

 \vskip.1in
 
As before, let $*$ denote the Hodge star operator $* :
\wedge^{2} E^{4} \rightarrow \wedge^{2} E^4$ induced from
the natural orientation and the canonical inner product of
$E^4$. Denote by $\wedge_{+}^{2}E^4$ and
$\wedge_{-}^{2}E^4$ the eigenspaces of $*$ with
eigenvalues $+1$ and $-1$, respectively. 
And denote by $S_{+}^2$ and $S_{-}^2$ the 2-spheres in
$\wedge_{+}^{2}E^4$ and  $\wedge_{-}^{2}E^4$
centered at the origin  with radius $1/{\sqrt{2}}$,
respectively. Then we have $D_{1}(2,4) = S_{+}^{2}\times
S_{-}^2$. Let 
 $$\pi_{+} : D_{1}(2,4) \rightarrow
S_{+}^{2},\,\,\,\,\,\, \pi_{-} : D_{1}(2,4) \rightarrow
S_{-}^{2}\leqno(3.2)$$ 

\noindent denote the natural projections. We
define as before the two maps $\nu_+$ and $\nu_-$ given
respectively by
 $$\nu_{+} = \pi_{+}\circ\nu \,\,\,\,\,{\rm
and}\,\,\,\,\, \nu_{-} =\pi_{-}\circ\nu.\leqno(3.3)$$ 

Suppose that the  slant immersion $f:N \rightarrow
{\bf C}^{2} =  (E^{4},J_{1})$ satisfying the
condition $rank(\nu) < 2.$ Then we also have
$rank(\nu_{\pm}) < 2$. Hence, by Lemma 4.2 of Chapter III,
$f(N)$ is a flat surface in $E^4$.

Furthermore, we
have [CT2]

\vskip.1in
{\bf Lemma 3.2.} {\it If $f$ is a  slant immersion
with $rank(\nu) < 2$, then $f(N)$ is a union of flat ruled
surfaces in $E^4$.}
\vskip.1in 

{\bf Proof.} Since the normal curvature $R^{D}=
0$, identically, we can choose local orthonormal frame
$\{e_{1},e_{2}\}$ such that with respect to it the second
fundamental form $\{h_{ij}^{r}\}$ is simultaneously
diagonalized, that is, we have
 $$(h_{ij}^{3}) =
\begin{pmatrix}b&0\\ 0&c\end{pmatrix}, \,\,\,\, (h_{ij}^{4}) =
\begin{pmatrix}d&0\\ 0&e\end{pmatrix}.\leqno(3.4)$$ 

Put
$$N_{1} = \{\, p\in N\,\mid\, H(p) \not= 0 \,\},\,\,N_{0} =
{\rm Interior \, of}\, (N-N_{1}).\leqno(3.5)$$ 

\noindent Then $$N =
N_{0}\cup \partial N_{1}\cup N_{1},$$ where $H$ is the
mean curvature vector of $N$ in ${\bf C}^2$.

 Since  $f(N)$ is flat and $H = 0$ on $N_{0}$,  $f(N_{0})$
is a union of portions of 2-planes in $E^4$ with the same
slant angle.

On $N_{1}$, we put $e_{3}=H/\|H\|$. Since $rank(\nu ) <
2$, we have $bc = 0$ and $d = e=0$. We may choose
$\{e_{1},e_{2}\}$ such that $b \not= 0, c=d=e=0$ on $N_1$.
From these we may prove that the integral curves of
$f_{*}e_2$ are open portions of straight lines and
therefore $f(N_{1})$ is a union of flat ruled surfaces.
Consequently, $f(N)$ is a union of flat ruled surfaces
possibly guled along $\partial N_1$. This proves the lemma.

\vskip.1in
For the local classifcation of flat ruled surfaces in
$E^4$, see
 [Sp1].

Now we give the proof of Theorem 3.1.

\vskip.1in
{\bf Proof of Theorem 3.1.} The first part of the theorem
is giving in Lemma 3.2. Now, we prove the remaining part
of the theorem.
\vskip.1in

{\bf Case (i):} If $f(N)$ is a slant cylinder, then we
may assume that $f(N)$ is of the form:
$$f(N) = \{c(s) + te \},\leqno(3.6)$$

\noindent where $e$ is a fixed unit vector in $E^4$ and $c(s)$ is a
curve parametrized by arclength which lies in the
orthogonal complement (up to sign), say $W \in G(3,4)$, of
$e$. Since $\{c'(s),e\}$ is a positive orthonormal basis
of $T_{c(s)+te}N,$ we obtain $ \,\cos\alpha =
<c'(s),-J_{1}e>$ by (1.4). Hence, $c(s)$ is a generalized
helix lies in the hyperplane $W (\equiv E^{3})$ whose
tangents make a constant angle $\alpha$ with $-J_{1}e\in
W$. \vskip.1in

{\bf Case (ii):} If $f(N)$ is a slant cone in $E^4$, then,
without loss of generality, we may assume that the vertex
of the cone is the origin of $E^4$. So we can write
$$f(N) = \{tc(s)\},\leqno(3.7)$$

\noindent where $c(s)$ is a curve in $S^3$ parametrized by
arclength. Since $\{c'(s), \eta(c(s))\}$ is a positive
orthonormal basis of $T_{tc(s)}N$, $\,\cos\alpha =
<c'(s),-J_{1}\eta(c(s))>$ for all s. Thus, by Lemmas 2.1
and 2.3, we conclude that $(\phi\circ c)(s)$ is a
generalized helix in $S^3$ with axis   ${\tilde
X}_1$ (cf. Definition 2.1).
\vskip.1in
{\bf Case (iii):} If $f(N)$ is a slant tangential
developable surface in $E^4$, the surface has the form:
$$f(N) = \{c(s)+(t-s)c'(s)\},\leqno(3.8)$$

\noindent where $c(s)$ is a curve parametrized by arclength. We put
$${\bf v}_{1}(s) = c'(s),\,\,\,\kappa_{1}(s) = \|{\bf
v}'_{1}(s)\|,\,\,\, {\bf v}_{2}(s) =
({{1}\over {\kappa_{1}(s)}}){\bf
v}'_{1}(s).\leqno(3.9)$$

\noindent  Note that $\kappa_{1} \not= 0$,
since $c(s)$ generates a tangential developable surface.
$\{{\bf v}_{2}(s) ,{\bf v}_{1}(s)\}$ forms a positive
orthonormal basis of $T_{c(s)+(t-s)c'(s)}N$, and so we
have 
$$\cos\alpha = <{\bf v}'_{1}(s)/\|{\bf v}'_{1}(s)\|,
-J_{1}{\bf v}_{1}(s)>$$ 

\noindent for all $s$. If we consider ${\bf
v}_{1}(s)$ as a curve in $S^3$, then (3.9) means that
$$\cos\alpha = <{\bf t}(s), -J_{1}\eta({\bf v}_{1}(s))>$$

\noindent where ${\bf t}(s)$ is the unit tangent of ${\bf
v}_{1}(s)$. So, as in Case (ii), $(\phi\circ {\bf
v}_{1})(s)$ is a generalized helix in $S^3$ with axis
  ${\tilde X}_1$. 

It is easy to verify that in each of the cases 
(i)-(iii), the converse is also true. 

This completes the proof of the theorem.

\vfill\eject
\noindent{}  

\vskip.3in
\noindent \S 4. SLANT SURFACES WITH CODIMENSION ONE
 \vskip.2in
In this section we want to classify  slant surfaces which
are contained in a hyperplane $W$ of  $E^4$. 

\vskip.1in
{\bf Lemma 4.1.} {\it Let $f:N \rightarrow {\bf C}^{2}
= (E^{4},J_{1})$ be a  slant immersion of an
oriented surface $N$ into ${\bf C}^2$. If $N$ is contained
in some hyperplane $W\in G(3,4)$, then }

(1) {\it $rank(\nu) < 2$  and}

(2) {\it The immersion $f$ is doubly slant with the same
slant angle.}

\vskip.1in
{\bf Proof.} We choose
a positive orthonormal $J_1$-basis
$\{e_{1},e_{2},e_{3},e_{4}\}$ such that $e_{1}, e_{2} \in
W\cap J_{1}W, \,\, e_{4}=J_{1}e_{3} = \eta_W$, where
$\eta_W$ is the positive unit normal vector of the
hyperplane $W$ in $E^4$. We put
$$G_{W} = G(2,4)\cap \wedge^{2}W \subset
\wedge^{2}E^{4}.\leqno(4.1)$$ 

\noindent Then $G_W$ is the unit
2-sphere in the 3-dimensional Euclidean space
$\wedge^{2}W$. 

For $\alpha \in [0,\pi]$ we put
$$G_{W,\alpha} = G_{J_{1},\alpha}\cap G_{W},\leqno(4.2)$$

\noindent where $G_{J_{1},\alpha}$ is the set of all 2-planes in
$E^4$ with slant angle $\alpha$ with respect to $J_1$.
We recall that a 2-plane $V$ was identified with a unit
decomposable 2-vector $e_{1}\wedge e_{2}$ in
$\wedge^{2}E^4$ with $\{e_{1}, e_{2}\}$ as a
positively oriented orthonormal basis  of $V$.  From the
proof of Lemma 2.2, we see that $G_{W,\alpha}$ is the
circle on $G_{W} = S^{2} \subset \wedge^{2}W$ defined by

$$G_{W,\alpha} = \{V \in G_{W}\,\mid\, <V,e_{1}\wedge
e_{2}> = \cos\alpha\}.$$

For each $J \in {\mathcal J}$, we denote by $\zeta_J$ the
2-vector which is the metrical dual of $-\Omega_J$ as
defined in Section 2. Let $\zeta : {\mathcal J} \rightarrow
\wedge^{2}E^4$ be the mapping defined by $\zeta(J) =
\zeta_J$. Then $\zeta$ gives rise to two bijections (cf.
Lemma 2.1 of this chapter):
 $$\zeta^{+} : {\mathcal J}^{+} \rightarrow
S_{+}^{2}\,\,\,\,\,{\rm and}\,\,\,\,\,\zeta^{-} ; {\mathcal
J}^{-} \rightarrow S_{-}^{2}.$$

For each oriented 2-plane $V \in G(2,4)$ we define two
complex structures $J_{V}^{+} \in {\mathcal J}^+$ and
$J_{V}^{-} \in {\mathcal J}^-$ by 

$$J_{V}^{+} =
(\zeta^{+})^{-1}(\pi_{+}(V))\,\,\,\,\,{\rm
and}\,\,\,\,J_{V}^{-}=(\zeta^{-})^{-1}(\pi_{-}(V)).$$

Let ${\hat J} = J_{e_{1}\wedge
e_{2}}^{-}$. Then we have
$$\pi_{+}(G_{W,\alpha}) = S_{{\hat J},\alpha}^{+}
\subset S_{+}^{2},\,\,\,\,\pi_{-}(G_{W,\alpha}) =
S_{{\hat J},\alpha}^{-} \subset S_{-}^{2},\leqno(4.3)$$

\noindent where $S_{J,\alpha}^{\pm}$ are the circles (possibly
singletons) on $S_{\pm}^2$, respectively, consisting of all
2-vectors which make constant angle $\alpha$ with
$\zeta_{J}$. If $f$ is $\alpha$-slant with
respect to ${\hat J}$ and $f(N) \subset W$, then $\nu(M)
\subset G_{W,\alpha}$. Therefore, $rank(\nu) < 2$ and, by
(4.3), $f$ is  $\alpha$-slant with respect to ${\hat J}$.

This proves the lemma.
\vskip.1in
We note here that if we identify $\wedge^{2}W$ with the
Euclidean
3-space $E^{3} \equiv W$ (where $W$ is  spanned by
$\{e_{1},e_{2},e_{3}\}$) via the isometry $X\wedge Y
\rightarrow X\times Y$, then $\nu:M \rightarrow G_{W}
\subset \wedge^{2}W$ is nothing but the classical
Gauss map $g:M \rightarrow S^{2} \subset E^3$.

 Since
$e_{1}\times e_{2} = e_{3} = -J_{1}\eta_{W},$ $f$ is
$\alpha$-slant if and only if 
$$g(M) \subset S_{\alpha}^{1} = \{Z \in S^{2}\,\mid\,
<Z,-J_{1}\eta_{W}> = \cos\alpha\} \subset S^{2}\subset
W.\leqno(4.4)$$
\vskip.1in

Now we give the following classification theorem [CT2].
\vskip.1in
{\bf Theorem 4.2.} {\it Let $f:N
\rightarrow {\bf C}{^2} = (E^{4},J_{1})$ be a proper slant
immersion of an oriented surface $M$ into {\bf C}$^2$. If
$f(N)$ is contained in a hyperplane $W$ of $E^4$, then $f$
is a doubly slant immersion and
$f(N)$ is a union of some flat ruled surfaces in $W$.
 Therefore, locally, $f(N)$ is a cylinder, a cone
or a tangential developable surface in $W$.
Furthermore,}

 (i) {\it A cylinder in $W$ is a proper slant surface
with respect to the complex structure $J_1$ on $E^4$ if and
only if it is a portion of a 2-plane.}

(ii) {\it A cone in $W$ is a proper slant surface
with respect to the complex structure $J_1$ on $E^4$ if and
only if it is a circular cone.}

(iii) {\it A tangential developable surface in $W$ is
a proper slant surface with respect to the complex
structure $J_1$ on $E^4$ if and only if it is a
tangential developable surface obtained from a
generalized helix in $W$.}

\vskip.1in

{\bf Proof.} Assume  $f:N \rightarrow 
{\bf C}^{2} = (E^{4},J_{1})$ is a proper slant immersion of
an oriented surface $N$ such that  $f(N)$ is contained in
some hyperplane $W\in G(3,4)$. The first part of Theorem
4.2 is given by Lemma 4.1. For the remaining part it
suffices to check the three cases of Theorem 3.1. 

 Suppose $f$ is properly slant
with slant angle $\alpha$. Denote by $\xi$ the local
unit normal of $f(N)$ in $W$. We put
$$e_{1} = t\xi / \|t\xi\|,\,\,\,\,e_{2} =
(\sec\alpha)Pe_{1},\,\,\,e_{3} =
(\csc\alpha)Fe_{1},\,\,\,e_{4} =
(\csc\alpha)Fe_{2},\leqno(4.5)$$ 

\noindent where $PX$ and $FX$
denote the tangential and the normal components of
$J_{1}X$, respectively, and $t\xi$ is the tangential
component of $J_{1}\xi$.  Then $\{e_{1},\cdots,e_{4}\}$
is an adapted orthonormal frame along
$f(N)$ and it satisfies
 $$e_{3} = \,{\rm unit\,\,
normal\,\, of}\, f(N)\,\,{\rm in}\,\, W, \,\,\,e_{4} \in
W^{\bot},\leqno(4.6)$$
$$te_{3} =
-(\sin\alpha)e_{1},\,\,\,te_{4}=-(\sin\alpha)e_{2},
\leqno(4.7)$$
$$fe_{3}=-(\cos\alpha)e_{4},\,\,\,fe_{4}=(\cos\alpha)e_{3},
$$ 

\noindent where $fe_3$ is the normal component of $J_{1}e_3$.
Since $e_4$ is a constant vector in $E^4$, Corollary 3.6 of
Chapter II implies that the second fundamental form
$(h_{ij}^{r})$ is of the following form:
 $$(h_{ij}^{3}) = \begin{pmatrix}b&0\\
0&0\end{pmatrix},\,\,\,\, (h_{ij}^{4}) = 0,\leqno(4.8)$$

\noindent which shows that the orthonormal frame
$\{e_{1},\cdots,e_{4}\}$ coincides with that chosen in the
proof of Lemma 3.2 (up to orientations). Since $J_{1}e_4$
is also a constant vector in $E^4$, from (4.7), we have 
$$
-\sin\alpha \nabla_{X}e_{2} - \cos\alpha A_{e_{3}}X =
0,\,\,\,{\rm for}\,\, X\in TM.\leqno(4.9)$$

\noindent Hence we get 
$$\omega_{2}^{1}(e_{1}) = -
b\cot\alpha,\,\,\,e_{2}b = b^{2}\cot\alpha.\leqno(4.10)$$

\vskip.1in
{\bf Case (i):} In this case, the curve $c(s)$ of (3.6)
lies in a 2-plane $W' = \{e\}^{\bot}\cap W \subset W$
which is perpendicular to $e$. So, $f(N)$ is totally real
with respect to  the complex structures $\pm J_{W'}^{\pm}$
defined above. If $\nu_{+}(N)$ is not
a singleton, then $J_1$ is one of the complex
structures $\pm J_{W'}^{\pm}$ according to Propostion 3.1
and Theorem 3.2 of Chapter III. Hence we get $\alpha =
\pi/2$, which contradicts to the assumption. So,
$\nu_{+}(N) $ is a singleton and hence $f(N)$ is minimal
(cf. Theorem 3.2 of Chapter III). Thus, by (4.8), $f(N)$ is
an open portion of an $\alpha$-slant 2-plane.
\vskip.1in

In Cases (ii) and (iii), we may assume $N =N_1$

 {\bf Case (ii):}
In this case the curve $c(s)$ in (3.7) lies in the unit
2-sphere $S^{2} = S^{3}\cap W$. Choose
$\{e_{1},\cdots,e_{4}\}$ according to (4.5) and let ${\bf
t},{\bf n},{\bf b},\kappa$, and $\tau$ be the unit
tangent vector, unit principal normal vector, the unit
binormal vector, the curvature, and the torsion of $c(s)$
in $W = E^3$, respectively. We want to show that $\tau
\equiv 0$.

Since
$$e_{1}(s,t) = {\bf t}(s) = {{1}\over
{t}}{{\partial}\over {\partial s}},\,\, e_{2}(s,t) =c(s) =
{{\partial}\over {\partial
t}},\leqno(4.11)$$
$$\,\,e_{3}(s,t) = e_{1}(s,t)\times
e_{2}(s,t),$$ 

\noindent where $\times$ denotes the vector product in $W$,
we have
 $$b = -({{\kappa}\over {t}})<{\bf
b},c>.\leqno(4.12)$$ 

\noindent From $\|c\| = 1$, we get 
$$\kappa
<{\bf n},c> = -1.\leqno(4.13)$$ 

\noindent Differentiating (4.13)
with respect to $s$, we get 
$$\kappa^{2}\tau<{\bf b},c> =
\kappa'.\leqno(4.14)$$ 

\noindent From (4.12) we obtain
$$-t\tau \kappa b = \kappa'.\leqno(4.15)$$

\noindent Differentiating (4.15) with respect to $t$ and using
(4.10) and (4.15), we obtain
$$\kappa'(\tau\kappa\tan\alpha - \kappa') =
0.\leqno(4.16)$$ 

\noindent By (4.12), (4.14) and $<{\bf t},c>=0$, we
find
$$\kappa^{2}\tau c = -\kappa\tau{\bf n} + \kappa'{\bf
b}.\leqno(6.17)$$ 

\noindent Since $\|c\|$ = 1, we also get
$$\tau^{2}\kappa^{4}=\tau^{2}\kappa^{2} +
(\kappa')^{2}.\leqno(4.18)$$

If $\kappa'(s_{0})=0$ at a point $s=s_0$, then, by (4.15),
we have $\tau(s_{0})=0,$ since $b(s,t) \not= 0$ by
assumption and also  $\kappa(s) \not= 0$ because
$c(s)$ is spherical.

If $\kappa'(s_{0})\not= 0$, we choose a neighborhood $U$
of $s_0$ on which $\kappa'$ never vanishes. By (4.16),
(4.18) and $\kappa \not= 0$, we get
$$(\tau(s))^{2}\{(\kappa(s))^{2}-1-\tan^{2}\alpha\} = 0
\,\,\,\,\,{\rm for}\,\,\forall s \in U.\leqno(4.19)$$

If $\tau(s_{0})\not= 0$ in addition, we choose another
neighborhood $U'$ of $s_0$ contained in $U$ on which
$\tau$ never vanishes. Then, by (4.19), we get
$$(\kappa(s))^{2}-1-\tan^{2}\alpha = 0$$

\noindent for all s in $U'$. By
continuity we get $\kappa(s)=constant$ on $U'$ which
contradicts to $\kappa'(s) \not= 0$ on $U'$. So, again we
have $\tau(s_{0})=0.$ Therefore, $\tau \equiv 0,$ which
means that $c(s)$ is a circle on $S^2$ and thus $f(N)$ is a
circular cone. According to the remark after 
 Lemma 6.1, the axis of the cone is given by
$-J_{1}e_4$. 
\vskip.1in

{\bf Case (iii):} We  assume the surface
is given by  (3.8) and $\{e_{1},\ldots\,e_{4}\},$ $ {\bf
t},$ ${\bf n}, 
{\bf b}, \kappa$ and $\tau$ are given as in Case
(ii). We have 
$$e_{1}(s,t)={\bf n}(s)={{1}\over
{(t-s)\kappa}}{{\partial}\over {\partial s}},\,\,
e_{2}(s,t) = {\bf t}(s) = {{\partial}\over {\partial
t}},\leqno(4.20)$$
$$\,\,e_{3}(s,t)=e_{1}\times e_{2} = -{\bf
b}(s).$$ 

\noindent Hence
$${\tilde {\nabla}}_{e_{1}}e_{1} = -{{1}\over
{(t-s)}}e_{2}-{{\tau}\over
{(t-s)\kappa}}e_{3}.\leqno(4.21)$$ So, by (4.8), we find
$$b=-{{\tau}\over {(t-s)\kappa}}.\leqno(4.22)$$

By (4.20) and (4.22), we obtain
$$e_{2}b={{\tau}\over {\kappa (t-s)^{2}}}.$$
This formula together with (4.10) and (4.21) imply that
${{\tau}\over {\kappa}} =\tan\alpha$ is a constant. This
means that the curve $c(s)$ is a generalized helix in $W$.
 The axis of the helix is given by
$-J_{1}e_4$.

In each of the cases (i), (ii) and (iii), the converse is
easy to verify. For example, if $f(N)$ is a circular cone
with the axis vector $e$ in a 3-plane $W$ perpendicular to
a unit vector $\eta$ in $E^4$, then, by picking the complex
structure $J$ so that $J = J_{e \wedge \eta }^+$, 
$f(N)$ is  properly slant with respect to $J$.

This completes the proof of the theorem.

\vskip.1in
{\bf Remark 4.1.} In the classifications of slant surfaces
given in Sections 2, 3 and 4 of this chapter, we avoid the
messy argument of glueing.

\vfill\eject
\vskip1in
\centerline {CHAPTER V}
\vskip.2in
\centerline {\bf TOPOLOGY AND STABILITY OF SLANT
SUBMANIFOLDS}
\vskip.5in

\noindent \S 1. NON-COMPACTNESS OF PROPER SLANT
SUBMANIFOLDS. 
\vskip.2in
Let $E^{2m}=(R^{2m},<\,,\,>)$ and {\bf C}$^{m} =
(E^{2m},J_{0})$ be the Euclidean 2m-space and the complex
Euclidean m-space, respectively,  with the canonical inner
product $<\,,\,>$ and the canonical (almost) complex
structure $J_{0}$ given by
$$J_{0}(x_{1},\ldots,x_{m},y_{1},\ldots,y_{m})=(-y_{1},
\ldots,-y_{m},x_{1},\ldots,x_{m}). \leqno(1.1)$$

Denote by $\Omega_{0}$ the Kaehler form of {\bf C}$^m$,
that is,
$$\Omega_{0}(X,Y)=<X,J_{0}Y>, \hskip.2in X,Y \in E^{2m},
\hskip.2in \Omega_{0} \in
\wedge^{2}(E^{2m})^*.\leqno(1.2)$$

For an immersion $f:N \rightarrow {\bf C}^{m}$, the Gauss
map $\nu$ of the immersion $f$ is given by 
$$\nu:N\rightarrow G(n,2m)\equiv
D_{1}(n,2m) \subset S^{K-1} \subset
\wedge^{n}(E^{2m}),\leqno(1.3)$$
$$ \nu (p) =
e_{1}(p)\wedge\ldots\wedge {e_{n}} (p),\,\,\, p
\in N,$$
\vskip.1in
\noindent where $n = dim N,\,K={{2m}\choose {n}},\,
D_{1}(n,2m)\,$ is the set of all unit decomposable
$\,n$-vectors in $\wedge^{n}E^{2m}$, identified with
the real Grassmannian $G(n,2m)$ in a natural way, and $
S^{K-1}$ is the unit hypersphere of $\wedge^{n}(E^{2m})$
centered at the origin, and $\{e_{1},\ldots,e_{2m}\}$
is a local adapted orthonormal tangent frame along $f(N)$.

\vskip.1in
Before we give the main result of this section,
we give the following lemmas.
\vskip.1in {\bf Lemma 1.1.}
{\it For $X_{1},\ldots,X_{2k} \in E^{2m}\,(k < m)$, we have
$$(2k)!\,\Omega_{0}^{k}(X_{1}\wedge\ldots\wedge
X_{2k})=\leqno(1.4)$$
$$\sum_{\sigma\in S_{2k}} sign(\sigma)
\Omega_{0}(X_{\sigma(1)},X_{\sigma(2)})\cdots\Omega_{0}
(X_{\sigma(2k-1)},X_{\sigma(2k)}),$$
where $S_{2k}$ is the permutation group of order $2k$, sign
denotes the signature 

\noindent of permutations and} $\Omega_{0}^{k}
\in \wedge^{2k}(E^{2m})^{*} \equiv
(\wedge^{2k}E^{2m})^{*}.$
\vskip.1in
{\bf Proof.} Let $e_{1},\ldots,e_m$ be an orthonormal
frame of $E^{2m}$ with its dual coframe given by
$\omega^{1},\ldots,\omega^{2m}$. Let
$$\Omega_{0}=\sum_{A,B=1}^{2m}\,\,\varphi_{AB}\omega^{A}\wedge
\omega^{B}.$$
Then by direct computation we have
$$\Omega_{0}^{k}(X_{1},\ldots,X_{2k}) $$
$$={1\over {(2k)!}}\sum_{\sigma}sign
\,(\sigma)(\sum\varphi_{_{A_{1}A_{2}}}\omega^{A_{1}}
(X_{\sigma(1)})\omega^{A_{2}}(X_{\sigma(2)})\dots$$
$$\ldots
(\sum\varphi_{_{A_{2k-1}A_{2k}}}\omega^{A_{2k-1}}
(X_{\sigma(2k
-1)})\omega^{A_{2k}}(X_{\sigma(2k)})).$$

\noindent From these we obtain  (1.4).
\vskip.1in
 {\bf Lemma 1.2.}  {\it Let $V \in
G(n,2m)$ and $\pi_{V} : E^{2m} \rightarrow V$ be the
orthogonal projection. If $V$ is $\alpha$-slant, that is,
$V$ is slant with slant angle $\alpha \not= \pi/2$,
 in {\bf C}$^{m}
= (E^{2m},J_{0})$, then the
linear endomorphism $J_V$ of $V$ defined by $$J_{V}= (\sec
\alpha)(\pi_{V}\circ J_{0}{}_{\mid}{} _{V})\leqno(1.5)$$
is a complex structure compatible with the inner product
$<\,\,,\,\,>_{\,|\, V}$. In particular, $n$ is even.}
\vskip.1in
{\bf Proof.} Let $$P=\pi_{V}\circ (J_{\,| V}): V
\rightarrow V,\leqno(1.6)$$
$$P^{\perp}=J_{\,| V}-P : V \rightarrow
V^{\perp}\leqno(1.7)$$ and $$Q =P^{2} : V \rightarrow
V.\leqno(1.8)$$ Then $$J_{\, |
V}=P+P^{\perp}.\leqno(1.9)$$ 

\noindent  By simple computation and
using (1.7), we have
$$<QX,Y>=<X,QY>\leqno(1.10)$$ and
$$<PX,Y>=-<X,PY>\leqno(1.11)$$
for any $X,Y \in V$. Since $V$ is assumed to be
$\alpha$-slant, $$\angle\,(JX,V)=\angle\,(JX,PX) =
\alpha$$ for any nozero vector $X\in V$. Hence we have
$$||PX||=\cos\alpha \, ||X||\leqno(1.12)$$
for any nonzero vector $X\in V$. By (1.10), $Q$ is a
self-adjoint endomorphism. Since $J_{0}^{2}=-I$,
(1.6)-(1.9) imply that each eigenvalue of $Q$ is equal to
$-\cos^{2}\alpha$ which
 lies in $[-1,0)$. Therefore, by using (1.5), we may prove
that $J_{V}^{2}=-I$ and
$$||J_{V}X||^{2}=\sec^{2}\alpha\,||PX||^{2}=||X||^{2}$$ for
any $X\in V$. This proves the lemma.
 \vskip.1in
Let ${\hat \zeta}_{0}$ be the metrical dual of
$(-\Omega_{0})^k$ with respect to the inner product
$<\,,\,>$ naturally defined on $\wedge^{2k}E^{2m}$, 
that is,
 $$<{\hat \zeta}_{0},\eta>
=(-1)^{k}\Omega_{0}^{k}(\eta)\hskip.4in {\rm
for\,\,any}\,\,\,\,  \eta \in
\wedge^{2k}E^{2m},\leqno(1.13)$$ then we have the following

\vskip.1in
{\bf Lemma 1.3.} {\it Let $\, V \in G(2k,2m).$ If $\,V$ is
$\alpha$-slant in {\bf C}$^m$ with ${\alpha \not= \pi
/2}$, then $$<{\hat \zeta}_{0},V> =\mu_{k}\cos^{k}\alpha,
\leqno(1.14)$$ where $\mu_k$ is a nonzero constant
 depending only on k.}
\vskip.1in

{\bf Proof.} Let $J_V$ be the complex structure on $V$
defined by Lemma 1.2. For a unit vector $X\in V$, we put
$Y=J_{V}X \in V$. Then we have
 $$\Omega_{0}(X,J_{V}X) = <-J_{V}Y,J_{0}Y> =
-\cos\alpha.\leqno(1.15)$$

If $X, Y \in V$ and $Z$ is perpendicular to $J_{V}X$, then
$$\Omega_{0}(X,Z)=\cos\alpha <X,J_{V}Z>=0.\leqno(1.16)$$

Therefore, if we choose an orthonormal $J_V$-basis
$\{e_{1},\ldots,e_{2k}\}$ on $V$, that is,
$$e_{2i}=J_{V}e_{2i-1},\,\,\,\,\,
i=1,\ldots,k,\leqno(1.17)$$ and $$V=e_{1}\wedge\ldots\wedge
e_{2k},\leqno(1.18)$$ via the natural identification of
$G(2k,2m)$ with $D_{1}(2k,2m)$, then we have
$$\Omega_{0}(e_{a},e_{b})=-\delta_{a^*}{}_b\cos\alpha \,\,
\,\,\,\, {\rm for} \,\,\,a<b,\leqno(1.19)$$ where
$$(2i)^{*}=2i-1,\,\,\,\, (2i-1)^{*}=2i \,\, \,\,\,\,{\rm
for}\,\,\, i=1,\ldots,k.\leqno(1.20)$$

By (1.18), Lemma 1.1, and (1.19), we find
$$(2k)!\,\Omega_{0}^{k}(V)=(2k)!\,\Omega_{0}(e_{1}
\wedge\ldots\wedge e_{2k})$$
 $$=\sum_{\sigma\in S_{2k}}
sign(\sigma)\,\,\Omega_{0}(e_{\sigma(1)},e_{\sigma(2)})
\cdots\Omega_{0}(e_{\sigma(2k-1)},e_{\sigma(2k)})$$
$$=\sum_{a_{1},\ldots ,a_{2k}=1}^{2k} \delta_{a_{1}\cdots
a_{2k}}^{12\cdots (2k)}
\Omega_{0}(e_{a_{1}},e_{a_{2}})\cdots
\Omega_{0}(e_{a_{2k-1}},e_{a_{2k}})$$ 
$$=\sum_{a_{1},\ldots\,a_{k}=1}^{2k}
\delta_{a_{1}a_{1}^{*}\cdots
{a_{k}a_{k}^{*}}}^{12\cdots\cdots
(2k)}\Omega_{0}(e_{a_{1}},e_{a_{1}^{*}})\cdots
\Omega_{0}(e_{a_{k}},e_{a_{k}^{*}})\leqno(1.21)$$
$$=2^{k}\sum_{a_{1}<a_{1}^{*}} \cdots
\sum_{a_{k}<a_{k}^{*}} \delta_{{a_{1}a_{1}^{*}}\cdots
 a_{k}a_{k}^{*}}^{12\cdots\cdots
(2k)}\Omega_{0}(e_{a_1},e_{a_{1}^{*}})\cdots\Omega_{0}
(e_{a_k},e_{a_k}^{*})$$
$$=2^{k}(-\cos\alpha)^{k}\sum_{a_{1}<a_{1}^{*}} \cdots
\sum_{a_{k}<a_{k}^{*}}\delta_{a_{1}a_{1}^{*}\cdots
a_{k}a_{k}^{*}} ^{12\cdots\cdots\cdot (2k)}$$
$$=2^{k}(-\cos\alpha)^{k}\, k!.$$ 

\noindent Hence, by  (1.13) and (1.21), we obtain (1.14) with
$\mu_{k}=2^{k}k!/(2k)!\,.$ \vskip.1in
{\bf Lemma 1.4.} {\it Let $f : N \rightarrow E^m$ be an
isometric immersion of an $n$-dimensional compact oriented
manifold $N$ into $E^m$. Then the Gauss map $\nu : N
\rightarrow \wedge^{n}(E^{m})$ is mass-symmetric in
$S^{K-1}$, $K={m\choose n}$, that is, the center of gravity
of $\nu$ coincides with the center of the hypersphere
$S^{K-1}$ in $\wedge^{n}(E^{m})$.}
\vskip.1in
{\bf Proof.} Let $e_{1},\ldots,e_n$ be an oriented
orthonormal local frame of $TN$ with its dual coframe
given by $\omega^{1},\ldots,\omega^n$. Then we have
$$dx=e_{1}\omega^{1}+\ldots+e_{n}\omega^{n}.\leqno(1.22)$$

By direct computation we have
$$dx\wedge\ldots\wedge dx = n!\,(e_{1}\wedge\ldots\wedge
e_{n})\,\omega^{1}\wedge\ldots\wedge \omega^{n}=n!\,\nu
(*1),$$ where $dx$ on the left-hand-side is repeated
$n$-times. Therefore, by applying the divergence theorem,
we have $$n!\int_N \nu*1=\int_N dx\wedge\ldots\wedge dx =
\int_N d(x \wedge dx\wedge\ldots\wedge dx) = 0.$$

This proves that the center of gravity of $\nu $ is the
origin of $\wedge^{n}(E^{m})$, that is, the Gauss map is
mass-symmetric in $S^{K-1}.$
 \vskip.1in
Now we give the following  [CT2]
\vskip.1in
{\bf Theorem 1.5.} {\it Let $f: N \rightarrow {\bf C}^m$ be
a  slant immersion of an $n$-dimensional differentiable
manifold $N$ into the complex Euclidean $m$-space {\bf
C}$^m$. If $N$ is compact, then $f$ is totally real.}

\vskip.1in
{\bf Proof.} Without loss of generality we may
assume that $N$ is oriented because otherwise
 we may simply replace $N$ by its two-fold covering. Assume
$f$ is $\alpha$-slant with $\alpha\not= \pi/2.$ Then, by
Lemma 1.2, $\,n$ is even. Put $n = 2k$. Since $N$ is
compact, Lemma 1.4 implies that the Gauss map $\nu$
is mass-symmetric in $S^{K-1}, K={2m\choose {2k}}$.
Therefore $$\int_{p\in N}\,<\nu(p),\zeta>*1
=0\leqno(1.23)$$ for any fixed $2k$-vector $\zeta \in
\wedge^{2k}(E^{2m})$, where $*1$ is the volume element
of $N$ with respect to the metric induced from the
immersion $f$. Let $\zeta={\hat \zeta}_{0}$, where ${\hat
\zeta}_0$ is defined by (1.13). Then Lemma 1.3 and (1.23)
imply $$\mu_{k}\, vol(N)\cos^{k}\alpha = 0.\leqno(1.24)$$
But this contradicts to the assumption $\cos\alpha \not=
0$. Hence $\alpha = \pi/2$ and $f$ is a totally real
immersion. \vskip.1in

We recall that a submanifold $N$ of an almost complex
manifold $(M,J)$ is said to be {\it purely real\/} [C2] if
every eigenvalue of $Q=P^2$ lies in $(-1,0]$. In fact, by
using a method similar to the proof of Theorem 1.5, we
may prove the following
\vskip.1in
{\bf Theorem 1.6.} {\it Let $f: N \rightarrow {\bf C}^m$
be a purely real immersion from an $n$-dimensional 
differentiable manifold $N$ into ${\bf C}^m$. If $N$ is
compact, then $f$ is totally real.}

 \vskip.1in {\bf Definition 1.1.} An
almost Hermitian manifold $(M,g,J)$ is called an {\it exact
sympletic manifold\/} if the sympletic form (or,
equivalently, the Kaehler form) ${\tilde \Omega}={\tilde
\Omega}_{J}$ of $(M,g,J)$ is an exact 2-form.  \vskip.1in 
For compact submanifolds in an exact sympletic manifold,
we have the following \vskip.1in
 {\bf Theorem 1.7.} {\it Every compact slant
submanifold $N$ in an exact sympletic manifold $(M,g,J)$
is totally real.} \vskip.1in
{\bf Proof.} Let $(M,g,J)$ be an exact sympletic manifold
with sympletic form ${\tilde \Omega}$. Then there exists a
1-form ${\tilde\varphi}$ on $M$ such that ${\tilde
\Omega}=d{\tilde\varphi}$.  Let $f: N\rightarrow M$ be
the immersion and we put 
$$\Omega =
f^{*}{\tilde\Omega},\,\,\,\varphi=f^{*}{\tilde\varphi}.
\leqno(1.25)$$
Then we have
$$\Omega=f^{*}{\tilde\Omega} =f^{*}d{\tilde\varphi} =
d(f^{*}\varphi)=d\varphi.\leqno(1.26)$$
If $N$ is either a proper slant submanifold or a complex
submanifold of $M$, then, by using Lemmas 1.1 and  1.2, we
may prove  that there is a nonzero constant $C$ (which
depends only on the dimesnion of $N$ and the slant angle)
such that $\Omega^{k}= C(*1)$ where $dim\,N=2k$ and $*1$
is the volume element of the slant submanifold $N$. From
(1.26) we know that $\Omega^k$ is exact. Therefore, by the
Stokes theorem, we have
 $$ vol(N) = \int_N (*1) =
C^{-1}\,\int_N \Omega^{k} = 0.$$

\noindent  This is a contradiction.

\vfill\eject

\noindent{}  

\vskip.3in
\noindent \S 2. TOPOLOGY OF SLANT SURFACES.
\vskip.2in
Let $N$ be an $n$-dimensional proper slant submanifold
with slant angle $\alpha$ in a Kaehlerian manifold $M$ of
complex dimension $m$. Then $n$ is even, say $n=2k$. Let
$e_1$ be a unit tangent vector of $N$. We put
$$e_{2} = (\sec\alpha)Pe_{1},\,\,\,\,e_{1^{*}}:=e_{n+1} =
(\csc\alpha)Fe_{1},\leqno(2.1)$$
$$e_{2^{*}}:=e_{n+2} = (\csc\alpha)Fe_{2}.$$

\noindent If $n>1$, then, by induction, for each $\, \ell =
1,\ldots,k-1,$ we  may choose a unit tangent vector
$e_{2\ell +1}$  of $N$ such that  $e_{2\ell + 1}$ is
perpendicular to $e_{1},e_{2},$ $\ldots,$ $ e_{2\ell
-1},e_{2\ell}$. We put 
 $$e_{2\ell+2} =
(\sec\alpha)Pe_{2\ell+1},\,\,\,\,e_{(2\ell+1)^{*}}:=
e_{n+2\ell+1}=
(\csc\alpha)Fe_{2\ell+1},\leqno(2.2)$$
$$e_{(2\ell+2)^{*}}:=
e_{n+2\ell+2}= (\csc\alpha)Fe_{2\ell+2}.$$

\noindent If $N$ is totally real in $M$, that is, if the slant angle
$\,\alpha={{\pi}\over {2}}$, then we can just choose
$e_{1},\ldots,e_{n}$ to be any local orthonormal frame of
$TN$ and put
$$e_{n+1}:=e_{1^{*}}=Je_{1},\ldots,e_{2n}:=e_{n^{*}}=Je_{n}.$$

If $m>n$, then at each point $p\in N$ there exist a
 subspace $\,\nu_p \,$ of the normal space $T_{p}^{\perp}N$
such that $ \, \nu_p \,$ is invariant under the action of
the complex structure $J$ of $M$ and $$T_{p}^{\perp}N=
F(T_{p}N)\oplus \nu_{p}, \,\,\,\,\nu_{p} \perp
F(T_{p}N).\leqno(2.3)$$ We choose a local orthonormal frame
$e_{4k+1},\ldots,e_{2m}$ of $\nu$ such that
$$e_{2n+2}=:
e_{(4k+1)^{*}}=Je_{4k+1},\ldots,e_{2m}=: e_{(2m-1)^{*}}=J
e_{2m-1},$$ that  is, $e_{2n+1},\ldots,e_{2m}$ is a
$J$-frame of $\,(\nu,J)$.

We call such an orthonormal frame
$$e_{1},e_{2},\ldots,e_{2k-1},e_{2k},e_{1{*}},e_{2^{*}},
\ldots,e_{(2k-1)^{*}},e_{(2k)^{*}},\leqno(2.4)$$
$$e_{4k+1},e_{(4k+1)^{*}},
\ldots,e_{2m-1},e_{(2m-1)^{*}}$$

\noindent an {\it adapted slant frame\/} of $N$ in $M$.

\vskip.1in
{\bf Lemma 2.1.} {\it Let $N$ be an $n$-dimensional proper
 slant submanifold of a Kaehlerian manifold $M$. If $N$ is
a Kaehlerian slant submanifold, then, with respect to an
adapted slant frame (2.4), we have
$$\omega_{i}^{j^*}=\omega_{j}^{i^*},\,\,\,or\,\,
equivalently,\,\,\,\, h^{j}_{ik}=h^{i}_{jk}\leqno(2.5)$$
\noindent for any $i,j,k=1,...,n$, where $\omega_{A}^{B}$
are the connection forms associated with the
adapted slant frame.}
\vskip.1in {\bf Proof.} Since $N$ is a
Kaehlerian slant submanifold,
  $\nabla P=0$ by definition. Thus, by applying  
by Lemma 3.5 of Chapter II, we
have $A_{FX}Y=A_{FY}X$ for any $X,Y$ tangent to $N$.
Therfore, we have (2.5) according to the definition of
adapted slant frame.
\vskip.1in
{\bf Remark 2.1.} If $N$ is a totally real submanifold of
a Kaehlerian ma-nifold, Lemma 2.1 was given in [CO1].
\vskip.1in
{\bf Corollary 2.2.} {\it If $N$ is a proper slant
surface of a Kaehlerian manifold $M$, then, with respect
to an adapted slant frame of $N$ in $M$, we have
$$\omega_{i}^{j^*}=\omega_{j}^{i^*},\,\,\,or\,\,
equivalently,\,\,\,\, h^{j}_{ik}=h^{i}_{jk}\leqno(2.6)$$
\noindent for any $i,j,k=1,2.$} 
\vskip.1in
{\bf Proof.} This Corollary follows immediately from Lemma 
2.1 and 
Theorem 3.4 of Chapter II.
\vskip.1in

For an $n$-dimensional proper slant submanifold $N$ of an
almost Hermitian manifold $M$, we define a canonical 
1-form $\Theta$ on $N$ by
 $$\Theta =
\sum_{i=1}^{n} \,\omega_{i}^{i^{*}}.\leqno(2.7)$$ 
\vskip.1in
{\bf Lemma 2.3.} {\it Let $N$ be an $n$-dimensional
proper slant submanifold with slant angle $\alpha$ in a
Kaehlerian manifold $M$. Then we have
$$\Theta=
\sum_{i} \,(tr\,h^{i^{*}})\omega^{i}\leqno(2.8)$$ and
$$\Theta(X)=-n(\csc\alpha)<tH,X>\leqno(2.9)$$
for any vector $X$ tangent to $N$, where $H$ denotes the
mean curvature vector of $N$ in $M$.}
\vskip.1in
{\bf Proof.} Equation (2.8) follows from
equation (1.9) of Chapter II, (2.7) and Corollary
2.2.

Since
$$n<tH,e_{j}>=-n<H,Fe_{j}>=-n(\sin\alpha)<H,e_{j^{*}}>$$
$$
=-(\sin\alpha) \,tr\,h^{j^{*}},$$
(2.8) implies (2.9).

\vskip.1in
Now we give
the following [CM3] \vskip.1in
{\bf Lemma 2.4.} {\it Let $N$ be an $n$-dimensional
proper slant submanifold of ${\bf C}^n$.
 If $N$ is Kaehlerian slant, then the
canonical 1-form $\Theta$ is closed, that is, $d\Theta =0$.
Hence, $\Theta$ defines a canonical cohomology class on
$N$:} 
$$[\Theta] \in H^{1}(N;{\bf R}).\leqno(2.10)$$
\vskip.1in
{\bf Proof.} Under the hypothesis, formula (1.8) 
 of Chapter II gives
$$d\Theta = -\sum_{i,j=1}^{n}\omega_{i}^{j}\wedge
\omega_{j}^{i^{*}} -\sum_{i,j=1}^{n}\omega_{i}^{j^{*}}
\wedge\omega_{j^{*}}^{i^{*}}.\leqno(2.11)$$ Since
$\omega_{i}^{j^{*}}=\omega_{j}^{i^{*}}$ by Lemma 2.1 and
$\omega_{i}^{j}=-\omega_{j}^{i},$
$\omega_{i^{*}}^{j^{*}}=-\omega_{j^{*}}^{i^{*}}$, we have
$$\sum_{i,j}\omega_{i}^{j}\wedge\omega_{j}^{i^{*}} =
\sum_{i,j}\omega_{i}^{j^{*}}\wedge\omega_{j^{*}}^{i^{*}}=
0.\leqno(2.12)$$ 
From (2.11) and (2.12) we obtain the lemma.

\vskip.1in
{\bf Remark 2.2.} In fact, this lemma holds without 
the condition that $N$ is Kaehlerian slant (see Theorem
3.1.). However, the proof for the general case is much
more complicated. 
 \vskip.1in  Now we give the main result of this
section [CM3]. \vskip.1in {\bf Theorem 2.5.} {\it Let $N$
be a proper slant surface in ${\bf C}^2$ with slant angle
$\alpha$. Put $\Psi =
(2\sqrt{2}\pi)^{-1}({\csc\alpha})\Theta.$ Then $\Psi$
defines a canonical integral class on $N$:} $$\psi =
[\Psi]\in H^{1}(N;{\bf Z}).\leqno(2.13)$$ \vskip.1in
{\bf Proof.} Let $N$ be a proper slant surface with slant
angle $\alpha$ in ${\bf C}^2$. Denote by $\nu$ the Gauss
map $\nu : N \rightarrow G(2,4) \cong S^{2}_{+}\times
S^{2}_{-}$ and by $\nu_+$ and $\nu_-$ the projections:
$\nu_{\pm} : N \rightarrow G(2,4) \rightarrow
S^{2}_{\pm}$ (cf. Section 3 of Chapter III). 

From formula
(4.3) of Chapter III, we have
$$(\nu_{-})_{*} = {1\over \sqrt{2}}\{(-\omega_{1}^{4}+
\omega_{2}^{3})\eta_{5} + (\omega_{1}^{3}+\omega_{2}^{4})
\eta_{6}\}.\leqno(2.14)$$

By  Lemma 2.1 we obtain $\omega_{1}^{4}=\omega_{2}^{3}$.
Thus, (2.14) implies
$$(\nu_{-})_{*}={1\over
\sqrt{2}}\,\Theta\,\eta_{6},\leqno(2.15)$$
where $\eta_{6}={1\over 2}(e_{1}\wedge e_{4}-e_{2}\wedge
e_{3}).$ Now because $\nu_-$ is given by
$$(\nu_{-})(p)={1\over 2}(e_{1}\wedge e_{2}-e_{3}\wedge
e_{4})(p)$$
for any point $p\in N$, the slantness of $N$ in ${\bf
C}^2$ implies that the image $\nu_{-}(N)$ lies in the
small circle $S^{1}_{\alpha}$ of $S^{2}_-$. Moreover, it is
easy to see that $\sqrt{2}\,\eta_6$ is a unit  vector
tangent to $S^{1}_{\alpha}$. Let $\omega =
dS^{1}_{\alpha}$ be the arclength element of
$S^{1}_{\alpha}$. Then, for any vector $X$ tangent to $N$,
we have 
$$((\nu_{-})^{*}\omega)(X)=\omega((\nu_{-})_{*}(X)) =
\omega({1\over {\sqrt{2}}}\Theta(X)\eta_{6}) ={1\over
2}\Theta(X).\leqno(2.16)$$  Hence we have
$$(\nu_{-})^{*}\omega = {1\over 2}\Theta.\leqno(2.17)$$
Therefore, for any closed loop $\gamma$ on $N$, we have
$$\int_{\gamma} {1\over 2}\Theta = \int_{\gamma}
(\nu_{-})^{*}\omega = \int_{\nu_{-}(S^{1}_{\alpha})}
\omega = (index \,\,of\,\,\nu_{-})\,vol(S^{1}_{\alpha})$$
$$= (\sqrt{2}\pi\sin\alpha)(index\,\, of\,\,\nu_{-}).$$

\noindent This implies that for any closed loop $\gamma$ in $N$,
$\int_{\gamma} \Psi \in {\bf Z}$. Thus by Lemma 2.4
we obtain (2.13). This completes the proof of the theorem.
\vskip.1in
As an application  we obtain the following
[CM3]
\vskip.1in
{\bf Theorem 2.6.} {\it Let $N$ be a complete, oriented,
proper slant surface in ${\bf C}^2$. If the mean
curvature of $N$ is bounded below by some positive
constant $c > 0$, then, topologically, $N$ is either a
circular cylinder or a 2-plane.}
 \vskip.1in
{\bf Proof.} Consider the map $\nu_{-} : N \rightarrow
S^{2}_-$. Assume the slant angle of  $N$ in ${\bf C}^2$ is
$\alpha$. Then $\nu_{-}(N) \subset S^{1}_{\alpha}$. Since
$$(\nu_{-})(X)={1\over {\sqrt{2}}}(\omega_{1}^{3}
+\omega_{2}^{4})(X)\eta_{6}={1\over {\sqrt{2}}}(\sum_{i}
h^{3}_{ii}\omega^{1}+\sum_{i}h^{4}_{ii}\omega^{2})$$
for any vector $X$ tangent to $N$ by Lemma 2.1, the
assumption on the mean curvature implies that the map
$\nu_-$ is an onto map. Furthermore, since the rank of
$\nu_-$ is equal to one, a result of Ehresmann implies
that $\nu_-$ is in fact a fibration. Because $N$ is not
compact by Theorem 1.5, topologically, $N$ is either the
product of a line and a circle or a 2-plane. 

\vfill\eject
\noindent{}  
\vskip.3in
\noindent \S 3. COHOMOLOGY OF SLANT SUBMANIFOLDS.
\vskip.2in
One of the purposes of this section is to improve Lemma 2.4
to obtain the following Theorem 3.1 of [CM3]. The other
purpose is to prove that every proper slant submanifold
in any Kaehlerian manifold is a sympletic manifold
(Theorem 3.4) with the sympletic structure  induced from
the canonical endomorphism $P$.
 \vskip.1in 
{\bf Theorem 3.1.} {\it Let $N$ be an $n$-dimensional
proper slant submanifold of ${\bf C}^n$. Then the
canonical 1-form $\Theta$ defined by (2.7) is closed, that
is, $d\Theta=0$. Hence, $\Theta$ defines a canonical
cohomology class on $N$:}
 $$[\Theta]\in H^{1}(N;{\bf
R}).\leqno(3.1)$$

\vskip.1in
If $N$ is an $n$-dimensional  proper slant submanifold of
${\bf C}^n$ with slant angle $\alpha$, then as we already
known the dimension $n$ is even. Let $n=2k$ and
$e_{1},\ldots,e_{n},e_{1^{*}},\ldots,e_{n^{*}}$  an
adapted slant (orthonormal) frame of $N$ in ${\bf C}^n$.
Then we have
$$e_{2}=(\sec\alpha)Pe_{1},\ldots,e_{2k}=(\sec\alpha)
Pe_{2k-1},\leqno(3.2)$$
$$e_{1^{*}}=(\csc\alpha)Fe_{1},e_{2^{*}}=(\csc\alpha)
Fe_{2},
\ldots,e_{(2k)^{*}}=(\csc\alpha)Fe_{2k}.$$ 

By direct
computation we also have
$$te_{i^{*}}=-(\sin\alpha)e_{i},\,\,\,\,\,i=1,...,2k,
\leqno(3.3)$$
$$Pe_{2j}=-(\cos\alpha)e_{2j-1},\leqno(3.4)$$
$$fe_{(2j-1)^{*}}=-(\cos\alpha)e_{(2j)^{*}},\,\,\,
fe_{(2j)^{*}}=(\cos\alpha)e_{(2j-1)^{*}},\,\,\,
j=1,\ldots,k.$$

\vskip.1in
In order to  prove Theorem 3.1, we need  the
following lemmas which can be regarded as generalizations
of Lemma 2.1. 
\vskip.1in 
{\bf Lemma 3.2.} {\it Let $N$ be an $n$-dimensional
($n=2k$) proper slant submanifold of ${\bf C}^n$. Then,
with respect to an adapted slant frame, we have

$$\omega_{2j-1}^{(2i-1)^{*}}-\omega_{2i-1}^{(2j-1)^{*}}
=\cot\alpha(\omega_{2i-1}^{2j}-\omega_{2j-1}^{2i}),
\leqno(3.5)$$
$$\omega_{2j}^{(2i-1)^{*}}-\omega_{2i-1}^{(2j)^{*}}
=\cot\alpha(\omega_{2i}^{2j}-\omega_{2i-1}^{2j-1}),
\leqno(3.6)$$
$$\omega_{2i}^{(2j)^{*}}-\omega_{2j}^{(2i)^{*}}
=\cot\alpha(\omega_{2i}^{2j-1}-\omega_{2j}^{2i-1}),
\leqno(3.7)$$
$$\omega_{2j-1}^{(2i-1)^{*}}-\omega_{2i-1}^{(2j-1)^{*}}
=\cot\alpha(\omega_{(2i-1)^{*}}^{(2j)^{*}}-
\omega_{(2j-1)^{*}}^{(2i)^{*}}),
\leqno(3.8)$$
$$\omega_{2i}^{(2j)^{*}}-\omega_{2j}^{(2i)^{*}}
=\cot\alpha(\omega_{(2i)^{*}}^{(2j-1)^{*}}-
\omega_{(2j)^{*}}^{(2i-1)^{*}}),
\leqno(3.9)$$
$$\omega_{(2i-1)^{*}}^{(2j-1)^{*}}-\omega_{2i-1}^{2j-1}
=\cot\alpha(\omega_{2i-1}^{(2j)^{*}}-\omega_{2i}^{(2j-1)^{*}})
\leqno(3.10)$$
$$\omega_{(2i-1)^{*}}^{(2j)^{*}}-\omega_{2i-1}^{2j}
=\cot\alpha(\omega^{2i-1}_{(2j-1)^{*}}-\omega_{2i}^{(2j)^{*}})
\leqno(3.11)$$
$$\omega_{(2i)^{*}}^{(2j)^{*}}-\omega_{2i}^{2j}
=\cot\alpha(\omega_{2i-1}^{(2j)^{*}}-
\omega_{2i}^{(2j-1)^{*}}),
\leqno(3.12)$$
$$\omega_{2j-1}^{(2i)^{*}}-\omega_{2i}^{(2j-1)^{*}}
=\cot\alpha(\omega_{(2i)^{*}}^{(2j)^{*}}-
\omega_{(2i-1)^{*}}^{(2j-1)^{*}}),
\leqno(3.13)$$
for any $i,j=1,\ldots,k.$}
\vskip.1in
{\bf Proof.} From the definition of adapted slant frames
we have
$$<Je_{2i-1},e_{2j-1}>=0,\hskip.3in
i,j=1,\ldots,k.\leqno(3.14)$$

By taking the derivative of (3.13) with respect to a
tangent vector $X$ of $N$ and applying (3.2), we have
$$0=<J{\tilde \nabla}_{X}e_{2i-1},e_{2j-1}>
+<Je_{2i-1},{\tilde\nabla}_{X}e_{2j-1}>$$
$$=-<\nabla_{X}e_{2i-1},Pe_{2j-1}>-<h(e_{2i-1},X),
Fe_{2j-1}>$$
$$+<Pe_{2i-1},\nabla_{X}e_{2j-1}>
+<Fe_{2i-1},h(e_{2j-1},X)>$$
$$=-(\cos\alpha)(<\nabla_{X}e_{2i-1},e_{2j}>
-<\nabla_{X}e_{2j-1},e_{2i}>)$$
$$+(\sin\alpha)(<e_{(2i-1)^{*}},h(X,e_{2j-1})>
-<e_{(2j-1)^{*}},h(X,e_{2i-1})>).$$

\noindent This implies
$$(\cot\alpha)(\omega_{2i-1}^{2j}-\omega_{2j-1}^{2i})(X)
\leqno(3.15)$$ $$=<A_{e_{(2i-1)^{*}}}e_{2j-1}-
A_{e_{(2j-1)^{*}}}e_{2i-1},X>.$$ 

Therefore, by applying
formula $(1.9)'$  of Chapter II, we get formula
(3.5).

Similarly, by taking the derivatives of the following
equations:

$$<Je_{2i-1},e_{2j}>=(\cos\alpha)\delta_{ij},\,\,\,
<Je_{2i},e_{2j}>=0,$$ 
\vskip.02in
$$<Je_{(2i-1)^{*}},e_{(2j-1)^{*}}>
=0,\,\,\,
<Je_{(2i)^{*}},e_{(2j)^{*}}>=0,$$
\vskip.02in
$$<Je_{2i-1},e_{(2j-1)^{*}}>=(\sin\alpha)\delta_{ij},\,\,\,
<Je_{2i-1},e_{(2j)^{*}}>=0,$$
\vskip.02in
$$<Je_{2i},e_{(2j)^{*}}>=(\sin\alpha)\delta_{ij},\,\,\,
<Je_{(2i)^{*}},e_{(2j)^{*}}>=0,$$

\noindent  we obtain formulas (3.6)-(3.13),
respectively. 

This proves the lemma.
\vskip.1in
{\bf Lemma 3.3.} {\it  Let $N$ be an $n$-dimensional
($n=2k$) proper slant submanifold of ${\bf C}^n$. Then,
with respect to an adapted slant frame, we have
\vskip.02in
$$\omega_{2i}^{(2j)^{*}}+\omega_{2i-1}^{(2j-1)^{*}}
=\omega_{2j}^{(2i)^{*}}+\omega_{2j-1}^{(2i-1)^{*}},
\leqno(3.16)$$
$$\omega_{(2i)^{*}}^{(2j)^{*}}-\omega_{(2i-1)^{*}}^{(2j-1)^{*}}
=\omega_{2i}^{2j}-\omega_{2i-1}^{2j-1},
\leqno(3.17)$$
$$\omega_{2j}^{2i-1}-\omega_{(2j)^{*}}^{(2i-1)^{*}}
=\omega_{2i}^{2j-1}-\omega_{(2i)^{*}}^{(2j-1)^{*}}
\leqno(3.18)$$
\vskip.02in
\noindent for any $i,j=1,\ldots,k.$}
\vskip.1in
{\bf Proof.} Formula (3.16) follows from (3.5) and (3.7).
Formula (3.17) follows from (3.10) and (3.12). And
formula (3.18) follows from (3.7) and (3.9).
\vskip.1in
Now we give the proof of Theorem 3.1.
\vskip.1in
{\bf Proof.} From the definition we have
$$\Theta=\sum_{\ell=1}^{2k}\,\omega_{\ell}^{{\ell}^{*}}.\leqno(3.19)$$
Thus from the structure equations we have
$$-d\Theta=\sum_{i,j=1}^{k}\omega_{2i}^{2j}
\wedge\omega_{2j}^{(2i)^{*}}
+\sum_{i,j=1}^{k}\omega_{2i}^{2j-1}
\wedge\omega_{2j-1}^{(2i)^{*}}$$
$$+\sum_{i,j=1}^{k}\omega_{2i}^{(2j)^{*}}
\wedge\omega_{(2j)^{*}}^{(2i)^{*}}
+\sum_{i,j=1}^{k}\omega_{2i}^{(2j-1)^{*}}
\wedge\omega_{(2j-1)^{*}}^{(2i)^{*}}\leqno(3.20)$$
$$+\sum_{i,j=1}^{k}\omega_{2i-1}^{2j}
\wedge\omega_{2j}^{(2i-1)^{*}}
+\sum_{i,j=1}^{k}\omega_{2i-1}^{2j-1}
\wedge\omega_{2j-1}^{(2i-1)^{*}}$$
$$+\sum_{i,j=1}^{k}\omega_{2i-1}^{(2j)^{*}}
\wedge\omega_{(2j)^{*}}^{(2i-1)^{*}}
+\sum_{i,j=1}^{k}\omega_{2i-1}^{(2j-1)^{*}}
\wedge\omega_{(2j-1)^{*}}^{(2i-1)^{*}}.$$

By using formula (3.17) we have

$$\sum_{i,j=1}^{k}\omega_{2i}^{(2j)^{*}}
\wedge\omega_{(2j)^{*}}^{(2i)^{*}}
+\sum_{i,j=1}^{k}\omega_{2i-1}^{2j-1}
\wedge\omega_{2j-1}^{(2i-1)^{*}}$$
$$+\sum_{i,j=1}^{k}\omega_{2i-1}^{(2j-1)^{*}}
\wedge\omega_{(2j-1)^{*}}^{(2i-1)^{*}}
+\sum_{i,j=1}^{k}\omega_{2i}^{2j}
\wedge\omega_{2j}^{(2i)^{*}}$$

$$=\sum_{i,j=1}^{k}\omega_{2i}^{(2j)^{*}}
\wedge\omega_{(2j)^{*}}^{(2i)^{*}}
+\sum_{i,j=1}^{k}\omega_{2i-1}^{2j-1}
\wedge\omega_{2j-1}^{(2i-1)^{*}}$$
$$+\sum_{i,j=1}^{k}\omega_{2i-1}^{(2j-1)^{*}}\wedge
(\omega_{2j-1}^{2i-1}-\omega_{2j}^{2i}+
\omega_{(2j)^{*}}^{(2i)^{*}})
+\sum_{i,j=1}^{k}\omega_{2i}^{2j}
\wedge\omega_{2j}^{(2i)^{*}}$$

$$=\sum_{i,j=1}^{k}\omega_{2i}^{(2j)^{*}}
\wedge\omega_{(2j)^{*}}^{(2i)^{*}}
+\sum_{i,j=1}^{k}\omega_{2i-1}^{(2j-1)^{*}}\wedge
\omega_{(2j)^{*}}^{(2i)^{*}}$$
$$-\sum_{i,j=1}^{k}\omega_{2i-1}^{(2j-1)^{*}}\wedge
\omega_{2j}^{2i}
+\sum_{i,j=1}^{k}\omega_{2i}^{2j}
\wedge\omega_{2j}^{(2i)^{*}}$$
$$=\sum_{i,j=1}^{k}(\omega_{2i-1}^{(2j-1)^{*}}
+\omega_{2i}^{(2j)^{*}})\wedge\omega_{(2j)^{*}}^{(2i)^{*}}
+\sum_{i,j=1}^{k}\omega_{2j}^{2i}\wedge(\omega_{2i-1}^{(2j-1)^{*}}
+\omega_{2i}^{(2j)^{*}}).$$

Since $\omega_{2i-1}^{(2j-1)^{*}}
+\omega_{2i}^{(2j)^{*}}$ 
is  symmetric in $i$ and
$j$ by Lemma 3.3 and $\omega_{2i}^{2j}$ and
$\omega_{(2j)^{*}}^{(2i)^{*}}$ are skew-symmetric in $i$
and $j$, we obtain
$$\sum_{i,j=1}^{k}\omega_{2i}^{(2j)^{*}}
\wedge\omega_{(2j)^{*}}^{(2i)^{*}}
+\sum_{i,j=1}^{k}\omega_{2i-1}^{2j-1}
\wedge\omega_{2j-1}^{(2i-1)^{*}}\leqno(3.21)$$
$$+\sum_{i,j=1}^{k}\omega_{2i-1}^{(2j-1)^{*}}
\wedge\omega_{(2j-1)^{*}}^{(2i-1)^{*}}
+\sum_{i,j=1}^{k}\omega_{2i}^{2j}
\wedge\omega_{2j}^{(2i)^{*}}=0.$$

Moreover, by (3.18), we have
$$\sum_{i,j=1}^{k}\omega_{2i}^{2j-1}
\wedge\omega_{2j-1}^{(2i)^{*}}
+\sum_{i,j=1}^{k}\omega_{2i}^{(2j-1)^{*}}
\wedge\omega_{(2j-1)^{*}}^{(2i)^{*}}$$
$$+\sum_{i,j=1}^{k}\omega_{2i-1}^{2j}
\wedge\omega_{2j}^{(2i-1)^{*}}
+\sum_{i,j=1}^{k}\omega_{2i-1}^{(2j)^{*}}
\wedge\omega_{(2j)^{*}}^{(2i-1)^{*}}$$

$$=\sum_{i,j=1}^{k}\omega_{2i}^{2j-1}
\wedge\omega_{2j-1}^{(2i)^{*}}
+\sum_{i,j=1}^{k}\omega_{2i}^{(2j-1)^{*}}
\wedge(\omega_{(2i-1)^{*}}^{(2j)^{*}}
+\omega_{2j-1}^{2i}
-\omega_{2i-1}^{2j})$$
$$+\sum_{i,j=1}^{k}\omega_{2i-1}^{2j}
\wedge\omega_{2j}^{(2i-1)^{*}}
+\sum_{i,j=1}^{k}\omega_{2i-1}^{(2j)^{*}}
\wedge\omega_{(2j)^{*}}^{(2i-1)^{*}}$$

$$=\sum_{i,j=1}^{k}\omega_{2i}^{2j-1}
\wedge\omega_{2j-1}^{(2i)^{*}}
-\sum_{i,j=1}^{k}\omega_{2j}^{2i-1}
\wedge\omega_{2i}^{(2j-1)^{*}}$$
$$-\sum_{i,j=1}^{k}\omega_{2i}^{(2j-1)^{*}}
\wedge\omega_{(2j)^{*}}^{(2i-1)^{*}}
+\sum_{i,j=1}^{k}\omega_{2i-1}^{(2j)^{*}}
\wedge\omega_{(2j)^{*}}^{(2i-1)^{*}}$$

$$=\sum_{i,j=1}^{k}(\omega_{2j}^{2i-1}-
\omega_{(2j)^{*}}^{(2i-1)^{*}}
)\wedge(\omega_{2i-1}^{(2j)^{*}}-
\omega_{2i}^{(2j-1)^{*}}).$$

Since $\omega_{2j}^{2i-1}-
\omega_{(2j)^{*}}^{(2i-1)^{*}}$ is symmetric in $i$ and
$j$ by Lemma 3.3 and
$\omega_{2i-1}^{(2j)^{*}}-\omega_{2i}^{(2j-1)^{*}}$ is
skew-symmetric in $i$ and $j$ by formula (3.6) of Lemma
3.2, we obtain $$\sum_{i,j=1}^{k}\omega_{2i}^{2j-1}
\wedge\omega_{2j-1}^{(2i)^{*}}
+\sum_{i,j=1}^{k}\omega_{2i}^{(2j-1)^{*}}
\wedge\omega_{(2j-1)^{*}}^{(2i)^{*}}\leqno(3.22)$$
$$+\sum_{i,j=1}^{k}\omega_{2i-1}^{2j}
\wedge\omega_{2j}^{(2i-1)^{*}}
+\sum_{i,j=1}^{k}\omega_{2i-1}^{(2j)^{*}}
\wedge\omega_{(2j)^{*}}^{(2i-1)^{*}}=0$$

From (3.20), (3.21) and (3.22) we obtain the theorem.

\vskip.1in
{\bf Remark 3.1.} If $N$ is not a slant submanifold of
${\bf C}^n$, then the 1-form $\Theta$ is not closed in
general. For example, if $N$ is the standard unit 2-sphere
$S^2$ in $E^{3}\subset E^4$, then $N$ is not slant with
respect to any compatible complex structure on $E^4$.

Now, assume that $S^2$ is parametrized by
$$x(\theta,\varphi)=(\sin\varphi \cos\theta,
\sin\varphi \sin\theta, \cos\varphi,0).\leqno(3.23)$$
Put
$$e_{1}=(-\sin\theta,\cos\theta,0,0),\leqno(3.24)$$
$$e_{2}=(\cos\varphi \cos\theta, \cos\varphi \sin\theta,
-\sin\varphi,0).$$

\noindent Then, with respect to the complex structure $J_0$, we have
$$Fe_{1}=(-{1\over 4}\sin 2\theta  \sin 2\varphi,
-{1\over 2}\sin^{2}\theta \sin 2\varphi,
-\sin\theta \cos^{2}\varphi, \cos\theta),\leqno(3.25)$$
and $$Fe_{2}=(\cos^{2}\theta \sin\varphi, \sin\theta
\cos\theta \sin\varphi,\cos\theta \cos\varphi, \sin\theta
\cos \varphi).\leqno(3.26)$$ 

Let
$$Fe_{1}=||Fe_{1}||e_{1^{*}},\,\,\,
Fe_{2}=||Fe_{2}||e_{2^{*}}.$$

\noindent Then, with respect to the local orthonormal  frame
$\,e_{1},e_{2},e_{1^{*}},e_{2^{*}}\,$ on the open subset
of $S^2$ on which $\theta\not\equiv {\pi\over 2}$ and
$\varphi\not\equiv 0$ $(mod \,\pi)$, we have
$$\Theta=\omega_{1}^{1^{*}}+\omega_{2}^{2^{*}}
=(1-\sin^{2}\theta\cos^{2}\varphi)^{-1}({1\over
2}\sin\theta \sin 2\varphi d\theta - \cos\theta
d\varphi).$$ It is easy to see that $d\Theta\not=0.$

\vskip.1in
For an $n$-dimensional ($n=2k$) proper slant submanifold
$N$ of a Kaehlerian manifold $(M,g,J)$, we put 
$$\Lambda(X,Y)=\, <X,PY>\leqno(3.27)$$
for any vectors $X, Y$ tangent to $N$, that is, $\Lambda =
\Omega_{J | TN}$. Then, from Lemmas 1.1 to 1.3, 
$\Lambda$ is a non-degenerate 2-form on $N$, that is,
$\Lambda^{k}\not= 0$. 

\vskip.1in
We recall that an even-dimensional manifold is called a
{\it sympletic ma-nifold\/} if it has a
non-degenerate closed 2-form.

Now we give  the second main result of this
section.
  \vskip.1in
{\bf Theorem 3.4.} {\it Let $N$ be an
$n$-dimensional proper slant submanifold of a Kaehlerian
manifold $M$. Then $\Lambda$ is closed, that is,
$d\Lambda=0.$ 
Hence $\Lambda$ defines a canonical
cohomology class of $N$:
$$[\Lambda]\in H^{2}(N;{\bf
R}).\leqno(3.28)$$
In particular, $(N,\Lambda)$ is a sympletic manifold.}
\vskip.1in 
{\bf Proof.} By definition of the exterior
differentiation we have
$$d\Lambda(X,Y,Z)={1\over 3}\{X\Lambda(Y,Z) +
Y\Lambda(Z,X)+Z\Lambda(X,Y)$$
$$-\Lambda([X,Y],Z)
-\Lambda([Y,Z],X)-\Lambda([Z,X],Y)\}.$$

\noindent Thus, by the definition of $\Lambda$, we obtain
$$d\Lambda(X,Y,Z)={1\over 3}\{
<{\nabla_{X}}Y,PZ>+<Y,{\nabla_{X}}(PZ)>$$
$$+<{\nabla_{Y}}Z,PX>+<Z,{\nabla_{Y}}(PX)>
+<{\nabla_{Z}}X,PY>$$ 
$$+<X,{\nabla_{Z}}(PY)>
-<[X,Y],PZ>-<[Z,X],PY>$$ 
$$-<[Y,Z],PX>\}$$
$$={1\over 3}\{<Y,{\nabla_{X}}(PZ)>+<Z,{\nabla_{Y}}(PX)>
+<X,{\nabla_{Z}}(PY)>$$
$$+<{\nabla_{X}}Z,PY>+<{\nabla_{Y}}X,PZ>+
<{\nabla_{Z}}Y,PX>\}.$$

\noindent Therefore, by the definition of ${\nabla P}$, we obtain
$$d\Lambda(X,Y,Z)={1\over 3}\{<X,(\nabla_{Z}P)Y>
+<Y,(\nabla_{X}P)Z>\leqno(3.29)$$
$$+<Z,(\nabla_{Y}P)X>\}.$$

\noindent Therefore, by applying formula (3.2) of Chapter
II and formula (3.29), we get
$$d\Lambda(X,Y,Z)={1\over 3}\{<X,th(Y,Z)>+<X,A_{FY}Z>$$
$$+<Y,th(Z,X)>+<Y,A_{FZ}X>\leqno(3.30)$$
$$+<Z,th(X,Y)>+<Z,A_{FX}Y>.$$

\noindent Consequently, by applying formulas (1.1), (1.4), (1.5)
of Chapter II and formula (3.30), we obtain (3.28). 

This
proves the theorem.
\vskip.1in
As an immediate consequence of Theorem 3.4 we obtain the
following 
\vskip.1in
{\bf Theorem 3.5.} {\it If $N$ is a compact
2k-dimensional proper slant submanifold of a Kaehlerian
manifold $M$, then 
$$H^{2i}(N;{\bf R})\not= 0\leqno(3.31)$$
for any $i = 1,\ldots,k.$}
\vskip.1in
In other words we have the following {\it non-immersion
theorem.} 
\vskip.1in
{\bf Theorem 3.5$'$.} {\it Let $\,N\,$ be a compact
2k-dimensional differentiable manifold such that
$H^{2i}(N;{\bf R})=0$ for some $i \in \{1,\ldots,k\}$.
Then $N$ cannot be immersed in any Kaehlerian manifold as a
proper slant submanifold.}

\vfill\eject

\vskip.3in
\noindent \S 4. STABILITY AND INDEX FORM.
\vskip.2in
The main purpose of this section is to present some
results conerning the stability and index of totally real
submanifolds, a special class of slant submanifolds, in a
Kaehlerian manifold.

Let $f:N \rightarrow M$ be an immersion from a compact
$n$-dimensional ma-nifold $N$ into an $m$-dimensional
Riemannian manifold $M$. Let $\{f_{t}\}$ be a 1-parameter
family of immersions of $N \rightarrow M$ with the
property that $f_{0}=f$. Assume the map $F: N\times
[0,1] \rightarrow M$ defined by $F(p,t)=f_{t}(p)$ is
differentiable. Then $\{f_{t}\}$ is called a {\it
variation\/} of $f$. A variation of
$f$ induces a vector field in $M$ defined along
the image of $N$ under $f$. We shall denote this field by
$\zeta$ and it is constructed as follows:

Let $\partial/\partial t$ be the standard vector field in
$N\times [0,1]$. We set
$$\zeta(p)=F_{*}({\partial\over {\partial t}}(p,0)).$$
Then $\zeta$ gives rise to cross-sections $\zeta^T$ and
$\zeta^N$ in $TN$ and $T^{\perp}N$, respectively. If we
have $\zeta^{T}=0$, then $\{f_{t}\}$ is called a {\it
normal variation\/} of $f$. For a given normal vector
field $\xi$ on $N$, $exp\,t\xi$ defines a normal
variation $\{f_{t}\}$ induced from $\xi$. We denote by
${\mathcal V}(t)$ the volume of $N$ under $f_t$ with respect
to the induced metric and by ${\mathcal V}'(\xi)$ and ${\mathcal
V}''(\xi)$, respectively, the values of the first and the
second derivatives of ${\mathcal V}(t)$ with respect to $t$,
evaluted at $t=0$. 

The following formula is well-known:
$${\mathcal
V}'(\xi)=-n\int_{N} <\xi,H> *1.\leqno(4.1)$$ 

For a compact minimal submanifold $N$ of  a Riemannian
manifold $M$, the second variation formula is given by
$${\mathcal V}''(\xi)=\int_{N} \{||D\xi||^{2}-{\bar
S}(\xi,\xi)- ||A_{\xi}||^{2}\}*1,\leqno(4.2)$$ 

\noindent  where ${\bar S}(\xi,\eta)$  is
defined by $${\bar S}(\xi,\eta) =\sum_{i=1}^{n} {\tilde
R}(\xi,e_{i},e_{i},\eta),\leqno(4.3)$$
  $e_{1},\ldots,e_n$  a local orthonormal frame of
$TN$ and ${\tilde R}$ 
the Riemann curvature tensor of the ambient manifold $M$.

Applying the Stokes theorem to the integral of the first
term of (4.2) (as Simons  did in [Si1]), we have
$$I(\xi,\xi)=:{\mathcal V}''(\xi)=\int_{N}
<L\xi,\xi>*1,\leqno(4.4)$$
 in which $L$ is a self-adjoint,  strongly elliptic
linear differential operator of the second order acting on
the space of sections of the normal bundle
  given by
$$L=-\Delta^{D}-{\hat A}-{\hat S},\leqno(4.5)$$
where $\Delta^D$ is the Laplacian operator
associated with the normal connection, $<{\hat 
A}\xi,\eta>=trace<A_{\xi},A_{\eta}>$, and $<{\hat
S}\xi,\eta>={\bar S}(\xi,\eta).$

The differential operator $L$ is called the {\it Jacobi
operator\/} of $N$ in $M$. The differential operator $L$
has discrete eigenvalues $\lambda_{1} < \lambda_{2} <\ldots
\,\nearrow \infty.$ We put $E_{\lambda}=\{\xi\in
\Gamma(T^{\perp}N)\, :\, L(\xi)=\lambda\xi \,\}.$ The
number of $\sum_{\lambda <0}dim(E_{\lambda})$ is called
the {\it index\/} of $N$ in $M$. A vector field $\xi$ in
$E_0$ is called a {\it Jacobi field.\/}

A minimal submanifold $N$ of $M$ is said to be {\it
stable\/} if ${\mathcal V}''(\xi) \geq 0$ for any normal
vector field $\xi$ of $N$ in $M$. Otherwise, $N$ is said
to be {\it unstable\/}. It is clear that $N$ in $M$ is
stable if and only if the index of $N$ in $M$ is equal to
0.

Concerning the  stability of totally real submanifold we
mention the following result of [CLN] obtained in
1980.

\vskip.1in
{\bf Proposition 4.1.}  {\it Let $N$ be a compact
n-dimensional minimal totally real submanifold in a real
2n-dimensional Kaehlerian manifold $M$. Then $N$ is
stable if and only if
$$I(JX,JX)=\int_{N} \{||\nabla X ||^{2}+S(X,X)-{\tilde
S}(X,X)\}*1\leqno(4.6)$$
is non-negative for every tangent vector field $X$ on
$N$, where $S$ and ${\tilde S}$ denote the Ricci forms
of $N$ and $M$, respectively.}
\vskip.1in
{\bf Proof.} Under the hypothesis, we have
$$D_{X}JY=J\nabla_{X}Y, \,\,\,A_{JX}Y=A_{JY}X,\,\,\,
{\tilde S}\cdot J = J\cdot {\tilde S}\leqno(4.7)$$
(see, pp.145-146 of [CO1]). From this we may obtain
$$ {\bar S}(X,X)={\tilde S}(X,X)- \sum_{i=1}^{n} {\tilde
R}(X,e_{i},e_{i},X).$$
From the equation of Gauss and (4.7) we may also obtain
$${\bar S}(JX,JX)={\tilde S}(X,X)-S(X,X)-||A_{JX}||^{2}$$
since $N$ is minimal. From these we obtain the result.
\vskip.1in
By using Proposition 4.1 we have the following results of
[CLN] also obtained in 1980.
\vskip.1in
{\bf Theorem 4.2.} {\it Let $N$ be a compact, totally
real  submanifold of a
Kaehlerian manifold $(M,g,J)$ with $dim_{\bf R}N =
dim_{\bf C}M$.
Then we have}

(1) {\it If $M$ has positive Ricci tensor and $N$ is
stable, then $H^{1}(N;{\bf R})=0$
and}

(2) {\it If $M$ has nonpositive Ricci tensor, then $N$ is
always stable.}
\vskip.1in
{\bf Proof.} Let $\varphi$ be the 1-form dual to a vector
field $X$ tangent to $N$. Then we have
$$\int_{N} \{||\nabla X ||^{2}+S(X,X)\}*1 = \int_{N} \{
{1\over 2} ||d\varphi ||^{2}+||\delta X ||^{2}\}*1,$$
where $\delta$ is the codifferential operator. Thus (4.6)
becomes $$I(JX,JX)=\int_{N} \{{1\over 2}||d\varphi
||^{2}+||\delta X||^{2}-{\tilde S}(X,X)\}*1$$
which implies the theorem.

  \vskip.1in {\bf Proposition
4.3.} {\it Let $N$ be a compact n-dimensional minimal
totally real submanifold in a real 2n-dimensional
Kaehlerian manifold $M$. Then $N$ is stable if $N$
satisfies condition (1) or (2) below and $N$ is unstable
if $N$ satifies condition (3):}

(1) {\it $i^{*}{\tilde S} \leq S$ where i is the inclusion:
$N \rightarrow M$.}

(2)  {\it $i^{*}{\tilde S} \leq 2S$ and the identity map
of $N$ is stable as a harmonic map.}

(3) {\it $i^{*}{\tilde S} > 2S$ and $N$ admits a nonzero
Killing vector field.}

 \vskip.1in
{\bf Proof.} Stability follows from (1) immediately by
Proposition 4.1 from  the (2) by Proposition 4.1 and the
fact that the second variation for the identity map is
$\int_{N} \{||{\nabla X} ||^{2}-S(X,X)\}*1 $ (cf. [S1]).
Sufficiency of (3) follows from Proposition 4.1 and the
formula $\int_{N} \{||\nabla X ||^{2}-S(X,X)\}*1=0$ for
Killing vector field $X$.
  \vskip.1in For two normal vector fields
$\xi,\eta$ to a minimal submanifold $N$ in $M$, their {\it
index form} is defined by
$$I(\xi,\eta)=\int_{N} <L\xi,\eta>*1.\leqno(4.8)$$

It is easy to see that the index form $I$ is a symmetric
bilnear form; $I: T^{\perp}N \times T^{\perp}N
\rightarrow {\bf R}$.
For a vector subbundle $V$ of the normal bundle
$T^{\perp}N$, we denote by $I_{V}$ the restriction of the
index form on $V$. Thus, $I_V$ is a symmetric bilinear
form on $V \,$; $I_{V}:V\times V \rightarrow {\bf R}$.

By the index of $I_{V}$, denoted by $index(I_{V})$, we
mean the number of negative eigenvalues of the index form
$I_V$.

The normal bundle of
 a totally real submanifold $N$ in a
Kaehlerian manifold $(M,g,J)$
 has the following orthogonal
decomposition:
$$T^{\perp}N = J(TN)\oplus\nu,\,\,\,\,J(TN)\perp \nu.$$

For totally real minimal submanifold of a Kaehlerian
manifold of higher codimension we have the following result
of the author and J. M. Morvan.
\vskip.1in
{\bf Proposition 4.4} {\it Let $N$ be a compact,
$n$-dimensional, minimal, totally real
submanifold of a Kaehlerian manifold of complex dimension
$n+p \,(p>0)$. If 
 $M$ has  non-positive holomorphic bisectional curvatures,
then the index form satisfies
$$I(\xi,\xi)+I(J\xi,J\xi) \geq 0\leqno(4.9)$$
for any normal vector field $\xi$ of $N$ in $M$.}
 \vskip.1in
{\bf Proof.} Let $N$ be a compact, $n$-dimensional,
minimal, totally real submanifold of a Kaehlerian
manifold $M$ of complex dimension $n+p$ with $p>0.$
 Then, for any normal vector field $\xi$ in the
normal subbundle $\nu$ and vector fields $X,Y$ tangent to
$N$, we have
$$<D_{X}\xi,JY>=-<{\tilde \nabla}_{X}J\xi,Y> =
<A_{J\xi}X,Y>.\leqno(4.10)$$

This implies
$$||D\xi ||^{2} \geq ||A_{J\xi}||^{2},\,\,\,\,
||DJ\xi ||^{2}\geq ||A_{\xi}||^{2},\leqno(4.11)$$
for any normal vector field $\xi$ in $\nu$.
 By using (4.2), (4.4), (4.10) and (4.11) we find
$$I(\xi,\xi)+I(J\xi,J\xi) \geq -\int_{N} \sum_{i=1}^{n}
\{{\tilde R}(\xi,e_{i},e_{i},\xi) +{\tilde R}(J\xi,
e_{i},e_{i},J\xi)\}*1.$$
Therefore, if $M$ has non-positive holomorphic
bisectional curvature, then, for any normal vector field
$\xi$ in $\nu$, we have

$$I(\xi,\xi)+I(J\xi,J\xi)\geq 0.$$ 

\noindent This proves the Proposition.
\vskip.1in
{\bf Example 4.1.} 
Let $N$ be any non-totally geodesic, minimal
hypersurface of an $(n+1)$-dimensional flat real torus
${\bf R}T^{n+1}$ which is imbedded in a complex
$(n+1)$-dimensional flat complex torus ${\bf C}T^{n+1}$ as
a totally geodesic, totally real submanifold. Denote by
$\xi$ a unit normal vector field of $N$ in ${\bf
R}T^{n+1}$. Then we have $$A_{\xi}\not= 0, \,\,\,D\xi=
DJ\xi = 0,\,\,\, A_{J\xi}=0.$$

\noindent Since ${\bf C}T^{n+1}$ is flat, (4.2) and
(4.3) yield
$$I(\xi,\xi) < 0,\,\,\, I(J\xi,J\xi)>0.\leqno(4.12)$$

\vskip.1in

For the index of the index form $I_{J(TN)}$ we have the
following result of the author and Morvan.
\vskip.1in
{\bf Theorem 4.5.} {\it Let $N$ be a compact,
$n$-dimensional, totally real, mi-nimal submanifold of a
Kaehlerian manifold $M$ of complex dimension $n+p$. If $M$
has non-positive holomorphic bisectional curvatures, then}
$$index(I_{J(TN)}) =0.\leqno(4.13)$$ 
\vskip.1in
{\bf Proof.} Under the hypothesis, let 
$e_{1},\ldots,e_{2n+2p}$ be a local orthonormal frame
along the submanifold $N$ such that
$$e_{n+1}=Je_{1},\ldots,e_{2n}=Je_{n},\,\,\,
e_{2n+p+1}=Je_{2n+1},\ldots,e_{2n+2p}=Je_{2n+p}.$$
Then $e_{n+1},\ldots,e_{2n}$ form a local orthonormal frame
of $J(TN)$ and $e_{2n+1},$ $\ldots,$ $e_{2n+2p}$ form a
local orthonormal frame of $\nu$.
By  applying
Lemma 3.5 of Chapter II and the equation of Gauss, we may
obtain $$\sum_{i=1}^{n} {\tilde R}(X,e_{i},e_{i},X) -
S(X,X)=\sum_{i=1}^{n}
<h(X,e_{i}),h(X,e_{i})>\leqno(4.14)$$ $$=||A_{JX}||^{2} +
\sum_{i=1}^{n} \sum_{r=2n+1}^{2n+2p}
(<h(X,e_{i}),e_{r}>)^{2}.$$

Combining (4.2) with (4.14) we may find
$$I(JX,JX)=\int_{N}
\{||DJX||^{2}+S(X,X)+\sum_{i,r}(<h(X,e_{i}),e_{r}>)^{2}
\leqno(4.15)$$
$$ - \sum_{t=1}^{2n} {\tilde
R}(X,e_{t},e_{t},X)\}*1.$$

As we did in [CLN] (see [C3, p.51]),  put
$$W=\nabla_{X}X + (div\,X)X$$
where $div\,X$ denotes the divergence of $X$. Let
$\varphi$ be the 1-form associated with $X$. Then, by
computing the divergence of $W$ and applying the
divergence theorem, we obtain
$$0=\int_{N} (div\,W)*1=\int_{N}\{S(X,X)+||\nabla X||^{2}
- {1\over
2}||d\varphi||^{2}-(\delta\varphi)^{2}\}*1.\leqno(4.16)$$
Combining (4.15) and (4.16) we find
$$I(JX,JX)=\int_{N}
\{||DJX||^{2}+\sum_{i,r}(<h(X,e_{i}),e_{r}>)^{2}
\leqno(4.17)$$
$$ - \sum_{t=1}^{2n} {\tilde
R}(X,e_{t},e_{t},X)-||\nabla X ||^{2} +{1\over
2}||d\varphi ||^{2} + (\delta\varphi)^{2}\}*1,$$
for  any vector field $X$ tangent to $N$.

For a totally real submanifold $N$ of a Kaehlerian
manifold, the Gauss and Weingarten formulas imply
$$D_{X}JY=J\nabla_{X}Y + fh(X,Y)\leqno(4.18)$$

\noindent  which
yields $||DJY|| \geq ||\nabla Y||$. Thus, by applying
 (4.17), (4.18) and the definition of the
index, we obtain the theorem.
 \vskip.1in
{\bf Remark 4.1.} An $n$-dimensional totally real
submanifold of a real $2n$-dimensional Kaehlerian
manifold is also called a {\it Lagrangian submanifold\/}
by some mathematicians. 

\vfill\eject

\vskip.3in
\noindent \S 5.  STABILITY OF TOTALLY GEODESIC 
SUBMANIFOLDS. 
\vskip.2in
In this section we would like to present a general method
introduced by the author, Leung and Nagano  [CLN]
obtained in 1980 for
 determining the stablity of totally geodesic
submanifolds in compact symmetric spaces. Since every
irreducible totally geodesic submanifold  of a Hermitian
symmetric space is a slant submanifold [CN1], the method
can be used to determine the stability of such slant
submanifolds.

Recall that the second variation formula of a compact
minimal subma-nifold $N$ in a Riemannian manifold is
given by formula (4.2) (or, equivalently, by formula
(4.4)).  If $N$ is
totally geodesic, then $A=0$. So the stability obtained
trivially when ${\bar S}$ is non-positive. For this
reason we are interested in the case $M$ is of compact
type. Also we assume $M$ is irreducible partly to
preclude tori as $M$. 

We need to fix some notations. Since $N$ is totally
geodesic, there is a finitely covering group $G_N$ of
the connected isometry group $G_{N}^o$ of $N$ such that
$G_N$ is a subgroup of the connected isometry group $G_M$
of $M$ and leaves $N$ invariant, provided that $G_{N}^o$
is semi-simple. Let $\mathcal P$ denote the orthogonal
complement of the Lie algebra ${\it g}_{_{N}}$
in the Lie algebra ${\it g}_{_{M}}$ with respect to the
bi-invariant inner product on ${\it g}_{_{M}}$ which is
compatible with the metric of $M$. Every member of ${\it g}_{_{M}}$ 
is though of as a Killing vector field because of the
action $\, G_M$ on $\, M$. Let $\hat P$ denote the space of
the vector fields corresponding to the member of $\mathcal P$
restricted to the submanifold $N$.
\vskip.1in
{\bf Lemma 5.1.} {\it To every member of $\, \mathcal P$ there
corresponds a unique (but not canonical) vector field $v
\in {\hat P}, v$ is a normal vector field and hence
${\hat P}$ is a $G_N$-invariant subspace of the
space, $\Gamma (T^{\perp}N),$ of the sections
of the normal bundle to $N$. Moreover, ${\hat P}$ is
homomorphic with $\mathcal P$ as a $G_N$-module.}
\vskip.1in
{\bf Proof.} Let $o$ be an arbitrary point of $N$. Let
$K_M$ and $K_N$ denote the isotropy subgroups of $G_M$
and $G_N$ at $o$, respectively. Then  ${\it g}_{_{M}}/
{\it k}_{_{M}}$ and ${\it g}_{_{N}}/
{\it k}_{_{N}}$ and ${\mathcal P}/({\mathcal P} \cap {\it
k}_{_{M}})$ are identified with $T_{o}M$, $T_{o}N$ and 
$T_{o}^{\perp}N$ by isomorphisms induced by the
evaluation of vector fields in ${\it g}_{_{M}}$ at $o$.
In particular, the value $v(o)$ of $v$ is normal to $N$.
This proves the lemma. \vskip.1in

Now we are ready to explain the method of [CLN].
\vskip.1in 
The group $G_N$ acts on sections in
$\Gamma(T^{\perp}N)$ and hence on the differential
opeators: $\Gamma(T^{\perp}N) \rightarrow
\Gamma(T^{\perp}N)$. $G_N$ leaves $L$ fixed since $L$ is
defined with $N$ and the metric of $M$ only. Therefore,
each eigenspace of $L$ is left invariant by $G_N$. Let $V$
be one of its $G_N$-invariant irreducible subspaces. We
have  a representation $\rho : G_{N} \rightarrow GL(V).$ We
denote by $c(V)$ or $c(\rho)$ the eigenvalue of the
corresponding  Casimir operator. 

 To define
$c(V)$ we fix an orthonormal basis $(e_{\lambda})$ for
${\it g}_{_{N}}$ and consider the linear endomorphism $C$
or $C_V$ of $V$ defined by 
$$C=-\sum\,\rho(e_{\lambda})^{2}.\leqno(5.1)$$
It is known that $C$ is $c(V)I_V$ (see Chapter 8 of
[Bo1]), where $I_V$ is the identiy map on $V$. In our case,
the {\it Casimir operator\/} $C_{V}=-\sum [e_{\lambda}
,[e_{\lambda},V]\,]$ for every member $v$ of $V$ (after
extending to a neighborhood of $N$).

Now, we state the following theorem of [CLN] which says,
modulo details, that {\it $N$ is stable if and only if
$c(V)\geq c({\mathcal P})$ for every $G_N$-invariant
irreducible space $V$.}
 \vskip.1in
{\bf  Theorem 5.2.} {\it A compact totally geodesic
submanifold $N \,(= G_{N}/K_{N})$ of a compact symmetric
space $M \,(=G_{M}/K_{M})$ is stable as a minimal
submanifold if and only if one has $c(V)\geq c(P')$ for
the eigenvalue of the Casimir operator of every simple
$G_N$-module $V$ which shares as a $K_N$-module some simple
$K_N$-submodule of the $K_N$-module $T_{o}^{\perp}N$ in
common with some simple $G_N$-submodule $P'$ of $\hat P$.}
\vskip.1in
{\bf  Proof.} 
Given a point $p$ of $N$, we choose a basis of ${\it
g}_{_{N}}$ given by 
$$(e_{\lambda})=
(\ldots,e_{i},\ldots,e_{\alpha},\ldots) $$  
and a finite
system $(e_{r})$ of vectors in ${\mathcal P} \subset {\it
g}_{_{M}}$ such that (1) $(e_{i}(p))_{1\leq i \leq n}$ is
an orthonormal basis for the tangent space $T_{p}N$, (2)
${\nabla e_{i}}=0,\, 1\leq i \leq n,$ at $p$, (3)
$e_{\alpha}(p)=0,\, n < \alpha \leq dim\,{\it g}_{_{N}}$,
and (4) $(e_{r}(p))$ is an orthonormal basis for the
normal space $T_{p}^{\perp}N$, which we can do as is
well-known.

An arbitrary normal vector field $\xi$ is written as
$\xi=\sum \xi^{r}e_r$ on a neighborhood of $p$ by Lemma 5.1.
Since $N$ is totally geodesic in $M$, we have
$D_{X}\xi={\tilde \nabla}_{X}\xi$ for $X$ tangent
to $N$ and $\xi\in \Gamma(T^{\perp}N)$, where ${\tilde
\nabla}$ is the Riemannian connection of $M$. As the
curvature tensor of a symmetric space we also have
${\tilde R}(X,Y)Z=-[[X,Y],Z]$  for $X, Y, Z \in {\it
m}=:{\it g}_{_{M}}/{\it k}_{_{M}}$. Therefore, by
 evaluating $L\xi$ and $C\xi$ at
$p$, we obtain 
$$L\xi= -\sum\,D_{e_{i}}D_{e_{i}}\xi -\sum\,{\bar
S}_{rs}\xi^{r}e_{s}$$ 
$$= -\sum\,{\tilde \nabla}_{e_{i}}{\tilde 
\nabla}_{e_{i}}\xi -\sum\,{\bar
S}_{rs}\xi^{r}e_{s},$$ 
where ${\bar S}_{rs}$ are the
components of ${\bar S}$ in (4.3) and 
$$C\xi= -\sum
\,[e_{\lambda},[e_{\lambda},\xi]\,]$$
$$=-\sum\,{\tilde \nabla}_{e_{i}}{\tilde \nabla}_{e_{i}}\xi
+\sum\, \xi^{r}{\tilde \nabla}_{e_{i}}{\tilde
\nabla}_{e_{r}}e_{i} -\sum \,\xi^{r}({\tilde
\nabla}_{e_{r}}e_{\alpha})^{s}{\tilde \nabla}_{e_{s}}
e_{\alpha}$$
$$=-\sum \,{\tilde \nabla}_{e_{i}}{\tilde \nabla}_{e_{i}}
\xi
 -\sum
\,\xi^{r}({\tilde \nabla}_{e_{r}}e_{\alpha})^{s}
{\tilde \nabla}_{e_{s}}e_{\alpha},$$ 
where for the vanishing of the
second term we use the fact that $e_i$ is a Killing vector
field. Thus we find 
$$(L-C)\xi=\sum\,
\xi^{r}({\tilde
\nabla}_{e_{r}}e_{\alpha})^{s}\,{\tilde \nabla}_{e_{s}}
e_{\alpha} -\sum\,
{\bar S}_{rs}\xi^{r}e_{s}$$ $$=\sum\,
(A_{\alpha})^{2}\xi-{\bar S}(\xi, *),$$ where  
$A_{\alpha}$ is the Weingarten map given by
the restriction of the operator: $X \rightarrow
-{\tilde \nabla}_{X}e_{\alpha}$ on $T_{p}N$ to the normal space
$T_{p}^{\perp}N$. This proves the following two
statements:

 {(a)} {\it The difference $\, L-C \,$ is an
operator of order one and}

 {(b)} {\it The difference $\, L-C \,$ is given by
a self-adjoint endomorphism ${\breve S}$ of the normal
bundle $T^{\perp}N$. }

Now, because both $L$ and $C$ are $G_N$-invariant,
statement (b) implies

 {(c)} {\it The endomorphism $\,{\breve S} \,$ is
$G_N$-invariant.}

The theorem follows from statements (a) through (c) easily
when the isotropy subgroup $K_N$ is irreducible on the
normal space $T_{p}^{\perp}N \cong {\it g}_{_{M}}/({\it
g}_{_{N}} \oplus {\it k}_{_{M}})$. In fact, ${\breve S}$ is then a
constant scalar multiple of the identity map of
$T^{\perp}N$; ${\breve S}=k\cdot I,$ by statements (a) through (c)
and Schur's lemma. $N$ is stable if and only if the
eigenvalues of $L$ are all non-negative. Since $L=
(c(P')+k)\cdot I=0$ on the normal Killing fields (see
Lemma 5.1), this is equivalent to say that $0 \leq
c(V)+k=c(V)-c(P')$ for every simple $G_N$-module $V$ in
$\Gamma (T^{\perp}N)$ (which is necessarily contained in
an eigenspace of $L$ by ${\breve S} = k\cdot I),$ and the
Bott-Frobenius theorem completes the proof.

In the general case, we decompose the normal space into
the direct sum of simple $K_N$-modules: $P'\oplus
P''\oplus \ldots \,.$ Accordingly we have
$T^{\perp}N=E'\oplus E'' \oplus \ldots,$ where $E',
E'',\ldots,$ etc. are obtained from $P', P'',\ldots,$
etc., in the usual way by applying the action of $G_N$ to
the vectors in $P', P'',\ldots,$ etc. Since $G_N$ leaves
invariant $E', E'',\ldots,$ the normal connection leaves
invariant the section spaces $\Gamma (E'), \Gamma
(E''),\ldots \, .$ Hence, $L$ and $C$ leave these spaces
invariant. (For this, the irreducible subspaces $P', P'',
\ldots$ must be taken within eigenspaces of the symmetric
operator ${\bar S}$ at the point o). In particular, the
projections of $\Gamma (T^{\perp}N)$ onto $E', E'',
\ldots,$ etc. commute with $C$ and $L$. Thus one can
repeat the argument for irreducible case to each of $E',
E'',\ldots $ to finish the proof of the theorem.
\vskip.1in
 Theorem 5.2 provides us  an {\it algorithm for
stability} which goes like this: One can compute $c(V)$ by
the Freudenthal formula (cf.  Chapter 8, p.120 of [Bo1])
once one knows the action $\rho$ of $G_N$ on $V$. So the
rest is to know all the simple $G_N$-modules $V$ in $\Gamma
(T^{\perp}N)$. This is done by means of the Frobenius
theorem as reformulated by Bott, which asserts in our
case that a simple $G_N$-module $V$ appears in $\Gamma
(T^{\perp}N)$ if and only if $V$ as a $K_N$-module
contains a simple $K_N$-module which is isomorphic with a
$K_N$-module of $T_{o}^{\perp}N$.
\vskip.1in
{\bf Remark 5.1.}  A reformulation of the method of [CLN]
was given by Y. Ohnita  in [Oh2]. 
\vskip.1in
{\bf Proposition 5.3.} {\it A compact totally geodesic
submanifold $N$ of a compact symmetric space $M$ is
unstable as a minimal submanifold if the normal bundle
admits a nonzero $G_N$-invariant section and if the
centeralizer of $G_N$ in $G_M$ is discrete.}
\vskip.1in
{\bf Proof.} Let $\xi$ be a nonzero $G_N$-invariant
normal vector field on $N$. We have $D\xi
=0$. In view of (4.2) we will show that ${\bar
S}(\xi,\xi)$ is positive. The sectional curvature of a
tangential 2-plane at a point $p\in N$ equals
$||[e,f]||^{2}$ if (i) $e$ is a member of ${\it
g}_{_{N}}$, (ii) $f$ is that of ${\it g}_{_{M}}$, (iii)
$e(p)$ and $f(p)$ form an orthonormal basis for the
2-plane, and (iv) ${\tilde \nabla}e={\tilde\nabla}f=0$ at
$p$. Therefore, ${\bar S}(\xi,\xi)$ fails to be positive
only if $[e,f]=0$ for every such $e$ and $f$ satisfying
$\xi(p)\wedge f(p)=0.$ Since the isotropy subgroup $K_N$
at $p$ leaves the normal vector $\xi(p)$ invariant, we
have $[e',f]=0$ for every member $e'$ of ${\it k}_{_{M}}$
and hence $[{\it g}_{_{N}},f]=0$ if ${\bar S}(\xi,\xi)=0$
at $p$. Such an $f$ generates a subgroup in the
centeralizer of ${\it g}_{_{N}}$ in ${\it g}_{_{M}}$. This
contradicts to the assumption. 
\vskip.1in
{\bf Example 5.1.} Let $N$ be the equator in the sphere
$M=S^n$. That $N$ is unstable follows from the
Proposition 5.3 if one consider a unit
$G_N$-invariant normal vector field to it. The centralizer
in this case is generated by the antipodal map: $x
\rightarrow -x.$ Its orbit space is the real projective
space $M'$. The projection: $M \rightarrow M'$ carries
$N$ onto a hypersurface $N'$. The reflection in $N'$ is a
member of $G_N$ by our general agreement on $G_N$ (if
$n>1$) and precludes the existence of non-vanishing
$G_N$-invariant normal vector field to $N'$. It is clear
 by Theorem 5.2 that $N'$ is stable.

\vskip.1in
{\bf Remark 5.2.} In general, if $N$ is a stable minimal
submanifold of a Riemannian manifold $M$ and $M$ is a
covering Riemannian manifold of $M'$, then the projection
$N'$ of $N$ in $M$ is stable too. The example above shows
the converse is false.

\vskip.1in
{\bf Definition 5.1.} For a compact connected symmetric
space $M$, there is a unique symmetric space $M^*$ of
which $M$ and every connected symmetric space which is
locally isomophic with $M$ are covering Riemannian
manifold of $M^*$. We call $M^*$ the bottom space of $M$.
Ig $M$ is a group manifold, $M^*$ is the adjoint group
$ad(M)$.
\vskip.1in
By applying Theorem 5.2 and Proposition 5.3 above, we may
obtain the following result of [CLN].

  \vskip.1in {\bf Theorem 5.4.} {\it A
compact subgroup $N$ of a compact Lie group $M$ is stable
with respect to a bi-invariant metric on $M$ if }
(a) $N$ {\it has the same rank as $M$ and}
(b) {\it $M=M^{*}$, that is, $M$ has no nontrivial
center.}
\vskip.1in
{\bf Proof.} The compact group manifold $M$ has
$G_{M}=M_{L}\times M_{R},$ where $M_L$ is the left
trranslation group $M\times \{1\}$ and $M_R$ the right
translation group; here $M_R$ acts ``to the left'' too,
that is, $(1,a)$ carries $x$ into $xa^{-1}$. Similarly
for $G_N$. $G_N$ is effective on every invariant
neighborhood of $N$ in $M$ by (b). We first consider the
case where $N$ is a maximal toral subgroup $T$ of $M$.
Let $A_T$ denote the subgroup $\{(a,a^{-1}):a\in T\}$ of
$G_N$. We have an epimorphism $\epsilon :
K_{N}\times A_{T} \rightarrow G_T$ by the
multiplication whose kernel is the subgroup of
elements of order 2. In order to use Theorem 5.2, we look
at an arbitrary simple $G_T$-module $V$ in $\Gamma(E')$
where $E'$ is, as before, the vector bundle $G_{N}P'$
defined from the simple $K_T$-submodule $P'$ of the
normal space. $P'$ is a root space corresponding to a
root $\alpha$ of ${\it g}_{_{M}}$. With $V$ we compare
the space $P'$, a simple $K_N$-module in $\Gamma(E')$
which is defined from the members of the Lie algebra of
$M_L$ taking values in $P'$ at a point of $N$. We want to
show $c(V)\geq c(P').$ Since $\alpha\not= 0$ by (a), both
$V$ and $P'$ have dimension 2 and these are isomorphic as
$K_T$-modules. The relationship between $V$ and $P'$ can
be made more explicitly. Namely, a basis for $P'$ is a
global frame of $E'$ and so the sections in $V$ are
linear combinations of the basis vectors whose components
are functions on $N$. These functions form a simple
$G_T$-module $F$ of dimension 2 and $V$ is a
$G_T$-submodule of $F\otimes P'$. By the Bott-Frobenius
theorem, $K_T$ acts trivially on a 1-dimensional subspace
of $F$. Every weight $\varphi$ of $F$ is a linear
combination of roots of ${\it g}_{_{M}}$ whose
coefficients are even numbers. In fact all the weights of
the representations of $G_T$ are linear combinations of
those roots over the integers by (a) and (b) and, since
$K_{T}\cap A_{T} \cong ker\,\epsilon$ is trivial on $F$,
the coefficients must be even.

On the other hand, if one looks at the definition of
Casimir operator, $C=-\sum\,\rho(e_{\lambda})^2$, one
sees that the eigenvalue $c(V)$ is a sort of average of
the eigenvalues of $-\rho(x)^2$, $||x||=1,$ or more
precisely, $c(V)=-(dim\,V)^{-1}\int trace(\rho(x))^2$,
where the integral is taking over the unit sphere of the
Lie algebra with an appropriately normalized invariant
measure. For this reason, showing $c(V)\geq c(P')$, or
equivalently, $(\varphi+\alpha)^{2}-\alpha^{2}\geq 0$
amounts to showing the inner product $<2\alpha + \varphi
,\alpha
>=<\varphi+\alpha,\varphi+\alpha>-<\alpha,\alpha>\,\,\geq
0$ in which we may assume that $\varphi$ is dominant (with
respect to the Weyl group of ${\it g}_{_{M}}$). This
proves the case $N=T$.

We trun to the general case $N\supset T.$ Assume $N$ is
unstable and will show this contradicts the stability of
$T$. There is then a simple $G_N$-module $V$ such that
the second variation (4.2) is negative for some member
$\xi$ of $V$. If we restrict $\xi$ to $T$ we still have a
normal vector field but the integrand in (4.2) for
$\xi_{|\,T}$ will differ from the restriction of the
integrand for $\xi$ by the terms corresponding to the
tangential directions to $N$ which are normal to $T$.
However, a remedy comes from the group action. First
(4.2) with $A=0$ is invariant under $G_N$ acting on $V$.
Second, every tangent vector to $N$ is carried into a
tangent vector to $T$ by some isometry in $G_N$. Third,
$N$ and $T$ are totally geodesic in $G$, but more
importantly the connection and the curvature restrict to
the submanifolds comfortably. And finally, the isotropy
subgroup $K_N$ acts irreducibly on the tangent space to
each simple or circle normal subgroup of $G_N$. From all
these it follows that (4.2) for $\xi$ is a positive
constant multiple of (4.2) for $\xi_{|\,T}$, as one sees
by integrating (4.2) for $g(\xi)_{|\,g(T)},\,g\in G$,
over the group $G$ and over the unit sphere of $V$. This
is a contradiction which completes the proof of the
theorem.
\vskip.1in
{\bf Remark 5.3.} Neither the assumption (a) or (b) can
be omitted from Theorem 5.4 as the examples of $M=SU(2)$
with $N=SO(2)$ and $M=G_2$ with $N=SO(2)$ show. Also the
Theorem will be false if $M$ is not a group manifold, a
counter-example being $M=M^{*}=GI$ with $N=S^{2}\cdot
S^2$ (local product).

\vskip.1in By applying Theorem 5.2,
Proposition 5.3 and Theorem 5.4 we may obtain the
following results of [CLN]. 
 \vskip.1in
{\bf Proposition 5.5.} 
(a) {\it Among the compact connected
simple Lie groups $M^{*}$, the only ones that have
unstable $M_{+}^*$ are $SU(n)^*$, $SO(2n)^*$ with n odd,
$E_{6}^*$ and $G_2$. }

(b) {\it The unstable $M_+$ are $G^{C}(k,n-k),\, 0<k<n-k,$
for $SU(n)^*$; $SO(2n)/U(n)^*$ for $SO(2n)^*$; $EIII^*$
for $E_{6}^*$; and $M_{+}^*$ for $G_2$.}

(c) {\it Every $M_-$ is stable for the group $M^*$.}

{\bf Comments on the Proof.} (I) The stability of $M_-$ is
immediate from Theorem 5.4 since $M_-$ has the same
rank as $M$ (cf. [CN1, II]). Otherwise the proof is based
on scrutinizing all the individual cases and omitted
except for a few cases to illustrate our methods.
(II) Take $M^{*}=SO(2n+1)$. Then $M_{+}=G_{+}/K_{+}=
G^{R}(k,2n+1-k),\, 0<k<n-k,$ the Grassmannians of the
unoriented $k$-planes in $E^{2n+1}$ by Table I in
[CN1, II]. The action of $G_{+}=SO(2n+1)$ on $P$ (in the
notation of Lemma 5.1) is the adjoint representation
corresponding to the highest weight ${\tilde \omega}_2$
in Bourbaki's notation [Bo1]. By Freudenthal's formula,
one finds that ${\tilde\omega}_1$ is the only
representation that has a smaller eigenvalue than
${\tilde\omega}_{2}; c({\tilde
\omega}_{1})<c({\tilde\omega}_{2}).$ But
${\tilde\omega}_1$ does not meet the Bott-Frobenius
condition simply because its dimension $2n+1$ is too
small. Therefore $M_+$ is stable by Theorem 5.2. (III) Take
$SU(n)^*$ for another example. We know
$M_{+}=G_{+}/K_{+}=G^{C}(k,n-k),$ the complex
Grassmann manifold. If $k\not= n-k,\, M_+$ is 1-connected
and hence $K_+$ is connected. On the other hand,
$M_{-}=K_{+}=S(U(k)\times U(n-k)),$ which contains a
circle group as the center. Therefore, $M_+$ admits a unit
$G_+$-invariant normal vector field. Moreover, the
centralizer  of $G_+$ in $G_M$ is trivial. Hence
Proposition 5.3 applies to conclude that $M_+$ is
unstable. This argument fails in the case $k=n-k$ and we
can conclude the stability of $M_+$ by Theorem 5.2 as in
(I). (IV) Unstability is established by means of
Proposition 5.3 except for the case of $G_2$. In this case
we have $c({\tilde\omega}_{1})<c({\tilde\omega}_{2})=c$
(the adjoint representation). This ${\tilde\omega}_1$
gives a monomorphism of $G_2$ into $SO(7)$ which restricts
to a monomorphism of $K_{+}=SO(4)$ into $SO(4)\times
SO(3)$ in $SO(7)$ and then projects to $SO(3)$. This
implies that ${\tilde\omega}_1$ appears in a space of
normal vector fields.
  
\vskip.1in 
{\bf Proposition 5.6.} {\it Let $M^*$ be a
compact symmetric space $G/K$ with $G$ simple. Then, among
the $M_+$ and $M_-$, the unstable minimal submanifolds are
$G^{R}(k,n-k),\,k<n-k,$ in $AI(n)^*$; $G^{H}(k,n-k),\,
k<n-k,$ in $AII(n)^*$; $SO(k)$  in $G^{R}(k,k)$ with k
odd; $M_{+} = M_{-}=SO(2)\times AI(n)$ in $CI(n)^*$;
$M_{+}=M_{-}=SO(2)\times AII({n\over 2})$ in $DIII^{*}=
SO(2n)/U(n)$ with n even; $G^{H}(2,2)$ in $EI^*$; $FII$ in
$EIV^*$; $AII(4)$ in $EV^*$; and $M_{+}=M_{-}= S^{2}\cdot
S^{2}$ in $GI$.}

{\bf Comments on the Proof.} (I) In some cases, one can
use another method to get the results quickly. For
instance, if $M^*$ is Kaehlerian, then it is well-known
that every compact complex submanifold is stable. (II)
Mostly, unstability is established by using Proposition
5.3. In the cases, $M_{+}=M_{-}=SO(2)\times L$, this
proposition does not literally apply but unstability is
proven in the same spirit. Consider, say $SO(2)\times
AI(n)$ in $CI(n)^*$. This space in $M^*$ is $U(n)/O(n)$.
The normal space is isomorphic with the space of the
symmetric bilinear forms on $E^n$ as an
$O(n)$-module. Therefore, there is a $U(n)$-invariant
unit normal vector field $\xi$ on $M_+$. We have ${\tilde
\nabla}\xi =0.$ We have to show ${\bar S}(\xi,\xi)>0$ in
view of (4.2). Since $M_{+}=M_-$ has the same rank as
$M$, there is a tangent vector $X$ in $T_{y}M_-$ such
that the curvature of the 2-plane spanned by $X$ and
$\xi(y)$ is positive. (III) The case of
$M_{+}=M_{-}=S^{2}\cdot S^2$ in $GI$. Precisely,
$M_{+}=M_-$ is obtained from $S^{2}\times S^{2} = $ (the
unit sphere in $E^{3}) \times$ (the unit sphere in
$E^{3}) \subset E^{3}\times E^3$ by
identifying $(x,y)$ with $(-x,-y)$. The group $G_-$ for
$M_{-}=G_{-}/K_-$ is the adjoint group but we have to
take its double covering group $SO(4)$ to let it act on a
neighborhood of $M_-$. The identity representation of
$SO(4)$ on $E^4$ restricts to the normal
representation of $K_{-}=SO(2)\times SO(2)$ as somewhat
detailed examination of the root system reveals.
Therefore, $M_-$ is unstable. Similarly for $M_+$ which is
congruent with $M_-$.
\vskip.1in
{\bf Remark 5.4.} From the known facts about geodesics,
one would not expect a simple relationship between
stability and homology. More specifically, we remark that
$M_+$ {\it is homologous to zero for a group manifold
$M^{*}$.\/} The proof may go like this. Consider the
quadratic map $f: x \mapsto s_{x}(o)$ on a symmetric
space $M=G/K$ for a fixed point $o$, where $s_x$ is the
symmetry at $x$. Assume $M$ is compact
and orientable. Then $f$ has a nonzero degree if and only
if the cohomology ring $H^{*}(M)$ is a Hopf algebra (cf.
M. Clancy's thesis, University of Notre Dame, 1980). On the
other hand, the inverse image $f^{-1}(o)$ is exactly $M_+$
and $\{o\}$. Since $H^{*}(M)$ is a Hopf algebra for a group
$M^*$, it follows that every $M_+$ is homologous to zero.
\vskip.1in

Finally, we give the folowing [CLN]

 \vskip.1in
{\bf Proposition 5.7.} {\it The minimal totally real
totally geodesic submanifold $G^{R}(p,q)$ is unstable
in $G^{C}(p,q)$.}
\vskip.1in
{\bf Proof.} Let $N= G_{N}/K_{N}$ be a  totally
real and totally geodesic submanifold of a compact
Kaehlerian symmetric space $M=G_{M}/K_M$. Then $N$ will
be unstable if we find $c(V)<c(P')$ as in Theorem 5.2.
For each simple $P'$ in ${\mathcal P}$, there is a simple
${\it g}_{_{N}}$-module $V$ in $\Gamma (E')$ whose
members are normal vector fields $\xi = JX$ for some
Killing vector field $X$ in ${\it g}_{_{N}}$. This is
obvious from the definition of a totally real
submanifold. In the case of $G^{R}(p,q)$ in $G^{C}(p,q)$,
${\mathcal P}$ is simple and $c(P)=c(2{\tilde
\omega}_{1})>c({\tilde \omega}_{2}) =c({\it g}_{_{N}}),$
where ${\tilde \omega}_{2}$ denotes the highest weight in
Bourbaki's notation (see, [Bo1]) and ${\tilde \omega}_{1}$
is the only representation that has a smaller eigenvalue
than ${\tilde \omega}_{2}$. This proves the Proposition.
\vskip.1in {\bf Remark 5.5.} If $p=1$, Proposition 5.7 was
due to
 [LS1].

\vskip.1in 
{\bf Remark 5.6.} The method of [CLN] was used by
several mathematicians in their recent studies. For 
results in this direction see, for instances, [MT1],
[MT2], [Oh2] and [Ta1]. 

{\bf Remark 5.7.} For the fundamental theory of
$(M_{+},M_{-})$ and some of its applications see
[C3], [C7], [CN1], [CN3] and [N2].

 Y. Ohnita (1987) improved the above algorithm to
include the formulas for the index, the
nullity and the Killing nullity of a
compact totally geodesic submanifold in a compact
symmetric space.

Let \ $f:N\to M$ \ be a compct totally geodesic
submanifold of a compact Riemannian symmetric
space. Then  $f:N\to M$ is expressed as follows:
There are compact symmetric pairs
$(G,K)$ and $(U,L)$ with $N=G/K,\, M=U/L$ so
that  $f:N\to M$ is
given by $uK\mapsto \rho(u)L$, where  $\rho:G\to
U$ is an analytic homomorphism
with $\rho(K)\subset L$ and the injective
differential $\rho:\mathfrak g\to\mathfrak u$ which
satisfies $\rho(\mathfrak m)\subset \mathfrak p$. Here
$\mathfrak u=\mathfrak l +\mathfrak p$ and $\mathfrak g=\mathfrak
k+\mathfrak m$ are the canonical decompositions of
the Lie algebras  $u$ and $g$, respectively.

Let $\mathfrak m^\perp$
denote the orthogonal complement of $\rho(\mathfrak
m)$ with $\mathfrak p$ relative to the
ad$(U)$-invariant inner product $(\;,\;)$ on
$\mathfrak u$ such that $(\;,\;)$ induces the
metric of $M$. Let $\mathfrak k^\perp$ be the
orthogonal complement of $\rho(\mathfrak k)$ in
$\mathfrak l$. Put $\mathfrak g^\perp=\mathfrak k^\perp
+\mathfrak m^\perp$. Then $\mathfrak g^\perp$ is the
orthogonal complement of $\rho(\mathfrak g)$ in
$\mathfrak u$ relative to $(\;,\;)$, and $\mathfrak
g^\perp$ is ad$_\rho (G)$-invariant. Let $\theta$
be the involutive automorphism of the symmetric
pair $(U,L)$. Choose an orthogonal decomposition\
$\mathfrak g^\perp=\mathfrak g_1^\perp\oplus\cdots
\oplus\mathfrak g_t^\perp$\ such that each $\mathfrak
g^\perp_i$ is an irreducible ad$_\rho
(G)$-invariant subspace with $\theta(\mathfrak
g^\perp_i)=\mathfrak g^\perp_i$. Then, by Schur's
lemma, the Casimir operator $C$ of the
representation of $G$ on each $\mathfrak g^\perp_i$ is
$a_i I$ for some $a_i\in{\bf C}$. 

Put ${\mathfrak m}_i^\perp ={\mathfrak m}\cap {\mathfrak
g}^\perp_i$ and let $D(G)$ denote the set of all
equivalent classes of finite dimensional
irreducible complex representations of $G$. For
each $\lambda\in D(G)$,
$(\rho_\lambda,V_\lambda)$ is a fixed
representation of $\lambda$. 

For each $\lambda\in
D(G)$, we assign a map $A_\lambda$ from
$V_\lambda\otimes \hbox{Hom}_K(V_\lambda,W)$ to
$C^\infty(G,W)_K$ be the rule $A_\lambda(v\otimes
L)(u)=L(\rho_\lambda(u^{-1})v)$. Here
Hom$_K(V_\lambda,W)$ denotes the space of all
linear maps $L$ of $V_\lambda$ into $W$ so that
$\sigma(k)\cdot L=L\cdot \rho_\lambda(k)$ for all
$k\in K$. 

Y. Ohnita's formulas for the index $i(f)$,
the nullity $n(f)$, and the Killing nullity
$n_k(f)$ are given respectively by:

(a) $i(f)=\sum_{i=1}^t\sum_{\lambda\in D(G),
a_\lambda <a_i} m(\lambda)d_\lambda$,

(b) $n(f)=\sum_{i=1}^t\sum_{\lambda\in D(G),
a_\lambda =a_i} m(\lambda)d_\lambda$,

(c) $n_k(f)=\sum_{i=1,\mathfrak
m_i^\perp\ne \{0\}}^t\dim \mathfrak g_i^\perp$,

\noindent where $m(\lambda)=\dim\,
\hbox{Hom}_K(V_\lambda,(\mathfrak m_i^\perp)^{\bf
C})$ and $d_\lambda$ denotes the dimension of the
representation $\lambda$

 By applying his formulas,
Ohnita  determined the  indices, the
nullities and the Killing nullities for all
totally geodesic submanifolds in compact
rank one  symmetric spaces.

 \vfill\eject

\vskip.6in

\centerline {\bf REFERENCES}
\vskip.4in
\noindent A. Bejancu

 {[Be1]} {\sl Geometry of CR-submanifolds,} D.
Reidel Publishing Co., Dor-
drecht-Boston-Lancaster-Tokyo, 1986
\vskip.1in

\noindent D. E. Blair

 {[Bl1]} {\sl Contact Manifolds in Riemannian
Geometry,} Lecture Notes in Math., {\bf 509},
Springer-Verlag, Berlin-New York, 1976.
 \vskip.1in

 \noindent D. E. Blair and B. Y.
Chen

 {[BC]} On CR-submanifolds of Hermitian manifolds,
{\sl Israel J. Math.,\/} {\bf 34} 
(1979), 353-363.
\vskip.1in
 
\noindent J. Bolton, G. R. Jensen, M. Rigoli and L. M.
Woodward

{[BJRW]} On conformal minimal immersion of $S^2$
into ${\bf C}P^n$, {\sl Math.\ Ann.,} 
{\bf 279} (1988),599-620.
 \vskip.1in
\noindent N. Bourbaki

{[Bo1]} {\sl Groupes et Algebres de Lie,
Chapters 7 \& 8}, Hermann, Paris, 1975.

 \vskip.1in 
\noindent E. Calabi

{[Ca1]} Isometric embeddings of complex
manifolds, {\sl Ann. of Math.,} {\bf 58} 
(1953), 1-23.

{[Ca2]} Metric Riemann surfaces, {\sl Ann. of
Math. Studies,} {\bf 58} (1953), 
1-23.

 \vskip.1in
\noindent  B. Y. Chen, 

{[C1]}{\sl Geometry of Submanifolds},
Mercel Dekker, New York, 1973.

 {[C2]} 
 Differential geometry of real
submanifolds in a Kaehler manifold, 
{\sl Monatsh. f\" ur
Math.,} {\bf 91} (1981), 257-274.

{[C3]} {\sl Geometry of Submanifolds and Its
Applications,} Science University of Tokyo, 1981.

{[C4.1]} $CR$-submanifolds of a Kaehler manifolds, I, {\sl J. Differential 
Geometry,} {\bf 16} (1981), 305-322.

{[C4.2]} $CR$-submanifolds of a Kaehler manifolds, II, {\sl J. Differential 
Geometry,} {\bf 16} (1981), 493-509.

{[C5]} {\sl Total Mean Curvature
and Submanifolds of  Finite Type}, World 
Scientific,
Singapore-New Jersey-London-Hong Kong, 1984.

 {[C6]}  Slant immersions, {\sl Bull.\ Austral.\
Math.\ Soc.,} {\bf 41} (1990), 135-147.

{[C7]} {\sl A New Approach to Compact Symmetric
Spaces and Applications: A report on joint work with
Professor T. Nagano,} Katholieke Universiteit Leuven, 1987.

\vskip.1in  \noindent B. Y. Chen, P. F. Leung and T. Nagano

{[CLN]} Totally geodesic submanifolds of
symmetric spaces, III, {\sl preprint,} 1980.
\vskip.1in  \noindent B. Y. Chen and G. D. Ludden

{[CL1]} Surfaces with mean curvature vector
parallel in the normal bundle, {\sl Nagoya Math. J.,}
{\bf 47} (1972), 161-167.

\vskip.1in 
\noindent B. Y. Chen, G. D. Ludden and S. Montiel
{[CLM]} Real submanifolds of a Kaehler
manifold, {\sl Algebras, Groups and Geometries,} {\bf 1}
(1984), 176-212.

\vskip.1in
\noindent B. Y. Chen, C. S. Houh and H. S. Lue

{[CHL]} Totally real submanifolds, {\sl J.
Differential Geometry,} {\bf 12} (1977), 473-480.

\vskip.1in
\noindent  B. Y. Chen and J. M. Morvan

{[CM1]} Propri\' et\' e riemanniennes
des surfaces lagrangiennes, {\sl C.R.\ Acad.\ Sc.\ Paris,\
Ser.} I, {\bf 301} (1985), 209-212.

 {[CM2]}  G\' eom\' etrie des
surfaces lagrangiennes de ${\bf C}^2$, {\sl J.\ Math.\
Pures et
Appl.,} {\bf 66}, (1987), 321-335.

{[CM3]} Cohomologie des sous-vari\'et\'es
$\alpha$-obliques, {\sl C. R. Acad. Sc. Paris,} {\bf 314}
(1992), 931--934.
\vskip.1in 

\noindent B. Y. Chen and T. Nagano

{[CN1]} Totally geodesic submanifolds of
symmetric spaces, I \& II, {\sl Duke Math. J.,} {\bf 44}
(1977), 745-755 \& {\bf 45} (1978), 405-425.

{[CN2]} Harmonic metrics, harmonic tensors and
Gauss maps, {\sl J. Math. Soc. Japan,} {\bf 36} (1984),
295-313.

{[CN3]} A Riemannian geometric invariant and
its applications to a problem of Borel and Serre, {\sl
Trans. Amer. Math. Soc.,} {\bf 308} (1988), 273-297.
 \vskip.1in  

\noindent  B. Y. Chen and K. Ogiue,

{[CO1]} On totally  real 
submanifolds, {\sl Trans.\ Amer.\ Math. Soc.,} {\bf
193} (1974), 257-266.
\vskip.1in

\noindent  B. Y. Chen and P. Piccinni

{[CP1]} Submanifolds with
finite  type  Gauss  map, {\sl Bull.\ Austral.\ Math.
\ Soc.,}
{\bf 35} (1987), 321-335.
\vskip.1in

\noindent  B. Y. Chen and Y. Tazawa

{[CT1]} Slant surfaces of 
codimension  two, {\sl Ann.\ Fac.\ Sc.\ Toulouse
Math.,} {\bf 11} (1990), 29--43.

{[CT2]}  Slant submanifolds in
complex Euclidean spaces, {\sl Tokyo J. Math.} {\bf 14}
(1991), 101--120. \vskip.1in

\noindent
B. Y. Chen and L. Vanhecke

{[CV]} Differential geometry of geodesic
spheres, {\sl J.\ Reine\ Angew. Math.,}
{\bf 325} (1981), 28-67.

\vskip.1in 
\noindent S. S. Chern 

{[Ch1]} {\sl Complex Manifolds without
Potential Theory,} 2nd Ed., Springer-Verlarg, Berlin-New
York, 1979.

  \vskip.1in  \noindent P. Dazord

{[D1]} Une interpr\' etation g\' em\'
etrique de la classe de Maslov-Arnold, {\sl J. 
Math.
Pures et Appl.,} {\bf 56} (1977), 231-150.
\vskip.1in

\noindent F. Dillen, B. Opozda, L. Verstraelen and L.
Vrancken

{[DOVV]} On totally real 3-dimensional
submanifolds of the nearly Kaehler 6-sphere, {\sl Proc.
Amer. Math. Soc.,} {\bf 99} (1987), 741-749.

\vskip.1in 
\noindent J. Eells and L. Lemaire

{[EL1]} A report on harmonic maps, {\sl Bull.\
London Math.\ Soc.,} {\bf 10} (1978), 1-68.

{[EL2]} Another report on harmonic maps, {\sl
Bull.\ London Math.\ Soc.,} 
{\bf 20} 
(1988), 385-524.

\vskip.1in 
\noindent C. Ehresmann

{[Eh1]} Sur la th\' eorie des vari\'et\'es
feuillet\'ees, {\sl Rend. di Mat.,\/} {\bf 10} (1951),
64-82.
\vskip.1in

\noindent J.  A. Erbacher

{[Er1]} Reduction of the codimension of an
isometric immersion, {\sl J. Differential Geometry,} {\bf
5} (1971), 333-340.

 \vskip.1in 
 \noindent L. Gheysens, P. Verheyen and L. Verstraelen

{[GVV]} Characterization and examples of Chen
submanifolds, {\sl Jour. of
Geometry,\/} 
{\bf 20} (1983), 47-72.

\vskip.1in 
\noindent S. Greenfield

{[G1]} Cauchy-Riemann equations in several
variables, {\sl Ann. della Scuola Norm. Sup. Pisa,} {\bf
22} (1968), 275-314.

\vskip.1in 
\noindent S. Helgason

{[He1]} {\sl Differential Geometry, Lie Groups,
and Symmetric Spaces,} 
Academic Press, New-York-San Francisco-London, 1978.
\vskip.1in

\noindent R. Harvey and H. B. Lawson

{[HL1]} Calibrated foliation, {\sl Amer. J.
Math.,} {\bf 104} (1982), 607-633.

 \vskip.1in
\noindent D. A. Hoffman and R. Osserman,

{[HO1]} The  Gauss  map  of surfaces
in $R^3$ and $R^4$, {\sl Proc.\ London\ Math.\ 
Soc.,} (3)
{\bf 50} (1985), 27-56.
\vskip.1in

\noindent C. S. Houh

{[Ho1]} Some totally real minimal surfaces in
$CP^2$, {\sl Proc. Amer. Math. Soc.} {\bf 40} (1973),
240-244.
 \vskip.1in  \noindent E. K\" ahler

{[Ka1]} \" Uber eine bemerkenswerte
Hermitische Metrik,  {\sl Abh. Math. Sem.
 Univ. Hamburg,} {\bf 9} (1933), 173-186.

\vskip.1in
\noindent S. Kobayashi

{[Ko1]} Recent results in complex differential
geometry, {\sl Jahresber. Deut-sch. Math.-Verein.,}
{\bf 83} (1981), 147-158.

{[Ko2]} {\sl Differential Geometry of Complex
Vector Bundles,}  Iwanami \&
Princeton University Press, 1987.
\vskip .1in

\noindent S. Kobayashi and K. Nomizu

{[KN1]} {\sl Foundations  of
Differential  Geometry,} Volume I, Johy Wiley 
and Sons, 1963.

{[KN2]} {\sl Foundations of Differential 
Geometry,} Volume  II, Johy Wiley 
and Sons,  1969.

\vskip.1in \noindent H. B.
Lawson

{[L1]} {\sl Lectures  on  Minimal
Submanifolds,} Publish or Perish, Berkeley, 
1980

\vskip.1in \noindent H. B. Lawson and J. Simons

{[LS1]} On stable currents and their
applications to global problems in real and complex
geometry, {\sl Ann. of Math.,} {\bf 98} (1973), 427-450.

\vskip.1in 
\noindent S. Maeda and S. Udagawa

{[MU1]} Surfaces with constant Kaehler angle all
of whose geodesics  are 
circles in a complex space form, preprint, 1990.
\vskip.1in 

\noindent K. Mashimo and H. Tasaki

{[MT1]} Stability of
maximal tori in compact Lie groups, {\sl Algebras, Groups
and Geometries,} to appear. 

{[MT2]} Stability of closed Lie subgroups in
compact Lie groups, {\sl Kodai Math. J.,} to appear.

\vskip.1in
\noindent Y. Miyaoka {[Mi1]}
Inequalities on Chern numbers, {\sl S$\hat u$gaku}, {\bf
41}  (1989), 193-207.

\vskip.1in 
\noindent J. D. Moore
{[Mr1]} Isometric immersions of Riemannian
products, {\sl J. Differential Geometry,} {\bf 5} (1971),
159-168.

\vskip.1in 
\noindent J. M. Morvan

{[Mo1]} Classe de Maslov d'une immersion
lagrangienne et minimaliti\' e, 
{\sl C.R. Acad. Sc. Paris,}
{\bf 292} (1981), 633-636.
\vskip.1in

\noindent T. Nagano

{[N1]} Stability of harmonic maps between
symmetric spaces, {\sl Lectures Notes in Math.,}
Springer-Verlag, {\bf 949} (1982), 130-137. 

{[N2]} The involutions of compact symmetric
spaces, {\sl Tokyo J. Math.,} {\bf 11} (1988), 57-79.

\vskip.1in 
\noindent K. Ogiue

{[O1]} Differential geometry of Kaehler
submanifolds, {\sl Advan. in 
Math.,} 
{\bf 13} (1974), 73-114.

{[O2]} Recent topics in submanifold theory, {\sl
S${\hat u}$gaku}, {\bf 39}
 (1987), 305-319.

\vskip.1in 
\noindent Y. Ohnita

{[Oh1]} Minimal surfaces with constant curvature
and Kaehler 
angle in 
complex space forms, {\sl Tsukuba J.\
Math.,} {\bf 13} (1989), 191-207.

{[Oh2]} On stability of minimal submanifolds in
compact symmetric spa-ces, {\sl Compositio Math.,} {\bf
64} (1987), 157-189.

 \vskip.1in 
\noindent B. O'Neill

{[On1]} {\sl Semi-Riemannian Geometry,}
Academic Press, New York-San Fra-ncisco-London, 1983.
\vskip.1in
\noindent B. Opozda

{[Op1]} Generic submanifolds in almost
Hermitian manifolds, {\sl Ann. Polon. Math.,} {\bf XLIX}
(1988), 115-128.
 \vskip.1in  
\noindent R. Osserman

{[Os1]} {\sl Survey of Minimal Surfaces,} Van
Nostrand Reinhold, New York, 
1969.

{[Os2]} Curvature in the eighties, {\sl Amer.\
Math.\ Monthly,} 1990.

\vskip.1in 
\noindent A. Ros and L. Verstraelen

{[RV1]} On a conjecture of K. Ogiue, {\sl J.
Differential Geometry,} {\bf 19}
(1984), 561-566.
 \vskip.1in
\noindent P. Ryan

{[Ry1]} Homogeneity and some curvature
conditions for hypersurfaces, {\sl Tohoku Math. J.,} {\bf
21} (1969), 363-388.

 \vskip.1in  \noindent J. A. Schouten
and D. van Dantzig

{[SD1]}  \" Uber unit\" are
 Geometrie, {\sl Math. Ann.,} {\bf 103}
(1930), 319-346.

{[SD2]} \" Uber  unit\" are
 Geometrie  konstanter Kr\" ummung, {\sl
Proc. Kon.
Nederl. Akad.
Amsterdam,} {\bf 34} (1931), 1293-1314. 
\vskip.1in

\noindent   J. Simons

{[Si1]} Minimal varieties in Riemannian
manifolds, {\sl Ann. of Math.,} {\bf 88} (1968), 62-105.
 \vskip.1in
\noindent  I. M. Singer and J. A. Thorpe

 {[ST1]} The curvature of 4-dimensional Einstein
spaces, {\sl Global Analysis,}
 Princeton University Press, 1969,
73-114.
 
\vskip.1in
\noindent R. T. Smyth

{[S1]} The second variation formula for
harmonic mappings, {\sl Proc. Amer. Math. Soc.,} {\bf 47}
(1975), 229-236.
 \vskip.1in  
\noindent B. Smyth

{[Sm1]} Differential geometry  of complex
hypersurfaces, {\sl Ann. of Math.,}
{\bf 85} (1967), 246-266.

\vskip.1in
\noindent M. Spivak

{[Sp1]} {\sl A Comprehensive Introduction to
Differential Geometry,\/} Vol. 4,
Publish or Perish, Berkeley 1979.
\vskip.1in

\noindent M. Takeuchi

{[Ta1]} Stability of certain minimal
submanifolds of compact Hermitian symmetric spaces, {\sl
Tohoku Math. J.,} {\bf 36} (1984), 293-314.

{[Ta2]} On the fundamental group of the group
of isometries of a symmetric space, {\sl J. Fac. Sci.
Univ. Tokyo, Sect. I,} {\bf 10} (1964), 88-123.
 \vskip.1in 

\noindent I. Vaisman

{[V1]} {\sl Sympletic Geometry and Secondary
Characterisstic Classes,} 
Birk-
h\" auser, Boston, 1987.
 \vskip.1in 

\noindent A. Weil

{[W1]} Sur la th\' eorie des formes
diff\' erentielles attach\' et\' e
analytique complexe, {\sl Comm. Math. Helv.,} {\bf 20}
(1947), 110-116.
\vskip.1in 
\noindent J. L. Weiner

{[We1]} The Gauss map for surfaces in 4-space,
{\sl Math. Ann.,} {\bf 269} (1984), 
541-560.
\vskip.1in 
\noindent A. Weinstein

{[Wn1]} {\sl Lectures on Sympletic Manifolds,}
Regional Conference Series in 
Mathematics,  No. {\bf 29}, American
Mathematical Society, Providence, 
1977.
\vskip.1in 
\noindent R. O. Wells, Jr.

{[Wl1]} {\sl Differential Analysis on Complex
Manifolds,} Springer-Verlag, Ber-lin-New York, 1980.

\vskip.1in 
\noindent W. Wirtinger, 

{[Wi1]} Eine Determinantenidentit\" at und
ihre Anwendung auf analytische Gebilde in Euclidischer und
Hermitischer Massbestimmung, 
{\sl Monatsch. Math.
Phys.,} {\bf 44} (1936), 343-365.
\vskip.1in

\noindent K. Yano

{[Y1]} {\sl Integral Formulas in Riemannian
Geometry,} Mercel Dekkker, New York, 1970.
 \vskip.1in
\noindent K. Yano and M. Kon

{[YK1]} {\sl Anti-invariant Submanifolds,}
Mercel Dekker, New York, 1976.

 {[YK2]} {\sl
CR-submanifolds of Kaehlerian and Sasakian Manifolds,}
Birkh\" au-ser, Boston, 1983.

\vfill\eject

\vskip.4in
\centerline {\bf SUBJECT INDEX}
\vskip.3in

\noindent $\alpha$-slant submanifold 36, 40, 46, 56

\noindent Adapted slant frame 29, 85

\noindent Algorithm for stability 108

\noindent Anti-holomorphic submanifold 36, 41

\noindent Autere submanifold 25, 26

\noindent Bott-Frobenius theorem 108

\noindent Casimir operator 106, 107

\noindent Canonical cohomolgy class 86, 87, 89, 95

\noindent Compatible complex structure 18, 35

\noindent Complex submanifold 8, 13

\noindent  $D_{1}(2,4)$ 33, 34, 77

\noindent Doubly slant submanifold 42, 43, 44, 45

\noindent Exact sympletic manifold 82
 
\noindent $f$ (canonical endomorphism on normal bundle) 14

\noindent $F$ (canonical 1-form on tangent bundle) 13,
25, 26, 27

\noindent Freudenthal formula 108

\noindent $g_+$, $g_-$ (spherical Gauss maps) 56, 57

\noindent $G_{J,a}$ 37

\noindent Gauss curvature 10, 25, 44

\noindent Gauss map $\nu$ 39, 44, 67, 78, 81, 

\noindent Harmonic map 101

\noindent Helices (spherical) 58

\noindent Helical cylinder 58, 59, 62

\noindent Holomorphic submanifold 36, 41

\noindent Index form 99, 100, 101

\noindent Index 99, 100, 101

\noindent Invariant subspace 106

\noindent  Jacobi field 99

\noindent  Jacobi operator 99, 101, 106-108

\noindent $J_{0}$ (a compatible complex structure) 18

\noindent $J_{1}, J_{1}^-$ (compatible complex
structures) 18, 54, 62

\noindent $J^{+}_{V}, J^{-}_{V}$ 38, 71

\noindent ${\mathfrak J},\,{\mathfrak J}^{+}, \,{\mathfrak J}^{-}$ 35,
36, 37, 39, 40, 41, 42, 55

\noindent ${\mathfrak J}_{V,a}^{+}, {\mathfrak J}_{V,a}^{-}$ 37

\noindent Kaehler form 35

\noindent Kaehlerian slant submanifold 13, 19, 22

\noindent Killing vector field 101, 105, 109

\noindent $\Lambda,\,\, [\Lambda]$ 95, 96

\noindent Lagrangian submanifold 104

\noindent $(M_{_{+}},M_{_{-}})$ 108, 109

\noindent Module 106

\noindent $\nu_+$, $\nu_-$ 39

\noindent Normal curvature 10, 25, 44

\noindent Normal subspace $\nu$ 26, 101

\noindent Normal variation 98

\noindent $\Omega_J$ (Kaehler form) 35, 77, 78

\noindent $P$ (canonical endomorphism on tangent bundle)
13, 20
 
\noindent $\phi$ (an orientation reversing isometry of
$S^3$) 53

 \noindent $\pi, \pi_{+}, \pi_{-}$ 34, 68

\noindent Parallel submanifold 17, 51

\noindent Proper slant submanifold 13

\noindent $\Psi$, $[\Psi]$ 87

\noindent Representation 106 

\noindent $S^{+}_{J,a}, S^{-}_{J,a}$ 37, 39

\noindent $S^{2}_{+}, S^{2}_{-}$ 34-39, 46 

\noindent Second variation formula 98, 99

\noindent Slant angle 13

\noindent Slant submanifold 9, 13, 46, 105

\noindent Stable submanifold 99-101, 105-109

\noindent Symmetric space 105-108

\noindent Symmetric unitary connection 8

\noindent Sympletic manifold 95, 96

\noindent $t$ (canonical 1-form on normal bundle) 14

\noindent $\Theta$, $[\Theta]$ 84-89

\noindent Totally geodesic submanifolds 14, 105-110

\noindent Totally real submanifold 8, 30, 32, 36,
82, 83, 99-101, 109, 110

\noindent highest weight 110

\noindent Whitnety immersion 42

\noindent Wirtinger angle 13

\noindent $\zeta$, $\zeta^+$, $\zeta^-$ 35, 71

\noindent $\zeta_J$ 35, 36, 37

\noindent ${\hat \zeta}_0$ 80, 82

\end{document}